\documentclass[11pt]{article}
\usepackage[margin=1in,headheight=14pt]{geometry}
\usepackage[T1]{fontenc}
\usepackage{lmodern}
\usepackage{amsmath,amssymb,amsthm,mathtools}
\usepackage{microtype}
\usepackage{titlesec}
\usepackage{needspace}
\usepackage{fancyhdr}
\usepackage[hidelinks]{hyperref}
\titleformat{\section}{\normalfont\normalsize}{\thesection.}{0.5em}{\itshape}
\titlespacing*{\section}{0pt}{2.5ex plus 1ex minus .2ex}{1.2ex plus .2ex}
\titleformat{\subsection}[runin]{\normalfont\normalsize}{\thesubsection.}{0.5em}{\itshape}[.]
\titlespacing*{\subsection}{0pt}{1.8ex plus .8ex minus .2ex}{0.6em}

\newtheorem{theorem}{Theorem}[section]
\newtheorem{proposition}[theorem]{Proposition}
\newtheorem{lemma}[theorem]{Lemma}
\newtheorem{corollary}[theorem]{Corollary}
\theoremstyle{definition}
\newtheorem{example}[theorem]{Example}

\hypersetup{pdftitle={Strong Borel--Cantelli laws and clusters of visits for random intermittent maps},pdfauthor={Qihang Wang},
pdfsubject={Strong Borel--Cantelli laws and excursion marks},
pdfkeywords={Borel--Cantelli, random LSV maps, sequential intermittency}}
\begin{document}
\thispagestyle{plain}
\begin{center}
{\Large\bfseries Strong Borel--Cantelli laws and clusters of visits\\[2pt]
for random intermittent maps\par}
\vspace{1.4em}
{\scshape Qihang Wang\par}
\vspace{0.5em}
{\itshape School of Mathematical Sciences, Peking University,
Beijing 100871, China\\
(e-mail: \href{mailto:2501110049@stu.pku.edu.cn}{2501110049@stu.pku.edu.cn})\par}
\end{center}
\vspace{1.2em}

{\small
\noindent\textit{Abstract.}
Near the neutral fixed point of an intermittent map, a single entrance
into a small target can be followed by a long run of consecutive visits,
so the number of visits need not follow its mean. We study targets
\([0,u_n]\) for iid random compositions of Liverani--Saussol--Vaienti
maps with a bounded parameter law that is regularly varying at its
minimum \(\alpha\in(0,1)\); this includes a single map. For
\(u_n=cn^{-g}L(n)\), where \(c>0\), \(0<g\leq1\) and \(L\) is
slowly varying, the visit count is almost surely asymptotic to its mean
if and only if \(\sum_n n^{-1}(nu_n)^{-(1-\alpha)/\alpha}<\infty\).
Under the stationary law, the count divided by its mean converges to one
in \(L^1\) exactly when \(nu_n\to\infty\). Thus clusters of visits
separate these two modes of convergence.
\par}
\vspace{1em}

\section{Introduction and main results}\label{sec:intro}

Let \(T\) preserve a probability measure \(\mu\), and let \((A_n)\) be
measurable sets with \(\sum_n\mu(A_n)=\infty\). The sequence \((A_n)\)
satisfies the strong Borel--Cantelli law if, for \(\mu\)-almost every
\(x\), the number of \(n\leq N\) with \(T^nx\in A_n\) is asymptotic to
\(\sum_{n\leq N}\mu(A_n)\). For independent events this follows from a
second-moment estimate. For dynamical systems with sufficiently fast
decay of correlations, variance arguments of the same kind apply to many
sequences of balls; see \cite{HNPV13} and the references there.

Intermittent maps behave differently. An orbit of a
Liverani--Saussol--Vaienti (LSV) map that comes close to the neutral
fixed point \(0\) leaves its neighbourhood only slowly, so one entrance
into a target near \(0\) may be followed by a long run of consecutive
visits. For a single LSV map with parameter \(\alpha\in(0,1)\), Kim
\cite[Proposition 4.2]{Kim07} showed that almost every orbit visits the
intervals \([0,n^{-1/(1-\alpha)})\) only finitely often, although their
invariant masses are of order \(1/n\). Gou\"ezel \cite{Gou07} proved that divergence of
the total Lebesgue length nevertheless forces infinitely many visits.
Whether the number of visits is also asymptotic to its mean is a finer
question, since a few long clusters may carry most of the count. When
the strong law fails, one would also like to describe the count itself.

In this paper we study both questions for iid random compositions of
LSV maps, with targets \([0,u_n]\) shrinking to the neutral fixed point.
The parameter law may be a point mass, so the results include a single
LSV map. Under the regular-variation assumptions stated below, the
strong law holds if and only
if an explicit integral converges, whereas convergence in mean requires
only \(nu_n\to\infty\) (Theorem~\ref{thm:neutral-counts}). At the
critical scale \(u_n\sim c\log\log n/n\), the almost-sure lower limit
of the normalized count is a constant in \((0,1)\), determined by a
scalar equation. When \(nu_n\) converges to a positive constant, the
clusters survive normalization, and the occupation paths converge to an
additive process of Poisson jumps (Theorem~\ref{thm:neutral-occupation}).
Away from the neutral fixed point, by contrast, a quenched strong law
holds for every sequence of intervals with divergent total length
(Theorem~\ref{thm:iid}). Here quenched means that the assertion holds
for almost every fixed parameter sequence.

We now describe the setting and the results precisely.
For \(s>0\), define the LSV map \cite{LSV99}

\begin{equation*}
T_s(x)=\begin{cases}
x(1+2^s x^s),&0\leq x\leq 1/2,\\
2x-1,&1/2<x\leq1.
\end{cases}
\end{equation*}

The endpoint convention is immaterial for the absolutely continuous laws
considered here. Fix \(B\in(0,\infty)\) and a Borel probability
\(\pi\) on \((0,B]\). On \(\Omega=(0,B]^{\mathbb Z}\) with law
\(\mathbb P=\pi^{\mathbb Z}\), let \(\sigma\) be the left shift
and write \(s_n(\omega)=\omega_n\). Set \(T_\omega^0=\mathrm{id}\)
and \(T_\omega^n=T_{s_{n-1}}\circ\cdots\circ T_{s_0}\).
Let \(\mathcal L_s\) be the transfer operator relative to Lebesgue
measure and \(\overline{\mathcal L}=\int\mathcal L_s\,d\pi(s)\).
Write \(\phi_s\) for the inverse of the left branch on \([0,1]\).
All return bases below use \(Y=[1/2,1]\), with nonnormalized Lebesgue
measure \(\operatorname{Leb}(Y)=1/2\).

\begin{theorem}[Sparse targets under iid forcing]\label{thm:iid}
The operator \(\overline{\mathcal L}\) has
an absolutely continuous stationary probability if and only if
\(\pi((0,1))>0\). Under this condition it is unique among absolutely
continuous stationary probabilities and equivalent to
Lebesgue measure; denote it by \(\overline\mu=h\,dx\). There is a
past-measurable probability family \(\mu_\omega=h_\omega\,dx\),
equivalent to Lebesgue measure for almost every \(\omega\), such that

\begin{equation*}
(T_{s_0})_*\mu_\omega=\mu_{\sigma\omega},\qquad
\mathbb E_{\mathbb P}h_\omega=h.
\end{equation*}

Fix \(d_0>0\) and a deterministic sequence of intervals
\(A_n\subset[d_0,1]\), \(n\geq1\), with
\(\sum_n\operatorname{Leb}(A_n)=\infty\). Put
\(p_n=\overline\mu(A_n)\) and

\begin{equation*}
M_N=\sum_{n=1}^N\overline\mu(A_n),\quad
Q_N(\omega)=\sum_{n=1}^N\mu_{\sigma^n\omega}(A_n),\quad
S_N(\omega,x)=\sum_{n=1}^N\mathbf1_{A_n}(T_\omega^n x).
\end{equation*}

For almost every \(\omega\),

\begin{equation*}
\frac{Q_N(\omega)}{M_N}\longrightarrow1,\qquad
\frac{S_N(\omega,x)}{Q_N(\omega)}\longrightarrow1
\quad\text{for Lebesgue-almost every }x.
\end{equation*}

The environment exceptional set may depend on the fixed target sequence.
No nesting, shrinking or rate of divergence is assumed.

\end{theorem}
The denominator is the expected count under the equivariant reference
family. It is not defined by evolving each arbitrary initial density.
On every fixed interval away from zero, \(h\) is bounded above and
below, so the length-divergence condition is equivalent to \(M_N\to\infty\).
Pointwise values of \(h\) and \(h_\omega\) refer to their continuous
representatives on \((0,1]\), constructed in Section~\ref{sec:fibres}.
The stationary-probability criterion concerns absolutely continuous laws;
the singular fixed-point mass at zero remains invariant in all regimes.
Sections~\ref{sec:inducing} and~\ref{sec:fibres} establish existence and the
criterion. Sections~\ref{sec:covariance} and~\ref{sec:mass}
complete the strong-law proof of Theorem~\ref{thm:iid}.

\subsection*{Targets at the neutral fixed point}

For these targets we use the following standing assumptions and notation.

\emph{Parameter law.} The law \(\pi\) is supported on \([\alpha,B]\),
where \(0<\alpha<1\), \(\alpha\leq B<\infty\), and
\begin{equation}\label{eq:RV}\tag{RV}
 \pi([\alpha,\alpha+t])\sim C_0t^{\kappa_0}L_0(1/t)\quad(t\downarrow0),
 \qquad C_0>0,\quad \kappa_0\geq0.
\end{equation}
Here \(L_0\) is positive, measurable and slowly varying at infinity,
with \(L_0\) and \(1/L_0\) locally bounded on a sufficiently large tail.
Put \(b(v)=v^{\kappa_0}/L_0(v)\), \(r=1/\alpha\),
\(\zeta=(1-\alpha)/\alpha\) and \(\overline c=h(1/2)/2\).

\emph{Visit counts.} Fix \(c>0\) and \(0<u_*<1/2\).
For a deterministic radius schedule \((u_n)\), put \(I_n=nu_n\) and
\begin{equation*}
 \begin{aligned}
 C_N(\omega,x)&=\sum_{n=1}^N\mathbf1_{[0,u_n]}(T_\omega^n x),&
 M_N&=\sum_{n=1}^N\overline\mu([0,u_n]),\\
 Q_N(\omega)&=\sum_{n=1}^N\mu_{\sigma^n\omega}([0,u_n]),&
 C_*&=\frac{h(1/2)}{2\,2^\alpha C_0\Gamma(\kappa_0+1)(1-\alpha)}.
 \end{aligned}
\end{equation*}
Here \(\Gamma\) denotes Euler's gamma function.

\emph{Regularly varying schedules.} Let \(L\) be positive and measurable,
slowly varying at infinity, with \(L\) and \(1/L\) locally bounded on a
sufficiently large tail. A regularly varying schedule is given by
\(u_n=\min(u_*,cn^{-g}L(n))\), \(0<g\leq1\), for sufficiently large
\(n\), with arbitrary positive initial radii. Put \(\xi=g(1-\alpha)\).
Define \(u(x)\) by the same formula for large real \(x\), and put
\begin{equation*}
 J(u)=\int^\infty\frac{dx}{x(xu(x))^\zeta}.
\end{equation*}
Only convergence at infinity matters.

\emph{The critical function.} Define
\begin{equation*}
 H_r(q)=\int_0^\infty\frac{1-e^{-qz}}{z}(1+z)^{-r}\,dz,\qquad
 h_r(x)=\sup_{q\geq0}\left\{H_r(q)-\frac{xq}{r-1}\right\},
 \quad 0<x<1.
\end{equation*}

The next theorem gives the count criterion, the lower envelopes and the
mean convergence for regularly varying neutral schedules.
Theorem~\ref{thm:neutral-log} makes these conclusions explicit for
logarithmic radii and near the strong-law threshold.

\begin{theorem}[Neutral counts]\label{thm:neutral-counts}
Assume \eqref{eq:RV}. For almost every \(\omega\),
\(\mu_\omega([0,u])\sim C_*u^{1-\alpha}b(\log(1/u))\) as
\(u\downarrow0\). Let \((u_n)\) be a regularly varying schedule.

\emph{(i) Normalization.} Almost surely,
\begin{equation*}
 \frac{Q_N}{M_N}\to1,\qquad
 M_N\sim\frac{C_*c^{1-\alpha}g^{\kappa_0}}{1-\xi}
 N^{1-\xi}L(N)^{1-\alpha}b(\log N).
\end{equation*}

\emph{(ii) Count law and upper fluctuations.} For almost every
environment and Lebesgue-almost every point,
\(C_N/Q_N\to1\) if and only if \(J(u)<\infty\).
For \(g<1\) this integral always converges. At \(g=1\), if
\(J(u)=\infty\) and \(I_n\to\infty\), then
\(\limsup C_N/Q_N=\infty\); if \(I_n\) does not tend to infinity,
then \(\liminf C_N/Q_N=0\).

\emph{(iii) Lower envelopes.} For almost every environment and
Lebesgue-almost every point, at \(g=1\) the lower envelope satisfies
\begin{equation*}
 \liminf_{N\to\infty}\frac{C_N}{Q_N}
 =\begin{cases}
 0,& I_n/\log\log n\to0,\\
 1,& I_n/\log\log n\to\infty.
 \end{cases}
\end{equation*}
More generally, at \(g=1\), the lower envelope is positive if
\(\liminf_nI_n/\log\log n>0\), and it is strictly less than one if
\(\limsup_nI_n/\log\log n<\infty\).
For every fixed schedule \(u_n\sim c\log\log n/n\),
\begin{equation*}
 \liminf_{N\to\infty}\frac{C_N}{Q_N}=d_\pi(c),\qquad
 \frac{c\overline c}{\alpha}h_r(d_\pi(c))=1,
\end{equation*}
where \(d_\pi(c)\) is the unique solution in \((0,1)\) of the second equation.
The function \(d_\pi\) is continuous and strictly increasing,
with limits zero as \(c\downarrow0\) and one as \(c\to\infty\).
After rescaling \(c\) by \(\overline c/\alpha\), this curve depends
only on \(\alpha\); the entrance coefficient \(\overline c\) still
depends on the parameter law.

\emph{(iv) Stationary mean convergence.} Under the stationary annealed physical law
\(\overline\Lambda=\overline\mu\times\pi^{\mathbb N_0}\), over both the
initial point and future symbols,
\(C_N/M_N\to1\) in \(L^1\) if and only if \(I_n\to\infty\).
The radii need not be monotone.

Each assertion concerns a fixed schedule;
its exceptional set may depend on the parameters of that schedule.
\end{theorem}

\begin{theorem}[Logarithmic and threshold-window schedules]\label{thm:neutral-log}
Assume \eqref{eq:RV}, and use the notation introduced before
Theorem~\ref{thm:neutral-counts}.

\emph{(i) Logarithmic schedules.} For \(u_n=\min(u_*,c(\log n)^\kappa/n)\), \(n\geq3\),
\(\kappa\in\mathbb R\), with arbitrary positive initial radii,
\begin{equation*}
 \frac{Q_N}{M_N}\to1\quad\text{almost surely},\qquad
 M_N\sim\frac{C_*c^{1-\alpha}}{\alpha}
       N^\alpha(\log N)^{\kappa(1-\alpha)}b(\log N).
\end{equation*}
Put \(\kappa_c=\alpha/(1-\alpha)\). For almost every environment
and Lebesgue-almost every point, the following classification holds:
\begin{align*}
 \kappa>\kappa_c &: \quad C_N/Q_N\to1;\\
 0<\kappa\leq\kappa_c &: \quad
       (\liminf C_N/Q_N,\limsup C_N/Q_N)=(1,\infty);\\
 -1\leq\kappa\leq0 &: \quad
       (\liminf C_N/Q_N,\limsup C_N/Q_N)=(0,\infty);\\
 \kappa<-1 &: \quad C_N\text{ is eventually constant}.
\end{align*}
For every \(\kappa>0\), \(C_N/M_N\to1\) also in \(L^1\) under the
stationary annealed physical law.

\emph{(ii) Strong-law threshold window.} Let \(\ell:[t_0,\infty)\to(0,\infty)\) be measurable and
slowly varying at infinity, with \(\ell\) and \(1/\ell\) locally bounded
on this sufficiently large tail. For sufficiently large \(n\), set
\begin{equation*}
 u_n=\min\{u_*,c(\log n)^{\kappa_c}\ell(\log n)/n\}
\end{equation*}
and choose an arbitrary positive initial prefix. Put
\begin{equation*}
 J_\ell=\int_{t_0}^{\infty}\frac{dt}{t\ell(t)^\zeta}.
\end{equation*}
Only convergence at infinity matters. Almost surely,
\begin{equation*}
 \frac{Q_N}{M_N}\to1,\qquad
 M_N\sim\frac{C_*c^{1-\alpha}}{\alpha}
 N^\alpha(\log N)^\alpha\ell(\log N)^{1-\alpha}b(\log N).
\end{equation*}
For almost every environment and Lebesgue-almost every point,
\begin{equation*}
 \begin{cases}
 C_N/Q_N\to1,&J_\ell<\infty,\\
 (\liminf C_N/Q_N,\limsup C_N/Q_N)=(1,\infty),&J_\ell=\infty.
 \end{cases}
\end{equation*}
For every such correction, \(C_N/M_N\to1\) in stationary annealed \(L^1\).
As in Theorem~\ref{thm:neutral-counts}, exceptional sets may depend on
the schedule.
\end{theorem}
For \(\ell(t)=(\log t)^\eta\), the strong law holds exactly when
\(\eta>\kappa_c\); at equality the lower and upper envelopes are
one and infinity. The schedule
\(u_n=cn^{-1}\exp(\sqrt{\log n})\) satisfies the strong law.
For \(u_n=cn^{-1}(\log\log n)^\eta\), \(\eta>0\), stationary
\(L^1\) convergence coexists with an infinite almost-sure upper limit.
The almost-sure lower limit is zero for \(0<\eta<1\),
\(d_\pi(c)\in(0,1)\) for \(\eta=1\), and one for \(\eta>1\).
Thus, even when convergence in mean holds, the lower envelope distinguishes
persistent lower bounds from subsequences with unusually few visits;
in both cases the normalized count can have arbitrarily large excesses.

Assumption~\eqref{eq:RV} requires
\(C_0L_0(v)\to\pi(\{\alpha\})\) when \(\kappa_0=0\);
it permits either a positive minimum atom or an atomless index-zero law.
For example, a CDF equal to \(C_0(\log(1/t))^{-\sigma}\) for small
\(t>0\), \(\sigma>0\), has \(\kappa_0=0\), no minimum atom,
and target mass \(C_*u^{1-\alpha}(\log\log(1/u))^\sigma\).
The uniform law on
\([\alpha,B]\), \(B>\alpha\), has \(\kappa_0=1\),
\(C_0=(B-\alpha)^{-1}\) and \(L_0=1\), even when \(B\geq1\).
The \(L^1\) assertions use the deterministic stationary mean
\(M_N\), not the random reference mean \(Q_N\). Both theorems are proved
in Section~\ref{sec:neutral}; the mean-convergence argument is in
Section~\ref{sec:neutral-l1}. Together with Theorem~\ref{thm:iid}, they
separate arbitrary sparse targets away from zero from regularly varying neutral
schedules whose clustering changes the count law.

The parameter law may be a point mass. For \(\pi=\delta_\alpha\),
Theorems~\ref{thm:neutral-counts}--\ref{thm:neutral-occupation} apply to
the single map \(T_\alpha\), with \(\kappa_0=0\), \(C_0=1\) and
\(L_0=1\), and the exceptional sets concern the initial point alone.
For this map and decreasing radii, Dedecker, Merlev\`ede and Rio
\cite{DMR} prove that \(nu_n\to\infty\) implies \(C_N/M_N\to1\) in
\(L^1\), and that \(\sum_nn^{-1}(nu_n)^{-(1-\alpha)}<\infty\) implies
the strong law. Within the regularly varying class of
Theorem~\ref{thm:neutral-counts}, that theorem shows that the first
condition is also necessary, and that the exponent \(1-\alpha\) in the
second may be replaced by \((1-\alpha)/\alpha\). Indeed, by the uniform
convergence theorem \cite[Theorem 1.2.1]{BGT}, \(J(u)<\infty\) is
equivalent to \(\sum_nn^{-1}(nu_n)^{-(1-\alpha)/\alpha}<\infty\), and
this series criterion is necessary and sufficient. For
\(u_n=c(\log n)^\kappa/n\) the earlier condition requires
\(\kappa>1/(1-\alpha)\), whereas the strong law holds exactly when
\(\kappa>\alpha/(1-\alpha)\).

The next theorem describes the corresponding occupation paths. We use the
Skorokhod space of real-valued c\`adl\`ag paths, that is, paths that are
right-continuous and have left limits.

\begin{theorem}[Neutral occupation limits]\label{thm:neutral-occupation}
Assume \eqref{eq:RV}, and take a regularly varying schedule with
\(g=1\). Put \(\overline\Lambda=\overline\mu\times\pi^{\mathbb N_0}\).
For every fixed \(T>0\), every fixed probability \(\rho\ll\overline\Lambda\),
and every deterministic subsequence \(N_k\to\infty\) such that
\(I_{N_k}=N_ku_{N_k}\to\iota\in[0,\infty]\),
\begin{equation}\label{eq:neutral-occupation-limit}
 \left(\frac{C_{\lfloor N_kt\rfloor}}{M_{N_k}}\right)_{0\leq t\leq T}
 \ \Longrightarrow\ (Z_\iota(t))_{0\leq t\leq T}
 \quad\text{in }D([0,T])\text{ with scalar }M_1\text{ topology}.
\end{equation}
The count at time zero is defined to be zero.
Here \(Z_0(t)=0\) and \(Z_\infty(t)=t^\alpha\).
For \(0<\iota<\infty\), put \(\beta_\iota=\iota\overline c/\alpha\).
Let \(\mathcal N_\iota\) be a Poisson random measure on \((0,\infty)^2\)
with intensity \(\iota\overline c\,ds/s\,r(1+w)^{-r-1}dw\), and set
\begin{equation}\label{eq:neutral-occupation-process}
 Z_\iota(t)=\frac{r-1}{\beta_\iota}
 \int_{(0,t]\times(0,\infty)}s^\alpha w\,\mathcal N_\iota(ds,dw).
\end{equation}
This nondegenerate process has independent increments,
is \(\alpha\)-self-similar, and satisfies
\begin{equation}\label{eq:neutral-occupation-transform}
 \mathbb EZ_\iota(t)=t^\alpha,\qquad
 \mathbb Ee^{-\theta Z_\iota(1)}
 =\exp\left[-\beta_\iota H_r\left(\frac{(r-1)\theta}{\beta_\iota}\right)\right],
 \quad\theta\geq0.
\end{equation}
The same limit holds with the single-horizon denominator \(Q_{N_k}\)
under the stationary two-sided law
\(m(d\omega,dx)=\mathbb P(d\omega)\mu_\omega(dx)\), and under every
fixed probability absolutely continuous with respect to \(m\).
\end{theorem}
These are annealed distributional laws, over the point and the environment.
The \(M_1\) topology allows the consecutive visits within a residence
cluster to merge into a single jump of the limiting path.
The additive process is classical: after the time change \(v=t^\alpha\),
its raw form has the scale-invariant Poisson construction of
\cite[Theorem 4.1]{PY22}, with generic jump tail
\((1+w)^{-r}\) and rate \(\beta_\iota\).
The theorem identifies this process as the limit of moving neutral counts.
It is proved in Section~\ref{sec:neutral-occupation}, using the induced
reward estimates rather than the \(L^1\) clause of Theorem~\ref{thm:neutral-counts}.

\subsection*{Irregular targets and other forcing sequences}

Regular variation is needed for the neutral count classification, but
not for a criterion for infinitely many visits.
For every fixed shrinking neutral schedule and every bounded iid law
with essential minimum in \((0,1)\), Theorem~\ref{thm:irregular} gives
such a criterion through a deterministic entrance envelope.
It measures the entrance opportunities created by future targets
during one excursion. Corollaries~\ref{cor:sampling}
and~\ref{cor:drift-classes} turn this envelope into an observation-gap
test and a comparison of recurrent schedules for two parameter laws.
In particular, summable radii can still give infinitely many visits.
These are hitting criteria, not strong laws for arbitrary neutral
counts.

Two further results clarify the distinction between entrances and
residence times without iid forcing. For arbitrary parameter sequences
in a compact subset of \((0,1)\), Theorem~\ref{thm:sequential} gives
a strong law for shrinking expanding-branch entrance intervals with
divergent total length. The reference densities are specified in
Section~\ref{sec:sequential}. Theorem~\ref{thm:neutral-sharpness}
contrasts this with a nested neutral schedule whose total length
diverges but whose normalized count has lower limit zero.
For ergodic forcing, Theorem~\ref{thm:marks} gives joint return and
residence marks and a compound-Poisson occupation limit. It assumes
parameters in a compact subunit interval and either a positive atom at
the minimum or a regularly varying minimum law with a quantitative
empirical-distribution bound. These limits concern a fixed target on
a rare-event time scale, not the moving schedule of
Theorem~\ref{thm:neutral-occupation}.

\subsection*{Related work}

For deterministic LSV maps, Frolov \cite{Fro21} obtained a strong law
for divergent interval schedules away from zero. His later paper \cite{Fro25} extends estimates
for base-supported observations to nonnegative BV observables supported
away from zero for a deterministic intermittent map.
Besides the counterexample recalled above, Kim
\cite[Proposition 4.1]{Kim07} proved the strong law for decreasing
intervals of divergent invariant mass, provided \(0\) does not lie in the
intersection of their closures. For decreasing neutral radii, Dedecker, Merlev\`ede and
Rio \cite{DMR} record the deterministic criterion for infinitely many
visits in terms of the total Lebesgue length, together with the
sufficient conditions for convergence in mean and for the strong law
compared above. For a single LSV map, Holland, Kirsebom, Kunde and Persson
\cite{HKKP24} show that almost every orbit eventually
always meets \([0,r_n]\) within its first \(n\) iterates when the
invariant mass of \([0,r_n]\) is at least \((c\log\log n/n)^{1-\alpha}\)
for a sufficiently large \(c\). This condition allows radii of order
\(\log\log n/n\). Nested targets with a fixed nonzero centre are treated
in \cite{KKP20}. These eventually-always-hitting results concern a target
frozen over the first \(n\) iterates, whereas the lower envelopes in
Theorem~\ref{thm:neutral-counts} concern cumulative visits to targets
that change at every step. For general event sequences, equivalent
characterizations of almost-sure finiteness and infinitely many
occurrences are given in \cite{BCNote}. The deterministic envelope
criterion in Section~\ref{sec:irregular} follows from the same entrance
geometry and \cite{Gou07}; the additional argument treats the
future-dependent entrance intervals produced by random compositions.

Haydn, Nicol, Persson and Vaienti \cite{HNPV13}
showed that, under an annulus regularity assumption on the invariant
measure, sufficiently fast polynomial decay of correlations gives
a strong law for balls of measure at least \(i^{-\gamma}\), \(0<\gamma<1\);
Kim's example shows that measure of order \(1/i\) does not suffice
\cite{HNPV13}. For nested intervals on Young towers, Gupta, Nicol and
Ott \cite{GNO} give lower bounds for the normalized visit count;
their application to LSV maps excludes the neutral fixed point.
Carney and Nicol \cite{CN} study the growth of fixed nonintegrable
observables, where deep excursions also play a central role.
Here the observables are indicators of moving targets, and the
normalization is their own cumulative mean.

The return geometry of random LSV maps has a well-developed theory.
Bahsoun, Bose and Duan \cite{BBD14} established polynomial mixing
for iid mixtures, with the faster map controlling the return tail.
Their results already allow a slower map with an infinite absolutely
continuous invariant measure.
Bahsoun and Bose \cite{BB16} obtained sharp inverse-endpoint and
return-tail asymptotics for two-map mixtures, limit laws for fixed
observables, and precise correlation asymptotics for observables
supported away from zero. Bahsoun, Bose and Ruziboev \cite{BBR} treated
quenched mixing and continuous parameter laws, including logarithmic
corrections to inverse endpoints. For sufficiently fast mixing
environments, Dragi\v{c}evi\'c, Hafouta and Lepp\"anen \cite{DHL} obtained quenched
memory loss and limit theorems under conditions on the essential infimum
of the random parameter. We use the inducing framework of
\cite{BB16,BBD14,BBR}.
For nonstationary compositions, Korepanov and Leppänen established
coupling estimates with uniform return-tail bounds \cite{KL21}, and
later allowed time-dependent bounds with ergodic applications \cite{KL25}.
Our sequential arguments use the uniform estimates.

There are related strong laws and distributional results for other
dynamical settings. Aimino, Nicol and Vaienti \cite{ANV} prove random
strong Borel--Cantelli laws under a spectral-gap assumption for nonnegative
observables with uniformly bounded norms in the underlying Banach space.
For sequential intermittent maps,
Nicol, T\"or\"ok and Vaienti \cite{NTV} treat uniformly \(C^1\)
observations. Their ratio law assumes parameters bounded above by a
number below \(1/2\) and a positive lower limiting mean; their result for
the full subunit range is normalized by time.
Holland and Todd \cite{HT25} obtain averaged Poisson recurrence counts
with a radius fixed over the observation horizon.
Rodriguez Sponheimer \cite[Theorems 4.2 and 4.4]{RS26} proves almost-sure
invariance principles for centred shrinking-target counts in Axiom A and
spectral-gap systems, with variance determining the Brownian time scale.
In our neutral-target problem, clusters remain
visible even after division by the expected count.

Occupation limits and deviation estimates also arise in infinite
ergodic theory. Sera \cite{Sera20,Sera26} studies occupation and waiting
times for fixed sets under an infinite invariant measure. Our counts
instead concern shrinking sets under a stationary probability, with
the target changing at each observation. For neutral holes, Bonanno
and Tikekar \cite{BT26} determine escape rates by first taking the
long-time limit for a fixed hole and then shrinking it. The lower
envelope here concerns cumulative visits to moving targets. It requires
conditional deviation estimates that are uniform over the changing
excursion rewards.

Sequential intermittent rare events are studied in \cite{FFV18} for
parameters below \(1/7\) and generic nonzero centres. Marked limits for
random intermittent systems are compared with \cite{ANT,AHV} in
Section~\ref{sec:marks}. Our prescribed right-branch entrances allow the
physical return time and the neutral residence time to be treated
jointly throughout a compact subunit parameter range.

\subsection*{Proof strategy and organization}

The proofs separate entrance frequency from residence duration by
inducing on the expanding branch. For targets away from zero, we average
the full-return operators and prove a covariance bound whose slowly
decaying term is proportional to the product of the two target masses.
The variance of the count is then bounded in terms of the cumulative
target mass rather than the elapsed time, and this gives
Theorem~\ref{thm:iid} even for arbitrarily slow divergence.
Covariance bounds of this product form, and their use in Borel--Cantelli
arguments, already occur in the deterministic renewal arguments of
\cite{Gou07,Fro21}; related fixed-observable asymptotics appear in
\cite{BB16}. Here we need uniform bounds for a bounded Borel iid parameter
law and a comparison of random and stationary cumulative target masses.

Direct second-moment control can fail for neutral targets even when
the strong law holds. Proposition~\ref{prop:blocks} shows that neither
deterministic thinning nor fixed nonnegative reweighting repairs this
failure under the stationary annealed law in a range of power schedules,
whereas block unions do satisfy a vanishing-relative-variance criterion.

For neutral targets, the proof has three further steps. First,
the proof of Lemma~\ref{lem:neutral-normalizer} describes compositions
of inverse left branches near zero, uniformly in the starting point and
in the number of inverse steps. This identifies the target masses, and
an integrable bound in the environment compares the random and
stationary means. Second, call the number of target visits during one
excursion away from the base its reward.
Lemmas~\ref{lem:neutral-rewards}--\ref{lem:neutral-long-rewards}
estimate these rewards, including conditional bounds, and give the
strong-law criterion together with its necessity. The slowly varying
correction from the parameter law enters both the target mass and the
reward tail, and cancels from the criterion. Third,
Lemma~\ref{lem:neutral-critical-rate} gives a matching conditional
estimate for the probability of unusually small counts, by analysing
ordered products of weighted return operators. Spectral perturbation
for rare events is classical \cite{KL09}; the difficulty here is to
control these products uniformly when the weights change with time.
Passing from the return-time scale to physical time then gives the
critical lower envelope in Section~\ref{sec:neutral-critical-constant}.
The same inverse-branch and weighted-operator estimates yield the
occupation limits in Section~\ref{sec:neutral-occupation}.

For irregular neutral schedules, Lemma~\ref{lem:irregular-envelope}
compares the expected future entrance envelope with one computed from
the averaged drift. The comparison is uniform over all deterministic
schedules, including laws with no atom at the minimum. A stopping-time
argument then supplies the divergent half of the hitting criterion.

Sections~\ref{sec:inducing}--\ref{sec:mass} construct the stationary
measures and prove the iid strong law. Section~\ref{sec:arbitrary}
gives the arbitrary-interval conclusion and the irregular neutral-target
criterion, followed by the comparison of thinning and block unions.
Section~\ref{sec:sequential} treats sequential entrances and the neutral
obstruction. Section~\ref{sec:neutral} proves the neutral count criterion,
the critical lower envelope and the occupation limits.
Section~\ref{sec:marks} treats joint excursion marks.

\section{Averaged return operators}\label{sec:inducing}

We use the inducing framework for random intermittent maps developed in
\cite{BB16,BBD14,BBR}. The verification below treats a bounded Borel
parameter law and gives the operator estimates needed for moving
intervals, keeping the return geometry explicit.

Throughout this section \(0<s\leq B<\infty\), and the parameters
are iid with Borel law \(\pi\). Choose \(a\in(0,1)\) such that
\(p_a=\pi((0,a])>0\), and put \(r=1/a>1\), \(\theta=r-1\).
All constants may depend on \(B,a,p_a\); none is asserted to remain
bounded as \(p_a\) tends to zero. On \(BV(Y)\) use the norm
\(\|f\|_1+\operatorname{Var}_Yf\), with complexification for
spectral arguments. Denote the averaged first-return transfers by
\(R_n\) and set \(R(z)=\sum_{n\geq1}z^nR_n\).
More explicitly, let \(\tau_\omega(y)=\inf\{n\geq1:
T_\omega^n y\in Y\}\). Extending base inputs by zero, define
\begin{equation*}
 R_n f=\mathbb E\left[\mathbf1_Y
 \mathcal L_{s_{n-1}}\cdots\mathcal L_{s_0}
 (\mathbf1_{\{\tau_\omega=n\}}f)\right].
\end{equation*}
Each \(\tau_\omega\)-partition cell is called a full first-return
atom: its return branch maps onto \(Y\). The same terminology will
be used for deterministic parameter sequences and ergodic forcing.

\begin{proposition}[Averaged renewal and stationary density]\label{prop:renewal}
Under the fast-band assumption above, \(R(1)\) has a simple fixed
probability density \(\widehat h\), which is positive and Lipschitz on
\(Y\), and the mean return \(\gamma=\int_Y\mathbb E\tau_\omega\,\widehat h\,dy\)
is finite. Put \(\overline Pf=\gamma^{-1}\widehat h\int_Yf\). The renewal
operators \(T_n=\mathbf1_Y\overline{\mathcal L}^{n}\mathbf1_Y\) on \(BV(Y)\) satisfy
\begin{equation}
\label{eq:k-renewal}
 T_n=c_n\overline P+E_n,\qquad
 |c_n-1|\leq C(1+n)^{-(r-1)},\qquad
 \|E_n\|_{BV\to BV}\leq C(1+n)^{-r}.
\end{equation}
The physical operator \(\overline{\mathcal L}\) has a unique absolutely
continuous stationary probability. Its density \(h\) satisfies
\(h\geq c>0\) and \(h(x)\leq Cx^{-a}\).
\end{proposition}
\begin{proof}
We verify the return geometry and spectrum before applying the
renewal theorems, and then construct the physical stationary density.

\subsection{Inverse distortion and an integrable fast-band tail}

For a left-inverse word \(F(u)=q_m(u)\), enumerate the inverses in
their application order:
\(q_0=u\), \(q_j=\phi_{v_j}(q_{j-1})\). Put
\(t_j=(2q_j)^{v_j}\) and \(e_j=uq_j'/q_j\). Differentiating
\(q_{j-1}=q_j(1+t_j)\) gives
\begin{equation*}
 e_j=e_{j-1}\frac{1+t_j}{1+(1+v_j)t_j},\qquad
 -\partial_u\log\phi_{v_j}'(q_{j-1}(u))
 =\frac{1+v_j}{u}\frac{1+t_j}{1+(1+v_j)t_j}(e_{j-1}-e_j).
\end{equation*}
Since \(e_0=1\) and \(e_j\) decreases to a nonnegative value,
\begin{equation}
\label{eq:k-distortion}
 0\leq-\partial_u\log F'(u)\leq\frac{1+B}{u}.
\end{equation}
A first-return inverse on \(Y\) has the form \((1+F)/2\),
contracts by at least \(1/2\), and has log-Jacobian Lipschitz constant
\(D=2(1+B)\). All positive points return for every fixed word:
below \(1/2\), \(T_s\geq T_B\), whose iterates escape every
positive starting point. Thus the full-return partitions cover the base
up to null endpoints.

Starting an inverse word at \(q_0=1/2\), a step \(v_j\leq a\)
satisfies
\begin{equation*}
 q_j^{-a}-q_{j-1}^{-a}
 =q_j^{-a}[1-(1+t_j)^{-a}]
 \geq a2^{-a-1}q_j^{-a}t_j\geq a/2.
\end{equation*}
All other increments are positive. If \(K_m\) counts the fast
symbols, then \(q_m\leq[2^a+(a/2)K_m]^{-1/a}\). Reading an iid
finite word backwards preserves its iid distribution, so \(K_m\)
is binomial with parameters \(m,p_a\). Exponential Markov gives
\begin{equation*}
 \mathbb P(K_m\leq p_am/2)
 \leq\exp[-p_a(1/2-e^{-1})m],\qquad
 \mathbb E q_m\leq C(1+m)^{-r}.
\end{equation*}
The first-return tail has Lebesgue mass \(q_{n-1}/2\). For its
inverse branch \(\psi:Y\to I\), distortion and the product rule
give \(\|\psi'(f\circ\psi)\|_{BV}\leq C|I|\|f\|_{BV}\).
Consequently
\begin{equation}
\label{eq:k-tail}
 \mathbb E\operatorname{Leb}\{\tau>n\}\leq C(1+n)^{-r},
 \qquad\sum_{n>N}\|R_n\|_{BV\to BV}\leq C(1+N)^{-r}.
\end{equation}
Parameter averages are defined by their jointly measurable branch formulas
in \(L^1\); their variation bounds follow from Fubini.
The series defining \(R(z)\) converges in operator norm on the closed
disk, and \(\sum_n n\|R_n\|<\infty\).

\subsection{The induced spectrum and aperiodicity}

After \(k\) returns, each inverse \(\psi_I\) contracts by
\(2^{-k}\), has log-Jacobian Lipschitz constant at most \(2D\),
and maps onto a partition cell \(I\) of length at most
\(\operatorname{Leb}(Y)2^{-k}\). Its total physical return length
\(r_I\) is at least \(k\). Independence at finite stopping times,
or equivalently positive concatenation of stopped finite words, gives
\begin{equation}
\label{eq:k-partition}
 R(z)^kf(u)=\mathbb E\sum_I z^{r_I}\psi_I'(u)f(\psi_I(u)).
\end{equation}
For each fixed environment the cells form a partition. We use that
partition before averaging; cells belonging to different environments
need not be disjoint. The product rule and change of variables yield
\begin{equation}
\label{eq:k-ly}
 \|R(z)^kf\|_1\leq\|f\|_1,\qquad
 \operatorname{Var}(R(z)^kf)
 \leq C2^{-k}\operatorname{Var}(f)+2D\|f\|_1,
 \quad |z|\leq1.
\end{equation}
Let \(\Pi_m\) average over a finite equal-cell partition of mesh
\(\ell_m\to0\). Its range is finite-dimensional,
\(\|f-\Pi_mf\|_1\leq C\ell_m\operatorname{Var}(f)\), and
\(\operatorname{Var}(\Pi_mf)\leq C\operatorname{Var}(f)\).
Applying \eqref{eq:k-ly} to \(f-\Pi_mf\) bounds the essential norm of
\(R(z)^k\) by \(C2^{-k}\). Hence
\(r_{\mathrm{ess}}(R(z))\leq1/2\).

We also need regularity before evaluating eigenfunctions at endpoints.
In \eqref{eq:k-partition}, replace \(f\circ\psi_I\) by its cell average
\(f_I\), obtaining \(g_{k,z}\). Distortion, partition lengths and
the sum of restricted variations give
\begin{equation*}
 \|R(z)^kf-g_{k,z}\|_\infty
       \leq C2^{-k}\operatorname{Var}(f),\qquad
 \|g_{k,z}\|_\infty+\operatorname{Lip}(g_{k,z})\leq C\|f\|_1.
\end{equation*}
Thus every eigenfunction with eigenvalue of modulus one has a Lipschitz
representative. Positive Cesàro averages of \(R(1)^k(2\mathbf1_Y)\)
give an invariant probability density \(\widehat h\), uniformly
positive and log-Lipschitz on \(Y\).

Suppose \(R(z)f=\lambda f\), \(|z|=|\lambda|=1\), and set
\(v=f/\widehat h\). Normalizing \eqref{eq:k-partition} by \(\widehat h\)
gives positive branch weights summing to one. Equality at a maximum of
\(|v|\) propagates maximum modulus along every positive branch. The
iterates of the common return-one inverse \(\psi_1(u)=(u+1)/2\) converge
to the fixed point \(1\); hence \(|v(1)|=\max|v|\). At \(1\), apply
the \(k\)-return equation for each integer \(k\). For almost every
environment all inverse images of \(1\) have maximum modulus.
Intersect these countably many sets of environments of full probability.
The return partitions
have vanishing mesh, so their inverse images of \(1\) are dense and
\(|v|\) is constant on \(Y\).

Equality in the one-return equation now gives
\(zv((u+1)/2)=\lambda v(u)\). At \(u=1\) this implies
\(\lambda=z\), and iteration implies that \(v\) is constant.
The positive return-two branch then forces \(z^2=\lambda=z\),
so \(z=\lambda=1\). All endpoint evaluations use continuous branch
extensions; the sums are uniformly absolutely convergent by distortion
and partition lengths. The \(L^1\) power bound excludes Jordan blocks.
It follows that \(R(1)\) has a spectral gap with simple projection
\(Pf=\widehat h\int_Yf\). For \(|z|<1\), the stronger bound
\(\|R(z)^kf\|_1\leq|z|^k\|f\|_1\) excludes eigenvalue one;
on the unit circle the preceding argument does so for \(z\ne1\).
Thus \(I-R(z)\) is invertible on the closed disk away from \(1\).

\subsection{Physical normalization and renewal}

Let \(U_\ell\) transfer from \(Y\) for \(\ell\) steps while
killing at the first return, and let \(U_0\) be extension by zero.
For \(\ell\geq1\), its explicit averaged formula is
\begin{equation*}
 U_\ell f=\mathbb E\left[
 \mathcal L_{s_{\ell-1}}\cdots\mathcal L_{s_0}
 (\mathbf1_{\{\tau_\omega>\ell\}}f)\right],
 \qquad U_\ell:BV(Y)\longrightarrow L^1([0,1]).
\end{equation*}
The strict survival inequality excludes paths already returned at time
\(\ell\); parameter averaging is part of the definition.
By \eqref{eq:k-tail},
\begin{equation}
\label{eq:k-density}
 \gamma=\sum_{\ell\geq0}\int U_\ell\widehat h
       =\mathbb E_{\widehat h}\tau<\infty,\qquad
 h=\gamma^{-1}\sum_{\ell\geq0}U_\ell\widehat h.
\end{equation}
Surviving pieces produce \(U_{\ell+1}\widehat h\), returning
pieces produce \(R_{\ell+1}\widehat h\), and \(R(1)\widehat h
=\widehat h\); thus \(\overline{\mathcal L}h=h\).
On the base \(h_Y=\widehat h/\gamma\); below it the time-one
term is \(\widehat h((1+x)/2)/(2\gamma)\). Hence \(h\geq c>0\)
on the whole interval. The physical rank-one operator is
\(\overline P=P/\gamma\), not \(P\).

All hypotheses of \cite[Theorem 2.2]{Gou07} now hold. The derivative
normalization is \(PR'(1)P=\gamma P\). For
\(T_n=\mathbf1_Y\overline{\mathcal L}^{n}\mathbf1_Y\) one
obtains a rank-one approximation with an operator-norm remainder of the
return-tail order. Separately, \cite[Theorem 1.1]{Gou04} bounds the
operator-norm distance to the normalized projection by the integrated-tail
order; its remainder is smaller in every regime \(r>1\). Subtracting
the first approximation from the second bounds the scalar correction,
since the normalized projection is nonzero. This proves \eqref{eq:k-renewal};
here \(T_0=\mathrm{id}\), \(c_0=1\), and
\(E_0=\mathrm{id}-\overline P\). The renewal generating function
equals \((I-R(z))^{-1}\) near zero by stopped-path decomposition.
Invertibility continues that inverse throughout the open disk; Cauchy
estimates give operator-norm convergence of its coefficient series at
every \(|z|<1\), as required by the 2004 theorem.

Split base inputs at the last return:
\begin{equation*}
 \overline{\mathcal L}^{n}(\mathbf1_Yf)
       =\sum_{\ell=0}^nU_\ell T_{n-\ell}f.
\end{equation*}
Since \(\|U_\ell\|_{BV\to L^1}\leq C(1+\ell)^{-r}\),
dominated summation and \eqref{eq:k-renewal} yield convergence to \(h\int f\)
in \(L^1\). A \(BV\) input supported in \([d,1]\), \(d>0\),
has first entry to \(Y\) within a fixed number \(K_d\) of steps,
because \(T_s\geq T_B\) on the left branch. Its finitely many
unnormalized partial-entry pieces have bounded variation by
\eqref{eq:k-distortion}, the product rule and their two boundary jumps.
Independence after entry reduces them to the base case. Density of such
inputs in \(L^1\), with contraction, proves
\begin{equation}
\label{eq:k-uniqueness}
 \overline{\mathcal L}^{n}f\longrightarrow h\int f
 \quad\text{in }L^1\quad(f\in L^1).
\end{equation}
This proves uniqueness of the absolutely continuous stationary probability.

We also need an upper bound for the stationary density near zero.
For inverse words starting
at \(x<1/2\), the fast-band argument gives
\(\mathbb E q_m(x)\leq[x^{-a}+cm]^{-1/a}+xe^{-c'm}\).
Summation bounds \(\sum_m\mathbb E q_m(x)\leq Cx^{1-a}\).
Concavity gives \(F'(x)\leq F(x)/x\); inserting this into the
killed formula in \eqref{eq:k-density} gives \(h(x)\leq Cx^{-a}\).
This is an annealed bound, not a uniform envelope for the random fibres.

\end{proof}

\section{Equivariant measures and the finiteness criterion}\label{sec:fibres}

This section constructs the equivariant reference densities
\(h_\omega\) used in Theorem~\ref{thm:iid}, and shows that an absolutely
continuous stationary probability exists exactly when
\(\pi((0,1))>0\).

\begin{proposition}[Probability fibres and finiteness]\label{prop:fibres}
If
\(\pi((0,1))>0\), the physical density has a past-measurable equivariant
probability disintegration \(h_\omega\), equivalent to Lebesgue measure,
with \(\mathbb E h_\omega=h\). If \(\pi((0,1))=0\), there is
no absolutely continuous stationary probability.

\end{proposition}
\begin{proof}
The probability fibres follow from a locally bounded density
martingale; the obstruction follows by inducing and Kac's formula.

Consider the mass-one class \(\mathcal C_B\) of decreasing
nonnegative densities for which \(x^{B+1}f(x)\) increases. Every such
density satisfies \(f(x)\leq1/x\) and is locally Lipschitz on
\([d,1]\) with constant at most \((B+1)/d^2\). For the
branch calculation set \(g=\phi_s(x)\), \(t=(2g)^s\), and
\(v=(1+x)/2\). Each
\(\mathcal L_s\), \(s\leq B\), preserves the class: in the
weighted left summand the factor
\((1+t)^{B+1}/[1+(1+s)t]\) increases, with derivative numerator
\((B-s)+B(1+s)t\geq0\); the right summand uses increasing \(x/v\).
Unweighted monotonicity and mass preservation complete the check.
Apply \eqref{eq:k-uniqueness} to the constant density. Local compactness and
closure give \(h\in\mathcal C_B\). Its stationarity extends to
continuous representatives on \((0,1]\).

Define \(h_\omega^{(m)}=\mathcal L_{s_{-1}}\cdots
\mathcal L_{s_{-m}}h\), with \(h_\omega^{(0)}=h\).
At each rational \(x>0\) these are martingales for the increasing
finite-past filtration, since the newly revealed symbol is independent
and \(\overline{\mathcal L}h=h\). They are bounded by \(1/x\).
On the intersection of the rational convergence sets, uniform local
Lipschitz bounds give local uniform convergence to \(h_\omega\).
Bounded convergence yields \(\mathbb E h_\omega(x)=h(x)\).
Fatou gives \(\int h_\omega\leq1\), whereas Tonelli gives
\(\mathbb E\int h_\omega=\int h=1\). Therefore
\(\int h_\omega=1\) almost surely, and Scheffé's lemma gives full
\(L^1\) convergence; no uniform integrability near zero is needed.
Passing to \(L^1\) limits
in \(h_{\sigma\omega}^{(m+1)}=\mathcal L_{s_0}h_\omega^{(m)}\)
gives equivariance on a shift-invariant full-measure set. Weighted monotonicity
and mass one imply \(h_\omega(1)>0\); decreasingness then implies
equivalence with Lebesgue measure. Past measurability makes the fibres
independent of future symbols, identifying the stationary annealed law
used for the pair estimates below.

For necessity suppose \(\pi((0,1))=0\). Then \(T_s\leq T_1\)
below \(1/2\). The inverse endpoints are bounded below by those of
\(T_1\), and \(q_j^{-1}-q_{j-1}^{-1}\leq2\) for that map.
Thus \(\mathbb E\operatorname{Leb}\{\tau>n\}\geq c/n\).
The preceding induced spectral proof requires only bounded parameters
and finite returns, not an integrable tail. It still gives a unique
positive induced density \(\widehat h\), including uniqueness in
\(L^1\) by approximation and contraction, with infinite mean return.
If an absolutely continuous stationary probability \(\mu\) existed,
then \(\mu(Y)>0\): every positive point eventually enters \(Y\)
for every word, and stationarity forbids all inverse entry sets having
zero mass. The stationary one-sided skew product
\(\pi^{\mathbb N}\times\mu\), induced on \(Y\), has marginal
base density invariant under \(R(1)\); this follows by independence
of the unused future at each finite return. Its normalized marginal
must be \(\widehat h\). Kac's formula gives mean return at most
\(1/\mu(Y)<\infty\), a contradiction. Singular invariant laws,
including the fixed-point atom, are not excluded.

\end{proof}

\section{Covariance estimates in physical time}\label{sec:covariance}

We now bound the correlation between visits to two targets away from
zero. The point is that the slowly decaying part of the bound carries
the product of the two target masses, not only a function of the time lag.

Use the stationary skew-product law
\(m(d\omega,dx)=\mathbb P(d\omega)\mu_\omega(dx)\), and put
\(X_n=T_\omega^n x\), \(p_n=\overline\mu(A_n)\).
Fix \(d_0>0\) and \(A_n\subset[d_0,1]\).

\begin{proposition}[Two-target covariance]\label{prop:covariance}
Let \(0<a<1\) satisfy \(\pi((0,a])>0\), and put \(r=1/a\) and
\(\theta=r-1\). For intervals \(A_n\subset[d_0,1]\) and \(i<j\),
\begin{equation}
\label{eq:k-pair}
 |\operatorname{Cov}_m(\mathbf1_{A_i}(X_i),\mathbf1_{A_j}(X_j))|
 \leq C_{d_0}p_ip_j(1+j-i)^{-\theta}
       +C_{d_0}p_j(1+j-i)^{-r}.
\end{equation}
The constant depends only on \(d_0\), the parameter law and \(a\);
it is independent of the target sequence and its cumulative mass.
\end{proposition}
\begin{proof}
We first bound killed occupation away from zero and then split
the earlier target at its first return to the base.

Use the auxiliary distance \(d=\min\{d_0,1/4\}\); constants denoted
by \(C_d\) therefore depend only on the fixed target distance \(d_0\)
and the parameter law. A surviving excursion presently in \([d,1/2)\)
returns within \(K_d\) further steps. If \(W_\ell\) is the
environment-dependent Lebesgue mass
\(W_\ell(\omega)=\operatorname{Leb}\{y\in Y:
\tau_\omega(y)>\ell,\ T_\omega^\ell y\in[d,1/2)\}\) and
\(r_j=\mathbb E\operatorname{Leb}\{\tau=j\}\), then
\begin{equation*}
 \sum_{\ell>N}\mathbb E W_\ell
 \leq K_d\sum_{j>N}r_j\leq C_d(1+N)^{-r}.
\end{equation*}
By \eqref{eq:k-distortion}, the killed inverse derivative on \([d,1/2]\)
is bounded by \(C_dW_\ell\). Consequently the operators \(U_\ell\)
from \(BV(Y)\) to \(L^\infty([d,1])\) admit norm bounds
\(a_\ell\) with \(\sum_{\ell>N}a_\ell\leq C_d(1+N)^{-r}\).
The occupation tail away from zero has exponent \(r\), not \(r-1\).
Insert \eqref{eq:k-renewal} into the last-return identity and subtract
\(h\int f=\sum_{\ell\geq0}U_\ell\overline P f\).
Splitting convolutions at \(n/2\) gives
\begin{equation}
\label{eq:k-physical}
 \|\overline{\mathcal L}^{n}(\mathbf1_Yf)-h\int f\|_{L^\infty([d,1])}
 \leq C_d|\textstyle\int f|(1+n)^{-\theta}
        +C_d\|f\|_{BV}(1+n)^{-r}.
\end{equation}
For an earlier interval \(A_i\subset[d,1]\), split
\(\mathbf1_{A_i}h\) at its first entry into \(Y\), including
time zero. There are at most \(K_d+1\) unnormalized densities
\(f_{i,k}\), with total mass \(p_i\) and total \(BV\) norm
at most \(C_d\). Indeed, each image is an interval; the pushed
function \(G'(h\circ G)\) has uniformly bounded supremum and
variation on its nonzero part, since \(h\) is locally Lipschitz
and \(G'\leq1\) has controlled logarithmic distortion. Its two
boundary jumps are bounded. Averaging preserves these bounds.

Apply \eqref{eq:k-physical} after each entry and test against \(A_j\).
Since \(h\geq c>0\), \(\operatorname{Leb}(A_j)\leq Cp_j\).
For the finitely many short lags use bounded local densities and increase
the constant. Since \(d\) depends only on \(d_0\), this proves
\eqref{eq:k-pair} under the stationary skew-product law.
\end{proof}

\section{From covariance estimates to strong laws}\label{sec:mass}

This section turns the covariance bound into the strong law of
Theorem~\ref{thm:iid}; the same two lemmas are used again for the
sequential strong law of Section~\ref{sec:sequential}.
We first isolate the elementary estimate shared by both strong-law proofs.

\begin{lemma}[Weighted lag estimate]\label{lem:weighted}
For \(\rho>0\), weights
\(a_h\in[0,1]\), and \(\sum_{h\geq1}a_h\leq M\),
\(M\geq1\),
\begin{equation}
\label{eq:mass-lag}
 \sum_{h\geq1}a_hh^{-\rho}\leq C_\rho G_\rho(M),\qquad
 G_\rho(M)=
 \begin{cases}
 M^{1-\rho},&0<\rho<1,\\
 1+\log M,&\rho=1,\\
 1,&\rho>1.
 \end{cases}
\end{equation}
\end{lemma}
\begin{proof}
Use \(a_h\leq1\) up to \(\lfloor M\rfloor\),
and the total-mass bound afterwards.
\end{proof}

In the covariance bound \eqref{eq:k-pair}, the slow term retains
both target masses and the second term is summable in the lag.
Summing \eqref{eq:k-pair} and applying the weighted lag lemma with
exponent \(r-1\) gives, on every contiguous block \(I\),

\begin{equation*}
\operatorname{Var}_m\sum_{n\in I}\mathbf1_{A_n}(X_n)
\leq C_{d_0}\left(1+\sum_{n\in I}p_n\right)^{2-\delta},
\qquad \delta=\min\{r-1,1/2\}>0.
\end{equation*}

We recall a standard variance criterion and give
the short proof needed here.

\begin{lemma}[Mass-level interpolation]\label{lem:interpolation}
Let \(Z_N\) be nonnegative and
nondecreasing, with \(e_N=\mathbb E Z_N\to\infty\),
\(0\leq e_{N+1}-e_N\leq1\), and, for some \(\delta>0\),
\(\mathbb E|Z_N-e_N|^2\leq C e_N^{2-\delta}\) when \(e_N\geq2\).
Then \(Z_N/e_N\to1\) almost surely.

\end{lemma}
\begin{proof}
Fix \(b>1\) and let \(N_k\) be the first index with
\(e_{N_k}\geq b^k\). Its overshoot is at most one. For every
\(\eta>0\), Chebyshev gives a summable upper bound
\(C_\eta b^{-k\delta}\) on
\(\mathbb P(|Z_{N_k}/e_{N_k}-1|>\eta)\). The ratio tends to one
at these indices. Monotonicity bounds the intermediate lower and upper
limits by \(1/b\) and \(b\). Intersect over a countable sequence
\(b\downarrow1\).
\end{proof}

\begin{proof}[Proof of Theorem~\ref{thm:iid}]
Equivariance
gives \(\mathbb E_m(S_N\mid\Omega)=Q_N\). Conditional
Jensen bounds \(\mathbb E_{\mathbb P}|Q_N-M_N|^2\) by
\(\operatorname{Var}_m S_N\). Apply the lemma to both \(S_N\) and
\(Q_N\), and disintegrate \(m\).
Equivalence of the reference fibres with Lebesgue measure yields the
pointwise assertion in the main theorem. This argument does not require
a separate environmental correlation estimate for \(Q_N\).
\end{proof}

\section{Arbitrary intervals and irregular neutral targets}\label{sec:arbitrary}

We first record what Theorem~\ref{thm:iid} gives for arbitrary
intervals. Section~\ref{sec:irregular} then gives a criterion for
infinitely many visits to irregular neutral targets, and
Section~\ref{sec:blocks} shows that, for a range of power schedules,
neither thinning nor fixed reweighting of the visits restores
second-moment control, whereas block unions do.

\begin{corollary}\label{cor:arbitrary}
Under the parameter-law hypotheses of Theorem~\ref{thm:iid}, with
\(\pi((0,1))>0\), every fixed
deterministic interval sequence \(A_n\subset[0,1]\) satisfying
\(\sum_n\operatorname{Leb}(A_n)=\infty\) is visited infinitely often
for almost every environment and Lebesgue-almost every initial point.

\end{corollary}
\begin{proof}
Let \(J_n=(A_{n+1}+1)/2\subset Y\). Its total length
diverges and a visit to \(J_n\) at time \(n\) forces a visit to
\(A_{n+1}\) at time \(n+1\) through the common right branch.
Theorem~\ref{thm:iid} gives infinitely many visits to \(J_n\), proving
the assertion.
\end{proof}

This corollary gives only an infinitely-often conclusion for general intervals.
Theorem~\ref{thm:iid} supplies the strong ratio for targets bounded away from zero.
Theorem~\ref{thm:neutral-counts} gives the integral criterion for the
stated regularly varying neutral schedules, and
Theorem~\ref{thm:neutral-log} gives the logarithmic and
strong-law-threshold classifications.

\subsection{A criterion for irregular neutral targets}\label{sec:irregular}
For the infinitely-often question, regular variation of the radii is
unnecessary. A target may be reached during an excursion that began much
earlier. The following envelope accounts for these possible entrances.

Suppose that \(\pi\) is supported on \([\alpha,B]\), where
\(0<\alpha<1\), \(B<\infty\), and \(\alpha\) is its essential
minimum. Define the averaged relative drift and its inverse scale by
\begin{equation}\label{eq:irregular-drift}
 a(x)=\int(2x)^s\,d\pi(s),\quad 0<x\leq\tfrac12,
 \qquad q(t)=a^{-1}\bigl((t+1)^{-1}\bigr),\quad t\geq0.
\end{equation}
The function \(a\) is continuous and strictly increasing from zero to
one, so \(q\) is well defined and \(q(0)=1/2\).

\begin{theorem}[Irregular neutral targets]\label{thm:irregular}
Under these iid parameter assumptions, fix deterministic radii
\(0\leq r_n\leq r_*<1/2\) with \(r_n\to0\), and put
\begin{equation}\label{eq:irregular-envelope}
 w_j=\sup_{h\in\mathbb N_0}\min\{r_{j+h},q(h)\},\qquad j\geq1.
\end{equation}
For almost every \(\omega\) and Lebesgue-almost every \(x\),
\begin{equation}\label{eq:irregular-criterion}
 T_\omega^n x\in[0,r_n]\ \text{infinitely often}
 \quad\Longleftrightarrow\quad \sum_{j\geq1}w_j=\infty.
\end{equation}
The environment exceptional set may depend on the fixed radius sequence.
No regular variation or atom assumption on the minimum of \(\pi\) is imposed.
\end{theorem}

For nonincreasing radii, \(w_j=r_j\). The deterministic version of this
special case is recorded in \cite{DMR}. More generally, for a single
LSV map the criterion follows from the interval theorem of \cite{Gou07}
and the inverse-branch entrance reduction used in \cite{Kim07}:
\(\phi_\alpha^h(r)\asymp\min\{r,(h+1)^{-1/\alpha}\}\), uniformly
for \(0<r\leq r_*\). In the random case the actual entrance interval
depends on future parameters, so the deterministic-target theorem cannot
be applied to it directly. We compare its expected length with
\eqref{eq:irregular-envelope}, then prove divergence by restarting
finite trials at entrances determined by the past.

Write
\(\Phi_{j,h,\omega}=\phi_{s_j(\omega)}\circ\cdots\circ
\phi_{s_{j+h-1}(\omega)}\), with \(\Phi_{j,0,\omega}=\mathrm{id}\),
and omit \(\omega\) when no confusion is possible. Let
\begin{equation}\label{eq:irregular-random-envelope}
 V_j(\omega)=\sup_{h\geq0}\Phi_{j,h,\omega}(r_{j+h}).
\end{equation}
These are measurable functions of the symbols from time \(j\) onwards.
The next estimate uses inverse-power and concentration arguments familiar
from random return-time estimates \cite{BB16,BBR,DHL}, but retains the
full drift \(a\) and is uniform over the radius sequence.

\begin{lemma}\label{lem:irregular-envelope}
There are constants \(c,C,w_0>0\), depending only on the parameter law
and \(r_*\), such that
\(cw_j\leq\mathbb E V_j\leq Cw_j\) whenever \(0<w_j<w_0\).
Moreover, one can choose a constant \(c_0\in(0,1/4)\) and deterministic
integers \(\ell_j\geq0\) so that, whenever \(0<w_j<w_0\), an
excursion started at any \(0\leq y\leq c_0w_j\) stays below
\(2c_0w_j\) through time \(\ell_j\)
with probability at least \(1/2\), and
\(2c_0w_j\leq r_{j+\ell_j}\).
\end{lemma}
\begin{proof}
First, for every fixed \(b>0\),
\begin{equation}\label{eq:irregular-drift-ratio}
 \frac{a(bx)}{a(x)}\longrightarrow b^\alpha\qquad(x\downarrow0).
\end{equation}
Indeed, under the tilted law \((2x)^s\,d\pi(s)/a(x)\), the mass of
\([\alpha+\varepsilon,B]\) is at most
\((2x)^{\varepsilon/2}/\pi([\alpha,\alpha+\varepsilon/2])\).
This tends to zero, and \(b^s\) is bounded on \([\alpha,B]\).

Choose \(\ell_j\) as the first integer for which
\(\min\{r_{j+\ell_j},q(\ell_j)\}\geq w_j/2\). Then
\begin{equation}\label{eq:irregular-witness}
 r_{j+\ell_j}\geq w_j/2,\qquad
 \ell_j+1\leq a(w_j/2)^{-1}.
\end{equation}
By \eqref{eq:irregular-drift-ratio}, choose \(c_0<1/4\) and then
\(w_0>0\) so that \(2a(2c_0w)/a(w/2)\leq1/2\) for \(w<w_0\).
Start at \(y\leq c_0w\) and stop at the first crossing \(\tau\) of
\(2c_0w\). Before crossing, the conditional mean increment is at most
\(2c_0w\,a(2c_0w)\). Summing the stopped increments, including the
possible overshoot, gives for \(h\leq a(w/2)^{-1}\)
\begin{equation}\label{eq:irregular-trial}
 \mathbb P(\tau\leq h)
 \leq\frac{\mathbb E(X_{h\wedge\tau}-y)}{c_0w}
 \leq2h\,a(2c_0w)\leq\tfrac12.
\end{equation}
Here \(X_t\) denotes the forward orbit started at \(y\) with fresh iid
symbols. Taking \(w=w_j\), \(h=\ell_j\) proves the trial assertion.
Taking \(y=c_0w_j\) also gives
\(\mathbb P(V_j\geq c_0w_j)\geq1/2\), hence the lower mean bound.

For the upper bound, put \(d_\alpha=\alpha2^{-\alpha-1}\).
Until an orbit from \(z>0\) reaches \(Y\),
\begin{equation*}
 X_t^{-\alpha}-X_{t+1}^{-\alpha}
 \geq d_\alpha X_t^{-\alpha}(2X_t)^{s_t}
 \geq d_\alpha z^{-\alpha}(2z)^{s_t}.
\end{equation*}
The last inequality uses \(s_t\geq\alpha\) and \(X_t\geq z\).
Survival below \(Y\) through time \(M\) therefore implies
\(\sum_{t<M}(2z)^{s_t}\leq d_\alpha^{-1}\).
Choose a fixed \(K>1\) with \(K^\alpha>16/d_\alpha\), and put
\(z=Kw\), \(M=\lfloor(2a(w))^{-1}\rfloor\).
For small \(w\), the iid summands \(Z_t=(2z)^{s_t}\) lie in
\([0,b_z]\), where \(b_z=(2z)^\alpha\), and their sum has mean
\(\mu=Ma(z)\geq K^\alpha/8>2/d_\alpha\).
The bounded-variable Chernoff estimate gives
\begin{equation*}
 \mathbb P\left(\sum_{t<M}Z_t\leq d_\alpha^{-1}\right)
 \leq \exp\{-\mu/(8b_z)\}
 \leq \exp\{-c_1w^{-\alpha}\}.
\end{equation*}
For example, convexity yields
\(\mathbb E e^{-tZ_t/b_z}\leq
1-(1-e^{-t})\mathbb EZ_t/b_z\); exponential Markov with
\(t=\log2\) gives the displayed bound.

If \(V_j>Kw_j\), some \(h\) has
\(\Phi_{j,h}(r_{j+h})>Kw_j\). Then \(r_{j+h}>w_j\), so
\eqref{eq:irregular-envelope} forces \(q(h)\leq w_j\) and
\(h+1\geq a(w_j)^{-1}\). The orbit from \(Kw_j\) remains on the
left through this horizon and hence through
\(M=\lfloor(2a(w_j))^{-1}\rfloor\), for small \(w_j\).
Thus every possible witness is covered by the same survival event, giving
\begin{equation*}
 \mathbb P(V_j>Kw_j)\leq e^{-c_1w_j^{-\alpha}},\qquad
 \mathbb EV_j\leq Kw_j+r_*e^{-c_1w_j^{-\alpha}}\leq Cw_j.
\end{equation*}
This completes the proof.
\end{proof}

\begin{proof}[Proof of Theorem~\ref{thm:irregular}]
Ignore the null set of points eventually mapped to zero. Every remaining
left excursion is finite, since \(T_s\geq T_B\) on the left.
An excursion landing at time \(j\) contains a scheduled hit exactly when
its landing point is at most \(V_j\). Indeed, invert the left branches
from the hit time back to \(j\). Conversely, such an inverse threshold
gives a hit before exit. A positive supremum in
\eqref{eq:irregular-random-envelope} is attained because its terms are
bounded by \(r_{j+h}\to0\). Discarding the finite initial excursion,
we obtain the pathwise equivalence
\begin{equation}\label{eq:irregular-entrance}
 \{X_n\leq r_n\ \text{infinitely often}\}
 =\{X_{j-1}\in[1/2,(1+V_j)/2]\ \text{infinitely often}\}.
\end{equation}

Suppose first that \(\sum_jw_j<\infty\). Start \(X_0\) with law
\(\overline\mu=h\,dx\), independently of the future symbols.
The variable \(X_{j-1}\) is independent of \(V_j\); the right branch
does not depend on \(s_{j-1}\). Boundedness of \(h\) on \(Y\) gives
\begin{equation*}
 \mathbb P\{X_{j-1}\in[1/2,(1+V_j)/2]\}
 =\mathbb E\int_{1/2}^{(1+V_j)/2}h(x)\,dx
 \leq C\mathbb EV_j\leq C'w_j
\end{equation*}
for all sufficiently large \(j\) with \(w_j>0\). If \(w_j=0\),
all later radii are zero. The first Borel--Cantelli lemma and
\eqref{eq:irregular-entrance} prove convergence. Equivalence of
\(\overline\mu\) with Lebesgue measure transfers this conclusion to
Lebesgue initial law.

Now suppose \(\sum_jw_j=\infty\), and use Lebesgue initial law.
Since \(w_j\leq\sup_{n\geq j}r_n\to0\), the lemma applies after a
finite prefix. The deterministic intervals
\(J_j=[1/2,(1+c_0w_j)/2]\) have divergent total length.
Theorem~\ref{thm:iid}, with a one-step index shift, gives
\(X_{j-1}\in J_j\) infinitely often almost surely.

Let \(\mathcal F_t=\sigma(X_0,s_0,\ldots,s_{t-1})\).
Choose the first sufficiently late such entrance as \(S_1\), and set
\begin{equation*}
 E_k=S_k+\ell_{S_k},\qquad
 S_{k+1}=\inf\{j\geq E_k+1:X_{j-1}\in J_j\}.
\end{equation*}
All these times are finite: each trial skips only finitely many indices
of an infinite entrance sequence. They are stopping times for
\((\mathcal F_t)\). In particular,
\(\{E_k\leq t\}=\bigcup_{j+\ell_j\leq t}\{S_k=j\}\);
the next start is determined by \(\{E_k\leq j-1\}\) and
\(X_{j-1}\), both known at time \(j-1\).
At time \(S_k\), the common right branch gives
\(X_{S_k}\leq c_0w_{S_k}\), and the unused symbols beginning at
\(S_k\) are iid and independent of \(\mathcal F_{S_k}\).
By \eqref{eq:irregular-trial}, the event \(H_k\) of staying below
\(2c_0w_{S_k}\) through \(E_k\) satisfies
\begin{equation*}
 H_k\in\mathcal F_{E_k},\qquad
 \mathbb P(H_k\mid\mathcal F_{S_k})\geq\tfrac12.
\end{equation*}
It implies a scheduled hit at \(E_k\). For deterministic \(a\leq b\),
successively conditioning the last failure, using \(E_k<S_{k+1}\), gives
\(\mathbb P(\bigcap_{k=a}^bH_k^c)\leq2^{-(b-a+1)}\).
Letting \(b\to\infty\) and then taking the countable union over
\(a\) proves that infinitely many trials succeed. Their endpoints
are strictly increasing. Finally, Fubini's theorem gives the claimed
almost-every-environment statements in both directions.
\end{proof}

The envelope has a simpler form when the target sizes decrease along a
prescribed sequence of observation times. Write \(\tau(b)=1/a(b)\).

\begin{corollary}[Prescribed observation times]\label{cor:sampling}
Under the parameter assumptions of Theorem~\ref{thm:irregular}, fix
deterministic integers \(0=n_0<n_1<n_2<\cdots\) and nonincreasing
radii \(0<b_k\leq r_*<1/2\) tending to zero. Put
\(d_k=n_k-n_{k-1}\). For almost every environment and
Lebesgue-almost every initial point,
\begin{equation}\label{eq:sampling-criterion}
 T_\omega^{n_k}x\in[0,b_k]\ \text{infinitely often}
 \quad\Longleftrightarrow\quad
 \sum_{k\geq1}b_k\min\{d_k,\tau(b_k)\}=\infty.
\end{equation}
\end{corollary}
\begin{proof}
Set \(r_{n_k}=b_k\) and \(r_n=0\) otherwise, ignoring the null set
of points eventually mapped to zero. For \(n_{k-1}<j\leq n_k\),
nonincrease of both \(b_k\) and \(q\) shows that the first future
observation attains the supremum in \eqref{eq:irregular-envelope}.
Consequently
\[
 \sum_jw_j=\sum_k F(b_k,d_k),\qquad
 F(b,d)=\sum_{h=0}^{d-1}\min\{b,q(h)\}.
\]
Uniformly over integers \(d\geq1\) and sufficiently small \(b\),
\begin{equation}\label{eq:sampling-comparison}
 F(b,d)\asymp b\min\{d,\tau(b)\}.
\end{equation}
Indeed, \(q(h)\geq b\) exactly when \(h+1\leq\tau(b)\).
Thus \(F=bd\) if \(d\leq\tau(b)\), while the first
\(\lfloor\tau(b)\rfloor\) terms give \(F\geq b\tau(b)/2\)
if \(d>\tau(b)\geq2\). For the upper bound, choose
\(\alpha<\beta<1\). The tilted-law concentration used to prove
\eqref{eq:irregular-drift-ratio} gives, for all small \(b\),
\(a(by)/a(b)\geq y^\beta/2\) uniformly in \(0<y\leq1\).
Hence, when \(h+1\geq\tau(b)\),
\[
 q(h)\leq b\left(\frac{2\tau(b)}{h+1}\right)^{1/\beta}.
\]
Since \(1/\beta>1\), summing this tail and the initial terms gives
\(F(b,d)\leq Cb\tau(b)\), independently of \(d\). Together with
\(F(b,d)\leq bd\), this proves \eqref{eq:sampling-comparison}.
Theorem~\ref{thm:irregular} now applies.
\end{proof}

The summand in \eqref{eq:sampling-criterion} grows with the preceding
gap only until that gap reaches the residence scale \(\tau(b_k)\).
For example, for an iid law on \(\{1/2,3/4\}\) with both masses
positive, \(\tau(b)\asymp b^{-1/2}\). Taking
\(n_k=k^2\) and \(b_k=\eta k^{-3/2}\), \(0<\eta<r_*\),
gives a divergent criterion comparable to \(\sum_k k^{-3/4}\),
although \(\sum_kb_k<\infty\). The positive physical-time schedule
\(r_{k^2}=b_k\), \(r_n=\eta2^{-n}\) at nonsquares therefore
also has summable radii and infinitely many hits.

We can also compare parameter laws through all their neutral-target
schedules. Fix \(0<r_*<1/2\). Let \(\mathcal R(\pi)\) be the
class of deterministic sequences \(\boldsymbol r=(r_n)\) in
\([0,r_*]\), tending to zero, for which the limsup in
\eqref{eq:irregular-criterion} has full Lebesgue measure for almost
every environment with parameter law \(\pi\). This assertion is made
separately for each schedule; its exceptional set may depend on the schedule.

\begin{corollary}[Comparison of parameter laws]\label{cor:drift-classes}
Let \(\pi_i\), \(i=1,2\), satisfy the parameter assumptions of
Theorem~\ref{thm:irregular}, with respective essential minima
\(\alpha_i\in(0,1)\), and let \(a_i\) be their averaged drifts.
Then
\begin{equation}\label{eq:drift-classification}
 \mathcal R(\pi_2)\subseteq\mathcal R(\pi_1)
 \quad\Longleftrightarrow\quad a_1(r)=O(a_2(r))\quad(r\downarrow0).
\end{equation}
In particular, \(\mathcal R(\pi_1)=\mathcal R(\pi_2)\) if and
only if \(a_1(r)\asymp a_2(r)\) near zero.
\end{corollary}
\begin{proof}
If \(a_1\leq Ca_2\) near zero, choose \(c\in(0,1)\) with
\(Cc^{\alpha_2}\leq1\). Since
\(a_2(cr)\leq c^{\alpha_2}a_2(r)\), we have
\(a_1(cr)\leq a_2(r)\) for small \(r\), and hence
\(q_1(h)\geq cq_2(h)\) for all sufficiently large integers \(h\).
Decrease \(c\) to cover the finitely many remaining indices. For every
schedule, its envelopes satisfy \(w_j^{(1)}\geq cw_j^{(2)}\).
Theorem~\ref{thm:irregular} proves the inclusion.

Conversely, write \(v_i(r)=r/a_i(r)\) and
\(\tau_i(r)=1/a_i(r)\). Choosing \(\beta_i\in(\alpha_i,1)\)
gives \(a_i(r)\geq\pi_i([\alpha_i,\beta_i])(2r)^{\beta_i}\),
so \(v_i(r)\to0\). If \(a_1\) is not \(O(a_2)\), we may
choose \(0<z_m\leq r_*\) strictly decreasing to zero such that
\[
 \frac{a_1(z_m)}{a_2(z_m)}\geq m^4,\qquad
 v_1(z_m)\leq2^{-m},\qquad m\geq2.
\]
Make a group of
\(M_m=\lceil1/(m^2v_1(z_m))\rceil\) observations of radius
\(z_m\), each preceded by a gap
\(\lceil\max_i\tau_i(z_m)\rceil\). Concatenating these finite
groups gives one deterministic schedule with nonincreasing observation
radii tending to zero. Corollary~\ref{cor:sampling} tests it under
law \(i\) by \(\sum_mM_mv_i(z_m)\). But
\[
 M_mv_1(z_m)\leq m^{-2}+2^{-m},\qquad
 M_mv_2(z_m)\geq m^2.
\]
With zero radii between observations, this schedule belongs to
\(\mathcal R(\pi_2)\) but not \(\mathcal R(\pi_1)\).
This proves necessity; applying the equivalence in both directions
proves the last assertion.
\end{proof}

The support and its minimum alone do not determine this recurrence class.
For a concrete example, take the uniform law \(\pi_1\) on
\([1/2,3/4]\) and the law \(\pi_2\) with density
\(32(s-1/2)\) on the same interval. Both are atomless. Direct
integration in \eqref{eq:irregular-drift} gives
\[
 a_1(r)\asymp\frac{r^{1/2}}{\log(1/r)},\qquad
 a_2(r)\asymp\frac{r^{1/2}}{(\log(1/r))^2}.
\]
The corresponding distribution-dependent logarithmic return scales
are already present in \cite{BBR}. Here they give opposite recurrence
outcomes for one schedule: at times \(n_k=k^3\), \(k\geq3\),
take \(b_k=ck^{-2}(\log k)^{-5}\), with \(c>0\) small enough
that \(b_k\leq r_*\). The gaps are asymptotic to \(3k^2\),
whereas \(\tau_i(b_k)\asymp k(\log k)^{5/2+i}\), \(i=1,2\).
The gaps exceed both residence scales eventually. The two series
in \eqref{eq:sampling-criterion} are respectively comparable to
\[
 \sum_{k\geq3}\frac1{k(\log k)^{3/2}}<\infty,
 \qquad
 \sum_{k\geq3}\frac1{k(\log k)^{1/2}}=\infty.
\]
Thus this same sequence of targets is hit only finitely often under
\(\pi_1\), but infinitely often under \(\pi_2\), in the quenched
almost-everywhere sense of Corollary~\ref{cor:sampling}.

\subsection{Thinning and block unions}\label{sec:blocks}

For neutral targets, thinning the observation times may fail to produce
counts with vanishing relative variance. The next proposition shows that,
in a range of power schedules, every fixed nonnegative weighting with
divergent mean fails this criterion, whereas suitable consecutive block
unions satisfy it.
Weighted Borel--Cantelli bounds were studied in \cite{FLS09}, and
abstract limitations of subsequence bounds in \cite{FL13}.
Here the obstruction comes from long residence near the neutral point.

\begin{proposition}[Weights and block unions]\label{prop:blocks}
Suppose that the iid law \(\pi\) is supported on \([\alpha,B]\),
where \(1/2<\alpha<1\), \(B<\infty\), and \(\alpha\) is its
essential minimum. Let \(\mathbf P\) denote the stationary annealed
law, with \(X_0\) distributed according to \(\overline\mu\)
independently of the future parameters, and write \(\mathbf E\) and
\(\operatorname{Var}_{\mathbf P}\) for its expectation and variance.
Fix \(1/(2\alpha)<g\leq1\), \(0<u_*<1/2\), and put
\(A_n=\{X_n\leq\min(u_*,n^{-g})\}\).
For every fixed deterministic sequence \(w_n\in[0,\infty)\), set
\(Z_N=\sum_{n\leq N}w_n\mathbf1_{A_n}\) and
\(e_N=\mathbf E Z_N\). If \(e_N\to\infty\), then
\begin{equation*}
 \frac{\operatorname{Var}_{\mathbf P}(Z_N)}{e_N^2}\longrightarrow\infty.
\end{equation*}
There also exist deterministic integers \(0=m_0<m_1<\cdots\)
with the following property. The unions
\(B_k=\bigcup_{m_{k-1}<n\leq m_k}A_n\) and their count
\(S_K=\sum_{k\leq K}\mathbf1_{B_k}\) satisfy
\begin{equation*}
 \mathbf E S_K\longrightarrow\infty,
 \qquad
 \frac{\operatorname{Var}_{\mathbf P}(S_K)}{(\mathbf E S_K)^2}
 \longrightarrow0.
\end{equation*}
\end{proposition}
\begin{proof}
Write \(H(u)=\overline\mu([0,u])\), \(\beta=1-\alpha\) and
\(\zeta=\beta/\alpha\). For every \(a\in(\alpha,1)\), the
fast-band probability is positive, so Proposition~\ref{prop:renewal}
gives \(H(u)\leq C_a u^{1-a}\). We also have
\(H(u)\geq c u^{1-\alpha}\) for small \(u\). Indeed,
stationarity, the lower bound for \(h\) on \(Y\), and
\(\phi_s\geq\phi_\alpha\) give
\begin{equation*}
 H(u)=\int H(\phi_s(u))\,d\pi(s)
       +\int_{1/2}^{(1+u)/2}h(v)\,dv
 \geq H(\phi_\alpha(u))+cu.
\end{equation*}
The inequality
\(x^{-\alpha}-T_\alpha(x)^{-\alpha}\leq\alpha2^\alpha\)
on the left branch implies
\(\phi_\alpha^j(u)\geq(u^{-\alpha}+\alpha2^\alpha j)^{-1/\alpha}\).
Iteration and summation prove the lower bound for \(H\).

There is a small \(d>0\) such that an orbit starting below
\(dN^{-1/\alpha}\) remains below \(2dN^{-1/\alpha}\) for the
next \(N\) steps, for every parameter word. To see this, use
\begin{equation*}
 x^{-\alpha}-T_s(x)^{-\alpha}\leq\alpha2^\alpha,
 \qquad s\geq\alpha,
\end{equation*}
and choose \(d\) with
\(d^{-\alpha}-\alpha2^\alpha\geq(2d)^{-\alpha}\).
Summing the inequality up to a hypothetical first crossing of
\(2dN^{-1/\alpha}\) rules out that crossing. For large \(N\),
this bound is smaller than \((2N)^{-g}\).
On the dyadic block \(I_b=\{2^b,\ldots,2^{b+1}-1\}\),
stationarity therefore gives
\begin{equation*}
 \sup_{n\in I_b}\mathbf P(A_n)\leq p_b:=C_a2^{-bg(1-a)},
 \qquad
 \mathbf P\Bigl(\bigcap_{n\in I_b}A_n\Bigr)
 \geq q_b:=c2^{-b\zeta}.
\end{equation*}
Choose \(a>\alpha\) close enough to \(\alpha\) that
\(2g(1-a)>\zeta\). Then \(\sum_b p_b^2/q_b<\infty\).
All subsequent sums may omit a fixed initial segment.

Put \(W_b(N)=\sum_{n\in I_b,\,n\leq N}w_n\).
Nonnegativity gives
\(\mathbf E Z_N^2\geq\sum_b q_bW_b(N)^2\), including a partially
completed last block. For every fixed sufficiently large \(b_0\),
Cauchy--Schwarz now yields
\begin{equation*}
 e_N\leq C_{b_0}
       +\left(\sum_{b\geq b_0}\frac{p_b^2}{q_b}\right)^{1/2}
        (\mathbf E Z_N^2)^{1/2},
\end{equation*}
where \(C_{b_0}\) is the finite expected contribution of the earlier
weights. Divide by \((\mathbf E Z_N^2)^{1/2}\geq e_N\), let
\(N\to\infty\), and then let \(b_0\to\infty\).
Thus \(e_N/(\mathbf E Z_N^2)^{1/2}\to0\), which proves the
variance assertion.

Since \(\sum_n\min(u_*,n^{-g})=\infty\),
Corollary~\ref{cor:arbitrary} gives
\(\mathbf P(A_n\text{ infinitely often})=1\).
Here we use absolute continuity of \(\overline\mu\) and Fubini's theorem
to pass to the stationary annealed law.
The necessary direction of \cite[Theorem~1.3]{BCNote} supplies the
deterministic block unions and their stated variance limit.
\end{proof}

Taking \(w_n\) to be the indicator of a deterministic set of times
shows that no subsequence with divergent expected count satisfies the
vanishing-relative-variance criterion. In fact, every fixed nonnegative
weighting with divergent mean makes the lower bound in
\cite[Theorem~1]{FLS09} asymptotically zero:
\(e_N^2/\mathbf E Z_N^2\to0\). This assertion concerns the original
hit indicators. The subevent criterion in \cite{DMR} also uses
consecutive block unions in its proof, before selecting smaller
subevents. That selection, as well as signed, path-dependent or
horizon-dependent weighting, is outside the proposition.

For example, take
\(\pi=\theta\delta_{3/4}+(1-\theta)\delta_{7/8}\),
\(0<\theta<1\), and \(g=3/4\). Theorem~\ref{thm:neutral-counts}
gives \(C_N/M_N\to1\) almost surely, while the proposition gives
\(\operatorname{Var}_{\mathbf P}(C_N)/M_N^2\to\infty\).
Rare long residences can therefore dominate second moments even when
the strong law holds. Block unions recover a vanishing-relative-variance
criterion for these recurrent events, although this conclusion supplies
no strong law for the original visit count.

\section{Shrinking entrances under sequential forcing}\label{sec:sequential}

For entrance targets one can dispense with independence and stationarity
of the maps. Throughout this section fix
\(0<\alpha\leq\beta<1\), a sequence \(s_n\in[\alpha,\beta]\),
and \(X_n=T_{s_{n-1}}\circ\cdots\circ T_{s_0}x\).
For a word starting at time \(n\), write
\(T_{n,h}=T_{s_{n+h-1}}\circ\cdots\circ T_{s_n}\) and
\(T_{n,0}=\mathrm{id}\).
For
\(A>\max\{2^\beta(\beta+2),3/(1-\beta)\}\), let
\(\mathcal C_{\beta,A}\) consist of continuous probability densities
on \((0,1]\) which are nonnegative and decreasing, for which
\(x^{\beta+1}h(x)\) is increasing, and which satisfy
\(h(x)\leq A x^{-\beta}\). Set
\(\mu_n=(X_n)_*\mu_0\), where \(\mu_0=h_0\,dx\) and
\(h_0\in\mathcal C_{\beta,A}\). The first lower bound on \(A\) is the
cone condition of \cite[(1.2)]{KL21}, following \cite{LSV99}; the second
is used for invariance below. Statements of almost-sure
convergence in this section refer to \(x\sim\mu_0\).

\begin{theorem}[Sequential entrance strong law]\label{thm:sequential}
For every sequence
\(0\leq\varepsilon_n\leq1\) such that
\(\varepsilon_n\to0\) and \(\sum_n\varepsilon_n=\infty\), put
\begin{equation*}
 B_n=[1/2,(1+\varepsilon_n)/2],\qquad
 Q_N=\sum_{n=1}^N\mu_n(B_n).
\end{equation*}
Then
\begin{equation*}
 \frac{\sum_{n=1}^N\mathbf1_{B_n}(X_n)}{Q_N}\longrightarrow1
 \quad\mu_0\text{-almost surely}.
\end{equation*}
No monotonicity, rate of shrinkage or frequency assumption on the map
sequence is required. The assertion concerns entrances through the
expanding branch, not all visits to the neutral interval.

\end{theorem}
\subsection{Regularity and target-mass mixing}

We first establish the regularity and coupling estimates used in the
proof. The strong law will follow for unions of complete return atoms;
inner and outer full-atom approximations then sandwich the actual interval.

Write \(\phi_s\) for the inverse of the left branch and set
\(g=\phi_s(x)\), \(t=(2g)^s\), \(v=(1+x)/2\). Then
\begin{equation*}
 \mathcal L_s h(x)=\frac{h(g)}{1+(1+s)t}+\frac{h(v)}2.
\end{equation*}
Both summands decrease. For weighted monotonicity the first summand becomes
\begin{equation*}
 [g^{\beta+1}h(g)]\frac{(1+t)^{\beta+1}}{1+(1+s)t}.
\end{equation*}
The logarithmic derivative of the second factor has numerator
\((\beta-s)+\beta(1+s)t\geq0\). The weighted right summand is
\(v^{\beta+1}h(v)(x/v)^{\beta+1}/2\), also increasing.
Since \(\int h=1\), decreasingness gives
\(h(v)/2\leq1/(1+x)\). Moreover, \(t\leq1\), \(x\leq2g\)
and \(s\leq\beta\) imply \(t\geq x^\beta\), while
\begin{equation*}
 1-\frac{(1+t)^\beta}{1+(1+s)t}
 \geq\frac{(1-\beta)t}{3}.
\end{equation*}
Consequently the gap between \(Ax^{-\beta}\) and the left summand
is at least \(A(1-\beta)/3>1\), which absorbs the right summand.
The cone is invariant.

Choose \(x_0>0\) with
\(Ax_0^{1-\beta}/(1-\beta)\leq1/2\). Weighted monotonicity
gives \(h(x)\leq h(1)x_0^{-(\beta+1)}\) for \(x\geq x_0\).
The mass of \([x_0,1]\) is at least \(1/2\), so \(h(1)\)
has a uniform positive lower bound. Thus every evolved density \(h_n\) satisfies
\begin{equation*}
 c\leq h_n|_Y\leq C,\qquad
 |\log h_n(x)-\log h_n(y)|\leq2(\beta+1)|x-y|
 \quad(x,y\in Y).
\end{equation*}

We use the uniform return geometry and coupling estimates of \cite{KL21}.
Fix a sufficiently large
log-Lipschitz constant \(K\) on \(Y\). Restricting a regular base
density to a full first-return atom and pushing it to \(Y\)
preserves this class after normalization. A regular base probability
evolved for any number of physical steps has either zero restriction to
\(Y\), or a restriction whose normalized density is in the same class.
This follows as well by splitting into positive components at the last
completed return. Countable positive mixtures preserve the logarithmic
Lipschitz bound; no regularity assertion is needed before a component
returns.

For completeness, the reference cone uses the following regularization
bound. On a first-entry atom \(I=(x_-,x_+]\) outside \(Y\), let
\(G:Y\to I\) be the inverse entry branch. Full-branch distortion
gives \(\sup G'\leq C|I|\) and geometry gives \(|I|\leq x_-\).
Hence
\begin{equation*}
 |\log h(G(y))-\log h(G(y'))|
 \leq (\beta+1)\frac{|G(y)-G(y')|}{x_-}
 \leq C|y-y'|.
\end{equation*}
Together with inverse-Jacobian distortion, this proves uniform regularity
of the normalized entry pushforward.

Put \(b=1/\beta>1\). The worst-word inverse endpoints are
\(O(n^{-b})\): for an inverse step \(q'=\phi_s(q)\),
\((q')^{-\beta}-q^{-\beta}\geq c>0\), uniformly on the compact
parameter range. Thus a regular base probability has return tail
\(O(n^{-b})\). Integration of the cone envelope gives the initial
entrance tail \(O(n^{-(b-1)})\). Both laws have the common latter
tail envelope required by the coupling theorem, which yields
\begin{equation*}
 \| (T_{s_{k+h-1}}\circ\cdots\circ T_{s_k})_*\nu
       -\mu_{k+h}\|_{\mathrm{TV}}
 \leq C h^{-(b-1)}
\end{equation*}
for every regular base probability \(\nu\) at time \(k\).
The densities of both restrictions to \(Y\) are uniformly bounded
and Lipschitz. For such functions the elementary interpolation inequality
\(\|f-g\|_\infty\leq C\|f-g\|_1^{1/2}\) follows by integrating
over an interval of length comparable to the maximum discrepancy;
increasing \(C\) covers large discrepancies. Therefore, for every
measurable \(J\subset Y\),
\begin{equation}
\label{eq:f-mixing}
 \left|\nu((T_{s_{k+h-1}}\circ\cdots\circ T_{s_k})^{-1}J)
            -\mu_{k+h}(J)\right|
 \leq C\operatorname{Leb}(J)h^{-\rho},
 \qquad \rho=(b-1)/2>0.
\end{equation}
This estimate is uniform in the starting time and the map sequence.

\begin{lemma}[Uniform cone and base coupling]\label{lem:coupling}
For a parameter sequence in the fixed interval \([\alpha,\beta]\subset(0,1)\)
and fixed cone and regular-base constants, uniformly in \(n,h\),
\begin{equation*}
 \|(T_{n,h})_*\eta-(T_{n,h})_*\eta'\|_{\mathrm{TV}}
       \leq C(1+h)^{-(b-1)}
\end{equation*}
whenever both probabilities are cone inputs, or one is a regular base
probability and the other a cone input. In the latter case their base
target probabilities satisfy \eqref{eq:f-mixing}.

\end{lemma}
\begin{proof}
The regular base return tail is \(O(n^{-b})\) and the
initial cone entry tail is \(O(n^{-(b-1)})\), as established above.
Applying the polynomial coupling estimate of
\cite{KL21} to the common
\(O(n^{-(b-1)})\) entry-tail envelope gives both comparisons; two
cone inputs each obey this envelope. The clause of \cite{KL21} used here
requires only a positive entry-tail exponent not exceeding the
return-tail exponent; the entry tail need not be summable. The interpolation argument above
gives the base target estimate.
\end{proof}

\subsection{Full-return covariance and cumulative-mass variance}

\begin{proof}[Proof of Theorem~\ref{thm:sequential}]
We begin with a complete-atom
schedule and then remove the alignment condition.

Let \(A_i\subset Y\) be a union of full first-return atoms at time
\(i\), with \(p_i=\mu_i(A_i)\). Conditional on \(X_i\in A_i\),
split into the normalized atom probabilities, with weights
\(\mu_i(a)/p_i\). Write \(\tau_i(a)\) for the return time and
\(h=j-i\). Atoms with \(\tau_i(a)\leq h/2\) have returned
and evolved for at least \(h/2\) further steps, so \eqref{eq:f-mixing} bounds
the discrepancy of their next-hit probability from \(p_j\) by
\(Cp_jh^{-\rho}\). For the other atoms, there is no base mass before
return. After return, positive mixtures preserve regularity of the base
densities, which bounds the next-hit probability by \(Cp_j\).
Thus, including zero-probability targets by
continuity,
\begin{equation}
\label{eq:f-pair}
 \left|\operatorname{Cov}_{\mu_0}
     (\mathbf1_{A_i}(X_i),\mathbf1_{A_j}(X_j))\right|
 \leq Cp_i p_jh^{-\rho}
       +Cp_j\mu_i(A_i\cap\{\tau_i>h/2\})
 \leq Cp_i p_jh^{-\rho}+Cp_j\min\{p_i,h^{-b}\}.
\end{equation}
No upper bound on the return times of the selected atoms is needed.

Set \(M_N=\sum_{i\leq N}p_i\). Summing the second term in the
later-time index gives
\begin{equation*}
 \sum_{j\leq N}p_j\sum_{i<j}\min\{p_i,(j-i)^{-b}\}
 \leq M_N\sum_{h\geq1}h^{-b}\leq CM_N.
\end{equation*}
For the first term apply the weighted lag lemma with \(a_h=p_{i+h}\), then
sum with weight \(p_i\). The diagonal terms are at most \(M_N\),
so \eqref{eq:f-pair} gives
\begin{equation*}
 \operatorname{Var}_{\mu_0}\sum_{i\leq N}\mathbf1_{A_i}(X_i)
 \leq CM_N[1+G_\rho(M_N)]
 \leq CM_N^{2-\delta},\qquad
 \delta=\min\{\rho,1/2\}>0,\quad M_N\geq1.
\end{equation*}
The mass-level interpolation lemma now proves the strong law for every
full-atom schedule with divergent \(M_N\). The variance criterion is also
given in \cite{Fro21}; the dynamical content here
is the full-return pair estimate and the boundary approximation below.

\subsection{Removing atom alignment}

For any finite inverse word \(F\), including the empty word, and
\(v\in[a_*,1/2]\), where \(a_*=\phi_\alpha(1/2)>0\), one has
\begin{equation}
\label{eq:f-boundary}
 F'(v)\leq C F(v)^{1+\alpha},\qquad
 \frac{F(v)}{F(w)}\leq\frac vw\quad(0<w\leq v\leq1/2).
\end{equation}
To prove the derivative bound, enumerate inverses in their application
order, \(q_0=v\), \(q_j=\phi_{u_j}(q_{j-1})\), and
\(t_j=(2q_j)^{u_j}\). Taylor's formula on the compact parameter range
gives
\begin{equation*}
 \log(1+(1+u)t)\geq(1+\alpha)\log(1+t)-Ct^2.
\end{equation*}
Since \(q_{j-1}=q_j(1+t_j)\leq2q_j\),
\begin{equation*}
 t_j^2\leq2^\beta q_j^{\alpha-1}(q_{j-1}-q_j)
 \leq2^{\beta+1-\alpha}
       \int_{q_j}^{q_{j-1}}z^{\alpha-1}\,dz.
\end{equation*}
The sum is uniformly bounded. Multiplying the inverse derivatives and
telescoping \(\sum_j\log(1+t_j)=\log(v/q_m)\) proves the first
inequality in \eqref{eq:f-boundary}. The second follows because inverse words
are increasing, concave and vanish at zero.

The boundary of \(B_n\), for small positive \(\varepsilon_n\),
cuts at most one return atom. Its image through the right branch is
\([F(a_s),F(1/2)]\), where \(a_s=\phi_s(1/2)\geq a_*\)
and \(F(a_s)\leq\varepsilon_n\leq F(1/2)\). Concavity implies
\(F(1/2)\leq C\varepsilon_n\), and the derivative bound gives
\begin{equation*}
 F(1/2)-F(a_s)\leq C\varepsilon_n^{1+\alpha}.
\end{equation*}
The base width is half this quantity. Let \(A_n^-\) be the union of
full atoms contained in \(B_n\), and \(A_n^+\) the union of full
atoms whose interiors meet \(B_n\). Uniform base-density bounds yield
\begin{equation*}
 A_n^-\subset B_n\subset A_n^+\quad\text{up to null endpoints},
 \qquad\mu_n(A_n^+\setminus A_n^-)
           \leq C\varepsilon_n^{1+\alpha}.
\end{equation*}
For \(\varepsilon_n=0\) both unions are empty. Since
\(c\varepsilon_n\leq\mu_n(B_n)\leq C\varepsilon_n\),
the relative error tends to zero along positive targets. Its cumulative
version gives
\(M_N^\pm/Q_N\to1\), where
\(M_N^\pm=\sum_{n\leq N}\mu_n(A_n^\pm)\); both sums diverge.
Apply the full-atom strong law to the two schedules and sandwich the
actual count between them. This proves the theorem.
\end{proof}

\subsection{Neutral targets: infinite visits without a strong count law}

The distinction between entrance counts and neutral visits is necessary
even for nested targets with divergent total length.

\begin{theorem}[Neutral strong-law obstruction]\label{thm:neutral-sharpness}
For each \(0<\alpha\leq\beta<1\) and the cone constant \(A\)
fixed above, there is a deterministic nonincreasing sequence
\(\varepsilon_n>0\), depending only on these constants, such that
\(\varepsilon_n\to0\) and \(\sum_n\varepsilon_n=\infty\).
For every word \(s_n\in[\alpha,\beta]\) and every
\(h_0\in\mathcal C_{\beta,A}\), put
\begin{equation*}
 C_N=\sum_{n=1}^N\mathbf1_{[0,\varepsilon_n]}(X_n),\qquad
 E_N=\sum_{n=1}^N\mu_n([0,\varepsilon_n]).
\end{equation*}
Then \(E_N\to\infty\) and, \(\mu_0\)-almost surely,
\begin{equation*}
 C_N\longrightarrow\infty,\qquad
 \liminf_{N\to\infty}\frac{C_N}{E_N}=0.
\end{equation*}
The null set may depend on the word and initial density.
\end{theorem}

\begin{proof}
Write \(H_n(u)=\mu_n([0,u])\). Cone invariance and the uniform
base bounds give constants \(c_Y,C_Y>0\) with
\begin{equation*}
 c_Y\leq h_n|_Y\leq C_Y,\qquad
 H_n(u)\leq \frac{A}{1-\beta}u^{1-\beta}.
\end{equation*}
The two inverse branches give the exact recursion
\begin{equation}
\label{eq:neutral-injection}
 H_{n+1}(u)=H_n(\phi_{s_n}(u))
             +\int_{1/2}^{(1+u)/2}h_n(v)\,dv.
\end{equation}
For \(v=\phi_s(u)\), \(s\geq\alpha\), and \(u\leq1/2\),
\begin{equation*}
 u-v=v(2v)^s\leq2^\alpha u^{1+\alpha}.
\end{equation*}
Set \(c_0=(4\,2^\alpha)^{-1}\) and
\(m(u)=\lfloor c_0u^{-\alpha}\rfloor\). Starting at \(u\),
the first \(m(u)\) coordinates of any backward inverse word remain
at least \(u/2\): their total decrease is at most
\(m(u)2^\alpha u^{1+\alpha}\leq u/4\).
Unrolling \eqref{eq:neutral-injection} for \(m(u)\) steps and
discarding the oldest nonnegative CDF therefore yields, for sufficiently
small \(u\),
\begin{equation}
\label{eq:neutral-warmup}
 H_n(u)\geq\frac{c_Ym(u)u}{4}
        \geq c_1u^{1-\alpha},\qquad
 n\geq m(u),\quad c_1=c_Yc_0/8.
\end{equation}
This finite-time bound is uniform in the word and initial cone density.

Consider a block \((a,a+L]\) with constant neutral radius \(u\),
and put \(B_u=[1/2,(1+u)/2]\). Since the left branch increases
the coordinate, a neutral visit in the block requires either
\(X_a\in[0,u]\) or \(X_j\in B_u\) for some
\(a\leq j<a+L\), up to null endpoints. A union bound gives
\begin{equation}
\label{eq:neutral-empty}
 \mu_0(\text{the block contains a neutral visit})
 \leq H_a(u)+\sum_{j=a}^{a+L-1}\mu_j(B_u)
 \leq C u^{1-\beta}+CLu.
\end{equation}

Set \(N_0=0\). After fixing the first \(k-1\) blocks, choose
\(u_k>0\) sufficiently small that
\begin{equation*}
 u_k\leq\min\{2^{-k},1/(2k),u_*\},\qquad
 L_k=\lfloor1/(ku_k)\rfloor\geq2m(u_k),\qquad
 \frac{c_1}{4k}u_k^{-\alpha}\geq k(N_{k-1}+1),
\end{equation*}
and \(u_k\leq u_{k-1}/2\) for \(k>1\). Here \(u_*\)
is a uniform threshold for \eqref{eq:neutral-warmup}.
Such choices are possible because \(\alpha<1\), so
\(L_k/m(u_k)\to\infty\) and \(u_k^{-\alpha}\to\infty\)
as \(u_k\to0\), with \(k,N_{k-1}\) fixed.
Put \(N_k=N_{k-1}+L_k\) and \(\varepsilon_n=u_k\)
for \(N_{k-1}<n\leq N_k\). The resulting schedule is
deterministic and nonincreasing, tends to zero, and satisfies
\(\sum_n\varepsilon_n=\sum_kL_ku_k\geq\sum_k1/(2k)=\infty\).

At least the last \(L_k/2\) times in block \(k\) exceed
\(m(u_k)\). Hence \eqref{eq:neutral-warmup} implies
\begin{equation*}
 E_{N_k}\geq\frac{c_1L_k}{2}u_k^{1-\alpha}
          \geq\frac{c_1}{4k}u_k^{-\alpha}
          \geq k(N_{k-1}+1).
\end{equation*}
Thus \(E_N\to\infty\). If block \(k\) is empty, then
\(C_{N_k}=C_{N_{k-1}}\leq N_{k-1}\), so
\(C_{N_k}/E_{N_k}\leq1/k\). By \eqref{eq:neutral-empty},
the events \(G_k=\{C_{N_k}/E_{N_k}\leq1/k\}\) satisfy
\(\mu_0(G_k)\geq1-Cu_k^{1-\beta}-C/k\to1\).
For each \(K\), \(\mu_0(\bigcup_{k\geq K}G_k)=1\);
continuity from above gives \(\mu_0(\limsup_kG_k)=1\).
This proves the lower-limit assertion.

Finally apply Theorem~\ref{thm:sequential} to the shifted entrance radii
\(\varepsilon_{n+1}\). There are almost surely infinitely many
\(n\) with \(X_n\in B_{\varepsilon_{n+1}}\); the common
right branch sends each such hit into \([0,\varepsilon_{n+1}]\)
at time \(n+1\). Thus \(C_N\to\infty\).
\end{proof}

The lower-limit obstruction uses only injection and cone bounds; the
entrance strong law is needed only to ensure infinitely many visits.
Long neutral residence makes the expected count large even though a block
with small total entrance mass is usually empty. Positive strong laws for
neutral schedules consequently require additional conditions on the targets.

\section{Neutral targets: strong laws and count oscillations}\label{sec:neutral}

Throughout this section, \(\pi\) satisfies the standing assumption
\eqref{eq:RV} of Section~\ref{sec:intro}.
Propositions~\ref{prop:renewal} and~\ref{prop:fibres} apply with fast
band \(a_0\in(\alpha,1)\). Throughout this section \(r=1/\alpha\) is the
exponent of the minimum parameter; it replaces the fast-band exponent
\(1/a\) of Section~\ref{sec:inducing}. Lemma~\ref{lem:neutral-normalizer}
first identifies the target mass and compares its two cumulative means.
We then estimate complete excursion rewards, including conditional bounds
given the terminal state of an induced return block. These estimates give
the strong law when the integral converges
and the required long rewards when it diverges. The critical lower
envelope uses a matching conditional deviation rate, followed by a
comparison of the induced and physical clocks. The final subsection
uses the same reward profiles to prove the occupation-process limits
directly; the stationary \(L^1\) criterion follows from those limits.
Comparison with smooth schedules removes the monotonicity assumption
without changing any of these thresholds.

Put \(r=1/\alpha\), \(l_u=\log(1/u)\),
\(a_u=u^{-\alpha}b(l_u)\), and \(w(u)=ua_u\). Extend \(b\)
to be bounded above and away from zero on finite argument intervals when needed.
The constants used in the inverse profiles are
\begin{equation*}
 A=\alpha^{\kappa_0+1}2^\alpha C_0\Gamma(\kappa_0+1),\qquad
 D_*=[\alpha2^\alpha C_0\Gamma(\kappa_0+1)]^{-1},\qquad
 f(t)=(1+t/D_*)^{-r}.
\end{equation*}
The uniform-convergence theorem, smooth-equivalent theorem and Potter
bounds for regular variation \cite[Theorems 1.2.1, 1.3.3 and 1.5.6]{BGT}
will be used below. Choose \(\widetilde L_0\sim L_0\) with
\(v\widetilde L_0'(v)/\widetilde L_0(v)\to0\), and set
\(\widetilde b(v)=v^{\kappa_0}/\widetilde L_0(v)\) and
\(\widetilde w(u)=u^{1-\alpha}\widetilde b(l_u)\). Then
\begin{equation*}
 \frac{u\widetilde w'(u)}{\widetilde w(u)}
 =1-\alpha-\frac{\kappa_0+o(1)}{l_u}>0
\end{equation*}
for small \(u\). Thus \(\widetilde w\sim w\) is continuous and
eventually increasing; no monotonicity of the measurable weight \(b\)
is assumed. Sharp inverse-endpoint asymptotics for two-map mixtures
are proved in \cite{BB16}; logarithmic corrections for continuous
parameter laws occur in \cite{BBR}. Here both the starting coordinate
and the length of the inverse word vary. Besides the resulting joint
profile, we need a common integrable environment bound to normalize
counts along shifted probability fibres.

\begin{lemma}[Endpoint normalization]\label{lem:neutral-normalizer}
For almost every \(\omega\),
\(\mu_\omega([0,u])/w(u)\to C_*\) as \(u\downarrow0\),
with \(C_*\) as in Theorem~\ref{thm:neutral-counts}. For every schedule
in that theorem, \(Q_N/M_N\to1\) almost surely and the stated
asymptotics for \(M_N\) hold.
\end{lemma}
\begin{proof}
We first express the fibre mass as a sum over backward entrances.
The main steps are to obtain an inverse profile uniformly in its
starting scale, dominate its tail by an integrable environment bound,
and then pass from the fibre asymptotic to the cumulative means.

Write \(F_\omega(u)=\mu_\omega([0,u])\), \(q_0=u\), and
\(q_j=\phi_{s_{-j}}(q_{j-1})\). Equivariance gives
\begin{equation*}
 F_\omega(u)=F_{\sigma^{-m}\omega}(q_m)+
 \sum_{j=0}^{m-1}\int_{1/2}^{(1+q_j)/2}h_{\sigma^{-j-1}\omega}(v)\,dv.
\end{equation*}
The older fibre is independent of the finite word determining \(q_m\),
so the expected residual is \(\mathbb E\overline\mu([0,q_m])\to0\).
The positive partial sums increase to \(F_\omega(u)\) almost surely,
first for rational \(u\), then for all \(u\) by continuity.

\emph{Inverse profiles.}
For an inverse step starting at \(x\), set \(z=x^{-\alpha}\).
Uniformly for \(s\in[\alpha,B]\), Taylor expansion at the actual
inverse coordinate, followed by its bounded-increment estimate, gives
\begin{equation*}
 \phi_s(x)^{-\alpha}-z
 =\alpha2^s z^{-(s-\alpha)/\alpha}+O(z^{-1}),\qquad
 0<\phi_s(x)^{-\alpha}-z\leq C.
\end{equation*}
The Laplace--Stieltjes identity for \(V=s-\alpha\) gives
\begin{equation*}
 \mathbb E e^{-tV}=t\int_0^\infty e^{-tv}\pi(V\leq v)\,dv
 \sim C_0\Gamma(\kappa_0+1)t^{-\kappa_0}L_0(t).
\end{equation*}
After substituting \(tv=y\), Potter bounds dominate the small-\(v\)
integrand by an integrable multiple of
\(e^{-y}(y^{\kappa_0-\epsilon}+y^{\kappa_0+\epsilon})\),
\(0<\epsilon<1\); the remaining part is exponentially negligible.
This also covers an atomless index-zero law. Hence the
conditional mean increment at coordinate \(z\) is
\(m(z)\sim A/b(\log z)\).

Fix integers \(J,C,d\geq1\) and let \(R\to\infty\) through integers.
For each terminal phase \(|k|\leq R^d\), read its inverse symbols
\(v_1,v_2,\ldots\) in application order and start at an integer
\(Z_0=z_0\in[1,CR]\). Set \(Z_i=q_i^{-\alpha}\), with
\(q_0=z_0^{-1/\alpha}\), and let \(\mathcal H_i\) reveal the first
\(i\) symbols. The centered increments are
\(Z_i-Z_{i-1}-m(Z_{i-1})\), whose sum is a bounded-increment
\(\mathcal H_i\)-martingale. Write \(Z_i^{k,z_0}\) for this grid trajectory.
The conditional exponential bound and union over its phase/initial grid give
\begin{equation*}
 \mathbb P\left(\max_{\substack{|k|\leq R^d,\ z_0\in\mathbb Z\cap[1,CR]\\i\leq JR}}
 \left|Z_i^{k,z_0}-z_0-\sum_{\ell=1}^i m(Z_{\ell-1}^{k,z_0})\right|
 >R^{3/4}\right)\leq C'R^{d+1}e^{-cR^{1/2}}.
\end{equation*}
Its sum over \(R\) is finite. Intersect these full-measure sets over the fixed
integers \(J,C,d\) and rational fixed bands used below. All subsequent
profile assertions use this common full-measure event.
Overlapping words need not be independent. A fixed band
\([\alpha,\alpha+\epsilon]\), chosen below one, has positive mass.
For each phase the probability that its length-\(i\) prefix contains
less than half the expected number of band symbols is at most
\(C_\epsilon e^{-c_\epsilon i}\).
Summing over \(|k|\leq R^d\) and \(i\geq K\log R\) bounds the
failure probability by \(C_\epsilon R^{d-c_\epsilon K}\).
Choose \(K=K_{\epsilon,d}\) so that this is summable in \(R\),
and include these countably many events in the common full-measure set.
The resulting prefix counts give
\(Z_i\geq c_\epsilon i^{\alpha/(\alpha+\epsilon)}\) for
\(i\geq K_{\epsilon,d}\log R\), and \(Z_i\leq z_0+Ci\).
For the compact profile, specialize to the initial grid points bracketing
\(u^{-\alpha}\) and choose an integer \(R\asymp a_u\).
For these points \(z_0\sim u^{-\alpha}\), and the preceding bounds imply
\(\log Z_i/\log z_0\to1\) uniformly on compact scaled windows,
and the martingale error is \(o(z_0)\). Bracketing a real initial coordinate
by its integer floor and ceiling changes it by at most one, hence gives
asymptotically equivalent radii and scales; monotonicity of the inverse
maps extends the grid limits to every real \(u\downarrow0\):
\begin{equation}\label{eq:neutral-profile}
 q_{\lfloor ta_u\rfloor}(u)/u\longrightarrow f(t)
 \quad\text{uniformly on compact }t\text{-intervals}.
\end{equation}
Here the terminal phase can vary with \(u\), with absolute value at most
\(a_u^d\), for each fixed \(d\); the phase is suppressed in the notation.
The scaled drift is \(A\alpha^{-\kappa_0}=1/D_*\).
For the deep profile take \(R=L\). The same argument, discarding
\(i<R^{1-\delta}\), gives
\begin{equation}\label{eq:neutral-deep-profile}
 q_L(u)\sim A^{-r}L^{-r}b(\log L)^r
 \quad\text{if }u^{-\alpha}=o(L/b(\log L)).
\end{equation}
Indeed the discarded cost \(O(R^{1-\delta})\) and martingale error
\(O(R^{3/4})\) are both \(o(R/b(\log R))\). On the surviving
steps the fixed-band bound gives
\(\log Z_i/\log R\in[\ell+o(1),1+o(1)]\), where
\(\ell=(1-\delta)\alpha/(\alpha+\epsilon)\). Uniform convergence
of regular variation therefore gives
\begin{equation*}
 A\leq\liminf\frac{Z_R b(\log R)}R
 \leq\limsup\frac{Z_R b(\log R)}R\leq A\ell^{-\kappa_0}.
\end{equation*}
The initial coordinate is negligible by the displayed condition.
Let \(\epsilon,\delta\downarrow0\). The same estimates apply to reversed
future words. In particular, on the same full-measure set the deep limit holds
along deterministic pairs \((u,L)\) satisfying its displayed condition
and terminal phases of absolute value at most \(L^d\), for any fixed \(d\).

\emph{An integrable environment bound.}
To control all residence tails, take \(u_m=e^{-m}\),
\(\beta_m=\alpha+1/m<a_0<1\),
\(p_m=\pi([\alpha,\beta_m])\sim C_0/b(m)\), and
\(t_m=\lceil e^{\alpha m}b(m)\rceil\).
Let \(G_m\) mean that every past prefix \(j\geq t_m\) has at
least \(p_mj/2\) such fast symbols. Bernoulli exponential bounds give
\begin{equation*}
 \mathbb P(G_m^c)\leq Cp_m^{-1}e^{-cp_mt_m}
 \leq Cb(m)\exp(-ce^{\alpha m}).
\end{equation*}
On \(G_m\), the inverse \(\beta_m\)-power increment yields
\(q_j(u_m)\leq Cu_m(j/t_m)^{-1/\beta_m}\) for \(j\geq t_m\).
The constant is uniform because \(e^{m(1-\alpha/\beta_m)}\) is bounded.
Thus \(\sum_jq_j(u_m)\leq Cu_mt_m\); its normalized tail beyond
\(Tt_m\) is at most \(CT^{1-1/a_0}\).

On the exceptional scales use the fixed band \([\alpha,a_0]\).
Let \(\xi_\ell\) be its successive geometric waiting times in the past,
\(A_\ell=\sum_{i=1}^\ell\xi_i\), and
\(\mathcal K_\omega=\sup_{\ell\geq1}A_\ell/\ell\).
For some \(\lambda>0\), \(\mathbb E e^{\lambda\xi_1}<\infty\), so for large \(t\),
\begin{equation*}
 \mathbb P(\mathcal K_\omega>t)
 \leq\sum_{\ell\geq1}(e^{-\lambda t}\mathbb E e^{\lambda\xi_1})^\ell
 \leq Ce^{-ct}.
\end{equation*}
Each fast step increases the inverse \(a_0\)-power by at least a fixed
constant. Abel summation, using \(A_\ell\leq\mathcal K_\omega\ell\), gives
simultaneously for all small \(u>0\),
\begin{equation*}
 \sum_jq_j(u)\leq\sum_{\ell\geq0}\xi_{\ell+1}(u^{-a_0}+c\ell)^{-1/a_0}
 \leq\mathcal K_\omega\sum_{\ell\geq0}(u^{-a_0}+c\ell)^{-1/a_0}
 \leq C\mathcal K_\omega u^{1-a_0}.
\end{equation*}
Consequently
\begin{equation}\label{eq:neutral-envelope}
 \frac{\sum_jq_j(u)}{w(u)}\leq\mathcal J_\omega,
 \qquad\mathbb E\mathcal J_\omega<\infty.
\end{equation}
One may take a constant multiple of
\(1+\mathcal K_\omega\sum_{m\geq m_0}
e^{(a_0-\alpha)m}b(m)^{-1}\mathbf1_{G_m^c}\).
Cauchy--Schwarz proves integrability without independence: the expectation
of its \(m\)th exceptional summand is at most
\(e^{(a_0-\alpha)m}b(m)^{-1}(\mathbb E\mathcal K_\omega^2)^{1/2}
\mathbb P(G_m^c)^{1/2}\), a summable sequence.
Adjacent exponential-radius grids have bounded weight ratios, so the
bound extends to all small real \(u\). The first Borel--Cantelli lemma makes all sufficiently
large scales good and gives the stated uniform tail control.
Since \(h_\omega\leq2\) on \(Y\), the injection identity also gives
\(F_\omega(u)\leq\sum_jq_j(u)\).

\emph{Fibre asymptotics.}
Local Lipschitz bounds give each injection as
\(c(\sigma^{-j-1}\omega)q_j+O(q_j^2)\), where
\(c(\omega)=h_\omega(1/2)/2\). The total error is at most
\(Cu\sum_jq_j\). Approximate the profile using
\eqref{eq:neutral-profile} and its uniform tail, then apply weighted
backward Birkhoff to the bounded coefficients \(c\). No independence
of trace coefficients and inverse words is used. The limit is
\(\mathbb Ec\int_0^\infty f(t)\,dt
=\mathbb EcD_* /(r-1)=C_*\).

\emph{Cumulative normalization.}

Set \(R_\delta(\omega)=\sup_{0<u\leq\delta}
|F_\omega(u)/\widetilde w(u)-C_*|\), taking the supremum over rationals.
These errors decrease to zero and have an integrable common envelope by
\eqref{eq:neutral-envelope}. First suppose that the radii are eventually
decreasing. Then \(\widetilde w(u_n)\) is eventually decreasing and regularly varying
with exponent \(-\xi\), \(\xi<1\). Birkhoff and Abel summation give
\(\sum_{n\leq N}\widetilde w(u_n)R_\delta(\sigma^n\omega)/\sum_{n\leq N}\widetilde w(u_n)
\to\mathbb ER_\delta\). Indeed, for \(a_n=\widetilde w(u_n)\)
and centered Birkhoff sums \(S_n=o(n)\),
\begin{equation*}
 \sum_{n\leq N}a_n(S_n-S_{n-1})
 =a_NS_N+\sum_{n<N}(a_n-a_{n+1})S_n
 =o\left(\sum_{n\leq N}a_n\right).
\end{equation*}
Discard the finite prefix \(u_n>\delta\), then let \(\delta\downarrow0\).
This proves \(Q_N\sim C_*\sum_{n\leq N}w(u_n)\).
Dominated convergence gives the same asymptotic for \(M_N\).

For general target \(L\), the smooth-equivalent theorem gives
\(\widetilde L\sim L\) with \(x\widetilde L'(x)/\widetilde L(x)\to0\).
Thus \(\widetilde u_n=cn^{-g}\widetilde L(n)\) is eventually decreasing.
For every rational \(\epsilon\in(0,1)\), the decreasing comparison radii
\(v_n^\pm=(1\pm\epsilon)\widetilde u_n\), with the same cutoff,
bracket \(u_n\) beyond a finite deterministic prefix. Let \(M_N^\pm,Q_N^\pm\)
be their stationary and reference means. Regular variation of the target
mass, and divergence of \(M_N\), give
\begin{equation}\label{eq:neutral-sandwich-means}
 M_N^\pm/M_N\longrightarrow(1\pm\epsilon)^{1-\alpha}.
\end{equation}
The decreasing-weight argument just proved gives \(Q_N^\pm/M_N^\pm\to1\).
Squeeze \(Q_N\) between these means with a bounded prefix error, intersect
the countably many full-measure sets, and let \(\epsilon\downarrow0\). This proves
\(Q_N/M_N\to1\) without monotonicity of the actual radii. Finally,
uniform convergence compares \(w(u_n)\) with the integral of \(w(u(t))\)
over \([n,n+1]\). Since \(w(u(t))\) is regularly varying with exponent
\(-g(1-\alpha)>-1\), Karamata's integral asymptotic
\cite[Proposition 1.5.8]{BGT} gives the displayed sum asymptotics, using
\(b(g\log N)/b(\log N)\to g^{\kappa_0}\).
\end{proof}

\subsection*{Excursion rewards and conditional deviations}

Let \(\mathcal F\) be the first-return skew map on \(Y\) with
unused one-sided iid noise. Its invariant law is
\(\nu=\widehat h(y)\,dy\times\pi^{\mathbb N_0}\); independence
of the unused suffix holds at every finite stopping time. Denote its
transfer relative to \(\nu\) by \(\mathcal P\). For a base function \(g\),
\begin{equation*}
 \mathcal Pg=R(1)(\widehat h g)/\widehat h.
\end{equation*}
Positivity and Lipschitz regularity of \(\widehat h\) transfer the
BV spectral gap to this normalized operator. Put
\(\mathsf T_0=0\), \(\mathsf T_j=\sum_{\ell=0}^{j-1}\tau\circ\mathcal F^\ell\).
Let \(y_j\) be the base coordinate of \(\mathcal F^j(y,\omega)\).
Excursion \(j\) contributes neutral rewards at physical times
\(\mathsf T_j+1,\ldots,\mathsf T_{j+1}-1\). Define
\(K_N=\max\{j:\mathsf T_j\leq N\}\); excursion \(K_N\) is the possibly
unfinished one at time \(N\). For a point initially outside \(Y\),
add its finite first-entry time to these clocks and treat the preceding
observations separately. We use the rewards
\begin{equation*}
 V_u=\sum_{\ell=1}^{\tau-1}\mathbf1_{[0,u]}(T_\omega^\ell y),\qquad
 m_u=\nu(V_u)=\gamma\overline\mu([0,u]),\quad r=1/\alpha>1.
\end{equation*}
The occupation identity \eqref{eq:k-density} proves the equality for
\(m_u\); neither excursion endpoint belongs to the neutral target.

\begin{lemma}[Conditional reward bounds]\label{lem:neutral-rewards}
There are bounded functions \(\chi_u\), with
\(\|\chi_u\|_\infty\leq Cm_u\), such that
\(d_u=V_u-m_u+\chi_u-\chi_u\circ\mathcal F\) satisfies
\(\mathcal P d_u=0\), \(d_u\geq-Cm_u\), and
\begin{equation}\label{eq:neutral-conditional}
 \begin{aligned}
 \|\mathcal P\mathbf1_{V_u>t}\|_\infty&\leq C\min(u,t^{-r}b(\log(\mathrm e+t))^r),\\
 \|\mathcal P\mathbf1_{|d_u|>t}\|_\infty&\leq C\min(u,t^{-r}b(\log(\mathrm e+t))^r)
 \quad(t\geq C_1m_u).
 \end{aligned}
\end{equation}
Here \(C_1\) is a fixed constant. The first bound holds for every \(t>0\). Also
for some \(\vartheta\in(0,1)\) and all \(k\geq1\),
\begin{equation}\label{eq:neutral-reward-gap}
 \|\mathcal P^k(V_u-m_u)\|_\infty\leq C\vartheta^{k-1}m_u.
\end{equation}
In particular, \(\|\mathcal P^kV_u\|_\infty\leq Cm_u\).
Moreover,
\(\mathcal P\mathbf1_{B_u}\leq Cu\) for
\(B_u=[1/2,(1+u)/2]\), and \(\mathcal F\) is ergodic.
\end{lemma}
\begin{proof}
The entrance set is \(B_u=[1/2,(1+u)/2]\), and
\(\{V_u>0\}=\{y\in B_u\}\), up to endpoints. Consider a
return atom with nonzero reward. Its right-image endpoints are
\(F_0(\phi_s(1/2))\) and \(F_0(1/2)\), where \(F_0\) is a
left-inverse prefix. Concavity, \(F_0(0)=0\), and
\(\phi_s(1/2)\geq\phi_\alpha(1/2)>0\) show that their ratio
is at most \(C_\alpha=1/(2\phi_\alpha(1/2))\). The entire atom
thus lies in \(B_{C_\alpha u}\), even if its boundary cuts \(B_u\).

On a return-\(n\) atom \(I\), the reward decreases with the entrance
coordinate, so its variation is at most \(n\). Distortion gives
\(\|\psi_I'\|_\infty+\operatorname{Var}\psi_I'\leq C|I|\).
For a fixed Lipschitz input \(f\), the branch product therefore has
BV norm at most \(C\|f\|_{\mathrm{Lip}}n|I|\).
From entrance coordinate \(v\), \(Cv^{-\alpha}\) symbols in
\([\alpha,\alpha+1/\log(1/v)]\) suffice to return: the extra
factor \(v^{-1/\log(1/v)}\) is bounded. Their probability is
comparable to \(1/b(\log(1/v))\), so a negative-binomial
waiting bound gives \(\mathbb E\tau(v)\leq Cv^{-\alpha}b(\log(1/v))\).
Partitioning before averaging proves for the weighted transfer
\(W_uf=\mathbb E\sum_I(V_u\circ\psi_I)(f\circ\psi_I)\psi_I'\)
that
\begin{equation*}
 \|W_uf\|_{BV}\leq C\|f\|_{\mathrm{Lip}}
 \int_0^{C_\alpha u}\mathbb E\tau(v)\,dv
 \leq C\|f\|_{\mathrm{Lip}}w(u).
\end{equation*}
For \(u\) in the remaining compact range, constants absorb the bound.
The same argument for the indicator of \(B_u\), whose branch variation
is at most one, gives \(\|R(1)(\mathbf1_{B_u}\widehat h)\|_{BV}\leq Cu\).
Stopping removes the consumed prefix, so both transfers are base-only.
Consequently \(\mathcal P\mathbf1_{V_u>0}\leq Cu\).
For the return tail at \(n\), use instead the band
\([\alpha,\alpha+1/\log n]\). Its first \(n\) symbols contain
\(cn/b(\log n)\) successes except with probability
\(Ce^{-cn/b(\log n)}\). For \(\beta_n=\alpha+1/\log n\),
the inverse-power bound on the good event is
\(q_{n-1}(1/2)\leq C[n/b(\log n)]^{-1/\beta_n}\).
This is at most \(Cn^{-r}b(\log n)^r\), since
\((r-1/\beta_n)\log[n/b(\log n)]\) stays bounded. Thus
\(\mathbb E q_{n-1}(1/2)\leq Cn^{-r}b(\log n)^r\).
On the exceptional event use \(q_{n-1}\leq1/2\).
Since \(V_u\leq\tau\), distortion and the positive stopping-transfer
formula give the first bound in \eqref{eq:neutral-conditional}.

By Lemma~\ref{lem:neutral-normalizer}, \(m_u\asymp w(u)=ua_u\).
The base-only function \(\mathcal P(V_u-m_u)=W_u\widehat h/\widehat h-m_u\)
has mean zero and BV norm at most \(Cm_u\). The induced spectral gap
gives, for some \(\vartheta\in(0,1)\),
\(\|\mathcal P^k(V_u-m_u)\|_\infty\leq C\vartheta^{k-1}m_u\).
Set \(\chi_u=\sum_{k\geq1}\mathcal P^k(V_u-m_u)\).
The transfer identity \(\mathcal P(f\circ\mathcal F)=f\) proves
centering. The bounds on \(\chi_u\) give the second tail bound and
\(d_u\geq-Cm_u\), as well as the asserted bound on \(\mathcal P^kV_u\).

For ergodicity, take a BV base function times a bounded indicator of a
cylinder in \(m\) future symbols. After \(m\) returns all these symbols
are consumed; its transfer is base-only and has finite BV norm by the
full-branch product bound. Subsequent transfers converge to its mean by
the gap. These test functions span a dense subspace of \(L^1(\nu)\).
Contraction extends convergence to \(L^1\), proving ergodicity. In
particular, physical return times satisfy \(\mathsf T_j/j\to\gamma\).
\end{proof}

\begin{lemma}[Two-sided and lower deviations]\label{lem:neutral-deviations}
Fix \(\varepsilon,\Lambda,H>0\) and put \(I=nu\), \(\zeta=r-1\).
For small \(u\), sufficiently large \(I\), and
\(\log I\leq H\log(1/u)\),
\begin{align}
 \nu\left(\left|\sum_{j<n}V_u\circ\mathcal F^j-nm_u\right|>
 \varepsilon nm_u\right)&\leq C_{\varepsilon,\Lambda,H}I^{-\zeta}+2I^{-\Lambda},
 \label{eq:neutral-upper}\\
 \nu\left(\max_{j<n}V_u\circ\mathcal F^j>\varepsilon nm_u\right)
 &\leq C_{\varepsilon,H}I^{-\zeta}.\label{eq:neutral-last}
\end{align}
For large \(n\) and small \(u\), without the logarithmic restriction,
the endpoint conditional bound is
\begin{equation}\label{eq:neutral-endpoint-lower}
 \left\|\mathcal P^n\mathbf1_{\{\sum_{j<n}V_u\circ\mathcal F^j
 <(1-\varepsilon)nm_u\}}\right\|_\infty
 \leq\exp(-c_\varepsilon nu).
\end{equation}
Consequently, for every deterministic start \(s\geq0\),
\begin{equation}\label{eq:neutral-lower}
 \nu\left(\sum_{j=s}^{s+n-1}V_u\circ\mathcal F^j<(1-\varepsilon)nm_u\right)
 \leq\exp(-c_\varepsilon nu).
\end{equation}
\end{lemma}
\begin{proof}
Let \(x=\varepsilon nm_u/2\) and \(t_c=x/q\), where \(q\) is a
large constant chosen below, independently of \(I\). Put
\(D=d_u\mathbf1_{|d_u|\leq t_c}\), and
\(e=D-(\mathcal PD)\circ\mathcal F\). Then \(\mathcal Pe=0\),
\(|e|\leq2t_c\), and \(\mathcal Pe^2\leq\mathcal PD^2\).
For large \(I\), \(t_c/a_u\asymp I/q\to\infty\),
so \(t_c\geq a_u\) and \(t_c\gg m_u\).
Fix \(\epsilon_0>0\) and put \(\rho=(\kappa_0+\epsilon_0)r\).
Potter bounds and \(\log a_u=\alpha l_u+o(l_u)\) give, for \(s\geq1\),
the conditional tail at \(t=a_us\) bounded by
\(Cu s^{-r}(1+\log s/l_u)^\rho\).
Write \(s_0=t_c/a_u\). The restriction \(\log I\leq Hl_u\)
bounds the logarithmic factor up to \(s_0\). Above it, put
\(s=s_0v\); the factor is at most \(C_H(1+\log v)^\rho\).
Since \(\int_1^\infty v^{-r}(1+\log v)^\rho\,dv<\infty\),
conditional integration, including \(t_c\mathcal P\mathbf1_{|d_u|>t_c}\), gives
\(\|\mathcal PD\|_\infty\leq C_Hua_u(I/q)^{1-r}\).
Its total bias divided by \(x\) is \(O_{q,H}(I^{-\zeta})\to0\).
Splitting the second-moment integral at \(a_u\) gives
\begin{equation*}
 \|\mathcal Pe^2\|_\infty\leq\sigma^2,\qquad
 \sigma^2\leq C_Hua_u^2\begin{cases}
 1,&r>2,\\
 1+\log I,&r=2,\\
 I^{2-r},&1<r<2.
 \end{cases}
\end{equation*}
The range below \(C_1m_u\) contributes only \(Cm_u^2\).
Let \(\mathcal G\) be the underlying sigma-algebra and
\(\mathcal G_j=\mathcal F^{-j}\mathcal G\) the decreasing filtration.
The transfer identity gives
\begin{equation*}
 \mathbb E_\nu[e\circ\mathcal F^j\mid\mathcal G_{j+1}]=0,\qquad
 \mathbb E_\nu[e^2\circ\mathcal F^j\mid\mathcal G_{j+1}]
 =(\mathcal Pe^2)\circ\mathcal F^{j+1}\leq\sigma^2.
\end{equation*}
Set \(b_0=2t_c\). For \(\lambda\geq0\) and \(|z|\leq b_0\),
\begin{equation*}
 e^{\lambda z}\leq1+\lambda z+
 \frac{e^{\lambda b_0}-1-\lambda b_0}{b_0^2}z^2.
\end{equation*}
For positive \(z\), this follows from the Taylor series; for negative
\(z\), the remainder is at most \(\lambda^2z^2/2\), which is bounded
by the displayed quadratic term. Apply this inequality to both signs of
\(e\). Conditional centering gives
\begin{equation*}
 \mathbb E_\nu[\exp(\pm\lambda e\circ\mathcal F^j)\mid\mathcal G_{j+1}]
 \leq\exp\left[\frac{\sigma^2}{b_0^2}
 (e^{\lambda b_0}-1-\lambda b_0)\right].
\end{equation*}
On a finite block, condition successively on
\(\mathcal G_1,\ldots,\mathcal G_n\), removing the earliest remaining
term at each step. If \(v=n\sigma^2>0\), choose
\(\lambda=b_0^{-1}\log(1+b_0x/v)\). With
\(\Psi(s)=(1+s)\log(1+s)-s\), exponential Markov yields
\begin{equation*}
 \nu\left(\left|\sum_{j<n}e\circ\mathcal F^j\right|>x\right)
 \leq2\exp\left[-\frac{v}{b_0^2}\Psi(b_0x/v)\right]
 \leq2\left(\frac{\exp(1)v}{b_0x}\right)^{x/b_0}.
\end{equation*}
The last inequality uses
\((1+1/s)\log(1+s)-1\geq\log s-1\).
If \(v=0\), the centered sum vanishes. In the three variance regimes,
\(v/(b_0x)\) is bounded by
\(C_{q,H}/I\), \(C_{q,H}(1+\log I)/I\), and
\(C_{q,H}I^{-\zeta}\), respectively.
Put \(\beta_0=\min(1,\zeta)/2\), and choose the fixed \(q\) so that
\(\beta_0q/2>\Lambda\). Since \(x/b_0=q/2\), increasing the cutoff for
\(I\) makes the last bound at most \(2I^{-\Lambda}\).
The union probability of a discarded difference is at most
\(Cn t_c^{-r}b(\log(\mathrm e+t_c))^r\leq C_{q,H}I^{-\zeta}\).
On the no-discard event, the reward sum differs from the centered
truncated sum only by its compensator and the telescoping coboundary.
Their combined size is at most \(x\) for large \(I\).
This proves \eqref{eq:neutral-upper}; the direct union bound for \(V_u\)
at threshold \(\varepsilon nm_u\) proves \eqref{eq:neutral-last}.

For the lower tail take \(\lambda=\theta/a_u\), \(\theta>0\)
small. For \(z\geq0\), \(\varphi(z)=e^{-z}-1+z\) has derivative
at most \(\min(z,1)\). Conditional layer-cake integration bounds
\(\mathcal P\varphi(\lambda d_u^+)\) by
\begin{align*}
 &C\lambda^2m_u^2+C\lambda^2\int_0^{1/\lambda}t\min(u,t^{-r}b(\log(\mathrm e+t))^r)\,dt\\
 &\hspace{2em}+C\lambda\int_{1/\lambda}^\infty\min(u,t^{-r}b(\log(\mathrm e+t))^r)\,dt.
\end{align*}
The bounded negative part contributes at most \(C\lambda^2m_u^2\).
The resulting bound is \(Cu\psi_r(\theta)\), where \(\psi_r\) is
\(\theta^2\), \(\theta^2(1+\log(1/\theta))^{\rho+1}\), or
\(\theta^r(1+\log(1/\theta))^\rho\)
in the three regimes. Centering gives the one-step bound
\begin{equation*}
 \mathcal P\exp(-\lambda d_u)\leq\exp(Cu\psi_r(\theta)).
\end{equation*}
Transferring the earliest factor by
\(\mathcal P[H(g\circ\mathcal F)]=(\mathcal PH)g\) gives, inductively,
\begin{equation*}
 \mathcal P^n\exp\left(-\lambda\sum_{j<n}d_u\circ\mathcal F^j\right)
 \leq\exp(Cnu\psi_r(\theta)).
\end{equation*}
On the lower event, the sum in this exponential is at most
\(-\varepsilon nm_u+2\|\chi_u\|_\infty\leq-\varepsilon nm_u/2\)
for large \(n\). Positivity bounds its endpoint probability by
\(\exp[-\lambda\varepsilon nm_u/2+Cnu\psi_r(\theta)]\).
Since \(m_u/(ua_u)\) is bounded below positively and
\(\psi_r(\theta)/\theta\to0\), choosing \(\theta\) sufficiently small
proves \eqref{eq:neutral-endpoint-lower}. Integration and invariance
give \eqref{eq:neutral-lower} for every deterministic start.
Only the negative exponential moment is used.
\end{proof}

\begin{lemma}[Stopped long rewards]\label{lem:neutral-long-rewards}
For \(u<1/2\), integer \(L\geq1\), put
\(A_{u,L}=\{V_u\geq L\}\) and \(p_{u,L}=\nu(A_{u,L})\).
Then
\begin{equation*}
 p_{u,L}\asymp\mathbb E q_{L-1}(u),\qquad
 \|\mathcal P\mathbf1_{A_{u,L}}\|_{BV(Y)}\leq Cp_{u,L}.
\end{equation*}
The transferred indicator depends only on the base coordinate.
For any deterministic choices \(u_j,L_j\), if
\(E_j=\mathcal F^{-j}A_{u_j,L_j}\) and \(p_j=\nu(E_j)\), then
\begin{equation*}
 |\nu(E_i\cap E_j)-p_ip_j|\leq Cp_ip_j\vartheta^{j-i-1},\qquad
 \operatorname{Var}\sum_{j=l}^m\mathbf1_{E_j}\leq C\sum_{j=l}^mp_j.
\end{equation*}
\end{lemma}
\begin{proof}
The right landing must be at most the inverse threshold \(q_{L-1}(u)\)
formed from the next \(L-1\) left symbols. Since \(\widehat h\)
is uniformly positive and bounded, integration over
\([1/2,(1+q_{L-1}(u))/2]\) proves the probability comparison.
All these \(L\) observations occur before return, and
\(\tau\geq L+1\) on the event. Hence the observable is determined
by the consumed word; testing its stopped transfer against an unused-suffix
cylinder separates the iid suffix factor and proves the base-only assertion.

For a fixed infinite word, all selected full return atoms use the same
threshold \(q_{L-1}(u)\). The endpoint-ratio argument in
Lemma~\ref{lem:neutral-rewards} puts their whole right-landing intervals
in \([0,C_\alpha q_{L-1}(u)]\). Their union therefore has length
at most \(Cq_{L-1}(u)\), before averaging the partition.
On each atom the event indicator has variation at most one. The branch
product bound, and division by the positive Lipschitz \(\widehat h\),
give the asserted BV estimate.
Transfer the earlier event once, apply the base gap for the remaining
\(j-i-1\) returns, and integrate the resulting uniform error against
the later event. Its mean supplies \(p_j\), even though that event
depends on future symbols. Summing covariances proves the variance bound.
\end{proof}

\subsection*{Moving Laplace transforms and lower deviations}
The exact critical lower envelope requires a matching conditional
deviation exponent, rather than only an exponential upper bound.
The weighted operators below use a different target at each return.
We control their ordered products and the resulting conditional
terminal laws; a spectral estimate for a single fixed target would
not supply this uniformity.
Put \(k_B=\widehat h(1/2)/2\); the occupation identity gives
\(\overline c=k_B/\gamma\).

\begin{lemma}[Critical lower deviations]\label{lem:neutral-critical-rate}
Suppose \(u_n=c\log\log n/n\) eventually, with a positive finite prefix.
Put \(I_N=Nu_N\), \(a_N=a_{u_N}\), and fix \(\lambda>0\), integers
\(a\geq1,b_0\), and deterministic \(j_N\geq0\) with \(j_N=o(N)\).
Define
\begin{equation*}
 v_j=u_{\lfloor\lambda j\rfloor+a},\qquad
 n_N=\max\{0,\lfloor N/\lambda\rfloor+b_0\},\qquad
 S_N=\sum_{j=j_N}^{n_N-1}V_{v_j}\circ\mathcal F^{j-j_N}.
\end{equation*}
The limits below hold uniformly outside a null set of terminal base coordinates
with respect to \(\widehat h(y)\,dy\); the stopped-transfer kernel
constructed in the proof gives their base-valued representatives.
For each fixed \(\theta\geq0\),
\begin{equation}\label{eq:neutral-critical-pressure}
 \frac1{I_N}\log\bigl(\mathcal P^{n_N-j_N}e^{-\theta S_N/a_N}\bigr)
 \longrightarrow\Phi_\lambda(\theta)
 =-\frac{k_B}{\lambda\alpha}H_r(\theta D_*).
\end{equation}
For every \(0<x<\gamma/\lambda\),
\begin{equation}\label{eq:neutral-critical-rate}
 \frac1{I_N}\log\bigl(\mathcal P^{n_N-j_N}\mathbf1_{S_N\leq xM_N}\bigr)
 \longrightarrow-J_\lambda(x),\qquad
 J_\lambda(x)=\frac{k_B}{\lambda\alpha}h_r(x\lambda/\gamma).
\end{equation}
\end{lemma}
\begin{proof}
We first estimate one weighted return operator in the strong BV norm.
We then control the ordered product for the changing radii, including
its terminal density. Finally, exponential tilting gives the matching
conditional lower bound.

\emph{The weighted return operator.}
For a BV base input define
\(\mathcal A_{u,\theta}f=\mathcal P(e^{-\theta V_u/a_u}f)\).
The multiplier uses only the consumed return word, so this operator is
base-valued. On each complete return cell with nonzero reward, the
multiplier minus one has sup norm and variation at most one. The geometry
in Lemma~\ref{lem:neutral-rewards} places these whole cells in
\(B_{C_\alpha u}\). Their total and maximum lengths are at most \(Cu\).
Fix the environment first and let \(\mathcal I_u\) be these disjoint
return cells. The branch product bound for an arbitrary BV input is
\(C|I|(\|f\|_\infty+\operatorname{Var}_I f)\). Consequently,
\begin{equation*}
 \begin{aligned}
 \sum_{I\in\mathcal I_u}|I|
       (\|f\|_\infty+\operatorname{Var}_I f)
 &\leq\Bigl(\sum_{I\in\mathcal I_u}|I|\Bigr)\|f\|_\infty
       +\Bigl(\sup_{I\in\mathcal I_u}|I|\Bigr)\operatorname{Var}_Y f\\
 &\leq Cu\|f\|_{BV}.
 \end{aligned}
\end{equation*}
Here disjointness is used only within the fixed-environment partition;
we average the resulting bound afterwards. Conjugating by the positive
Lipschitz density \(\widehat h\), bounded away from zero, changes only
the constant. Thus the estimate is in the strong BV norm and gives
\begin{equation}\label{eq:neutral-weighted-perturbation}
 \|\mathcal A_{u,\theta}-\mathcal P\|_{BV\to BV}\leq Cu,
 \qquad\theta\geq0.
\end{equation}
The compact inverse profile \eqref{eq:neutral-profile} and its integrable
tail bound give, under \(\nu\),
\begin{equation*}
 \frac{\nu(V_u>ta_u)}u\to k_B(1+t/D_*)^{-r}\quad(t>0),\qquad
 \frac{m_u}{ua_u}\to\frac{k_BD_*}{r-1}.
\end{equation*}
Indeed the entrance interval for \(V_u>ta_u\) has length
\(q_{\lfloor ta_u\rfloor}(u)/2\), up to a one-step endpoint change;
\(\widehat h\) is Lipschitz at \(1/2\). Averaging the profile and
integrating its uniformly integrable tail proves both assertions.
Let \(W\) have tail \((1+t/D_*)^{-r}\), and set
\(\ell(\theta)=\mathbb Ee^{-\theta W}-1\). Layer cake and the mean
bound imply
\(u^{-1}\nu(e^{-\theta V_u/a_u}-1)\to k_B\ell(\theta)\), uniformly
on compact \(\theta\)-intervals.

The base gap in Proposition~\ref{prop:renewal} and the strong-norm bound
permit a resolvent Neumann series on fixed separating contours.
Writing \(\Pi f=\nu(f)1\) on base functions, the resulting eigenvalue
and projection satisfy, with projection error in \(BV\to BV\),
\begin{equation*}
 \eta_{u,\theta}=1+\nu(e^{-\theta V_u/a_u}-1)+O(u^2),\qquad
 \Pi_{u,\theta}=\Pi+O(u).
\end{equation*}
For the first formula normalize the right eigenfunction
\(g_{u,\theta}\) to have \(\nu\)-mean one; it equals \(1+O(u)\) in BV.
Since \(\nu\circ\mathcal P=\nu\) and \(\nu(g_{u,\theta}-1)=0\), applying
\(\nu\) to the eigen-equation gives
\begin{equation*}
 \eta_{u,\theta}
 =\nu(\mathcal A_{u,\theta}1)
  +\nu\bigl((\mathcal A_{u,\theta}-\mathcal P)(g_{u,\theta}-1)\bigr)
 =1+\nu(e^{-\theta V_u/a_u}-1)+O(u^2),
\end{equation*}
by \eqref{eq:neutral-weighted-perturbation}.
The remaining powers are bounded by \(C\rho^k\), \(\rho<1\), uniformly.
Since \(\Pi_{u,\theta}1=1+O(u)>0\), the relative error is
\(O((\rho/\eta_{u,\theta})^k)\). Hence, when \(ku\to\infty\),
\begin{equation*}
 (ku)^{-1}\log(\mathcal A_{u,\theta}^k1)\to k_B\ell(\theta)
\end{equation*}
uniformly at endpoints and on compact tilt intervals.

\emph{Ordered products for moving rewards.}
For \(\mathcal B_j f=\mathcal P(e^{-\theta V_{v_j}/a_N}f)\),
the transfer identity gives the ordered product
\begin{equation*}
 \mathcal P^{n_N-j_N}e^{-\theta S_N/a_N}
 =\mathcal B_{n_N-1}\cdots\mathcal B_{j_N}1.
\end{equation*}
Fix a finite partition of \([\delta,1]\). On the indices belonging to
each partition interval, bracket the moving rewards by the radii at the
two endpoints of that interval, keeping the single
denominator \(a_N\). For a frozen radius at fraction \(t\),
\(a_{u_{\lfloor tN\rfloor}}/a_N\to t^\alpha\), and block length
times radius tends to infinity. For each block product \(\mathcal T_i\),
frozen transfer powers bound \(L_i1\leq\mathcal T_i1\leq U_i1\).
Positivity then gives
\begin{equation*}
 \Bigl(\prod_iL_i\Bigr)1\leq\mathcal T_m\cdots\mathcal T_1 1
 \leq\Bigl(\prod_iU_i\Bigr)1.
\end{equation*}
Thus the factors need not commute. First let \(N\to\infty\) for fixed
\(\delta\) and partition, then refine that partition. This
gives \((k_B/\lambda)\int_\delta^1\ell(\theta t^\alpha)\,dt/t\).
For \(m=\lfloor\delta N/\lambda\rfloor\), conditional Jensen and
\eqref{eq:neutral-reward-gap} bound the early transform below by
\begin{equation*}
 \mathcal P^{m-j_N}\exp\left(-\frac\theta{a_N}
 \sum_{j=j_N}^{m-1}V_{v_j}\circ\mathcal F^{j-j_N}\right)
 \geq\exp\left(-\frac{C\theta}{a_N}\sum_{j<m}m_{v_j}\right).
\end{equation*}
Its upper bound is one. Karamata makes this logarithmic error, divided
by \(I_N\), at most \(C\theta\delta^\alpha+o(1)\). Here each mean is
transferred by \(\mathcal P^{m-j}\) with \(m-j\geq1\), so the bound
\eqref{eq:neutral-reward-gap} applies. The same holds for \(j_N=o(N)\).
Let \(\delta\downarrow0\). Tonelli gives
\(\int_0^1\ell(\theta t^\alpha)\,dt/t=-H_r(\theta D_*)/\alpha\),
proving \eqref{eq:neutral-critical-pressure}.

\emph{Conditional kernels and exponential tilting.}
For completeness the stopped transfer defines an actual conditional
probability kernel. With \(K\) returns, sum over total consumed lengths
\(n\geq K\), integrate their words against \(\pi^n\), and sum over
complete \(K\)-return inverse cells \(\psi_I\), using weights
\begin{equation*}
 \frac{\widehat h(\psi_I(y))|\psi_I'(y)|}{\widehat h(y)}.
\end{equation*}
They sum to \(\mathcal P^K1(y)=1\). The unused iid suffix is independent,
so every Borel event of the consumed rewards has this base-valued kernel.

Write \(X_N=S_N/(I_Na_N)\) and \(\Phi=\Phi_\lambda\). Since \(H_r\) is
concave and \(H_r'\) decreases strictly from \(1/(r-1)\) to zero, \(\Phi\)
is convex and \(s_\theta=-\Phi'(\theta)\) decreases continuously from
\(s_0=k_BD_*/[\lambda\alpha(r-1)]\) to zero.
For the upper bound, exponential Markov gives, for every
\(\theta\geq0\) and \(s>0\),
\begin{equation*}
 \mathcal P^{n_N-j_N}\mathbf1_{X_N\leq s}
 \leq e^{\theta sI_N}\,\mathcal P^{n_N-j_N}e^{-\theta S_N/a_N},
\end{equation*}
so \eqref{eq:neutral-critical-pressure} yields
\(\limsup I_N^{-1}\log\mathcal P^{n_N-j_N}\mathbf1_{X_N\leq s}
\leq-\sup_{\theta\geq0}\{-\theta s-\Phi(\theta)\}\).
For the lower bound fix \(\theta>0\) and \(\epsilon>0\), and let
\(\mathbb Q_N\) be the kernel tilted by
\(e^{-\theta S_N/a_N}/\mathcal P^{n_N-j_N}e^{-\theta S_N/a_N}\).
For \(0<z<\theta\), exponential Markov under \(\mathbb Q_N\) and
\eqref{eq:neutral-critical-pressure} at \(\theta\) and \(\theta-z\) give
\begin{equation*}
 \limsup_{N\to\infty}\frac1{I_N}\log\mathbb Q_N(X_N>s_\theta+\epsilon)
 \leq\Phi(\theta-z)-\Phi(\theta)-z(s_\theta+\epsilon),
\end{equation*}
and similarly, for \(z>0\),
\begin{equation*}
 \limsup_{N\to\infty}\frac1{I_N}\log\mathbb Q_N(X_N<s_\theta-\epsilon)
 \leq\Phi(\theta+z)-\Phi(\theta)+z(s_\theta-\epsilon).
\end{equation*}
Since \(\Phi(\theta\mp z)-\Phi(\theta)=\pm zs_\theta+o(z)\), both bounds
are negative for small \(z\); hence
\(\mathbb Q_N(|X_N-s_\theta|\leq\epsilon)\to1\) uniformly in the terminal
coordinate. Undoing the tilt on this event,
\begin{equation*}
 \mathcal P^{n_N-j_N}\mathbf1_{X_N\leq s_\theta+\epsilon}
 \geq e^{\theta(s_\theta-\epsilon)I_N}\,
 \bigl(\mathcal P^{n_N-j_N}e^{-\theta S_N/a_N}\bigr)\,
 \mathbb Q_N(|X_N-s_\theta|\leq\epsilon),
\end{equation*}
so \(\liminf I_N^{-1}\log\mathcal P^{n_N-j_N}\mathbf1_{X_N\leq s_\theta+\epsilon}
\geq\Phi(\theta)+\theta s_\theta-\theta\epsilon\). By convexity,
\(-\Phi(\theta)-\theta s_\theta=\sup_{\theta'\geq0}\{-\theta' s_\theta-\Phi(\theta')\}\).
Given \(s\in(0,s_0)\) and small \(\epsilon>0\), choose the unique
\(\theta>0\) with \(s_\theta=s-\epsilon\). Letting \(\epsilon\downarrow0\)
and using continuity of the rate in \(s\) shows that the two bounds agree,
with rate \(\sup_{\theta\geq0}\{-\theta s-\Phi(\theta)\}\) on \((0,s_0)\).
All Laplace arguments remain nonnegative; no positive exponential moment
of \(W\) is required. Finally
\(M_N/(I_Na_N)\to k_BD_*/[\gamma\alpha(r-1)]\), so the threshold
is \(s=xk_BD_*/[\gamma\alpha(r-1)]\). Substituting
\(q=\theta D_*\) in the supremum gives \(J_\lambda(x)\), precisely
for \(0<x<\gamma/\lambda\). Bracket the convergent deterministic
threshold and use continuity of the rate to obtain
\eqref{eq:neutral-critical-rate}.
\end{proof}

\begin{proof}[Proof of Theorems~\ref{thm:neutral-counts} and~\ref{thm:neutral-log}]
Lemma~\ref{lem:neutral-normalizer} proves all normalization assertions.
It remains to estimate the physical counts.
First assume eventual monotonicity of the radii. We transfer concentration
to the physical clock, use long rewards and survival for the two failure
regimes, and separate \(L^1\) convergence from the finer lower envelopes.
The matching conditional rate identifies the critical constant. Finally,
comparison with smooth schedules removes monotonicity from the general
assertions, and continuity transfers the critical constant to
asymptotically equivalent schedules.

\subsection{Strong laws and the physical clock}\label{sec:neutral-strong-proof}
Fix \(a>1\) and \(N_k=\lfloor a^k\rfloor\).
Uniform convergence for \(I(x)=xu(x)\), of index \(1-g\), shows that
\(J(u)<\infty\) if and only if
\begin{equation}\label{eq:neutral-critical-series}
 \sum_k I_{N_k}^{-\zeta}<\infty.
\end{equation}
For \(g<1\), Potter bounds imply convergence. In general, convergence
implies \(I_n\to\infty\), first on the grid, then between its points by
uniform convergence. Choose the deviation exponent \(\Lambda\geq\zeta\):
the remainder \(I_{N_k}^{-\Lambda}\) is dominated by the same summable series.
At any fixed return-index multiple and neighboring endpoint radius,
\(nu\asymp I_{N_k}\) and \(\log(nu)/\log(1/u)\) tends to
\((1-g)/g\), or zero at \(g=1\). Thus a fixed \(H\) permits
\eqref{eq:neutral-upper}--\eqref{eq:neutral-last} on all these grids.
Throughout the proof, \(\lambda\) and \(\Lambda_\pm\) denote numbers of
physical steps per return, to be compared with the mean return \(\gamma\).
Apply the deviation bounds at every return index \(\lfloor N_k/\lambda\rfloor\),
\(\lambda>0\) rational,
and both neighboring endpoint radii, on countably many grids \(a\downarrow1\).
First Borel--Cantelli gives centered-sum and maximum-reward errors tending
to zero relative to \(N_km_u\). For rational \(\Lambda_-<\gamma<\Lambda_+\),
the return clock sandwiches the random
return index \(K_{N_k}\) between \(N_k/\Lambda_+\) and \(N_k/\Lambda_-\).
The unfinished excursion
is bounded by one further reward using \eqref{eq:neutral-last}.
Letting \(\Lambda_\pm\to\gamma\) proves for the frozen physical count
\(H_N(u)=\sum_{n=1}^N\mathbf1_{[0,u]}(T_\omega^n x)\) that
\(H_{N_k}(u)\sim N_k\overline\mu([0,u])\).
For each block \((N_k,N_{k+1}]\), these estimates hold at both times
\(N_k,N_{k+1}\) for each radius \(u_{N_k},u_{N_{k+1}}\),
on the countable grids just specified.

Every positive initial point enters \(Y\) in finite time since
\(T_s\geq T_B\) below \(Y\). An absolutely continuous initial law
has an absolutely continuous base-entry distribution and unused iid
future, so it inherits the \(\nu\)-full event. Include
its finite entry offset in \(\mathsf T_j\); the frozen count changes by only
finitely many observations. Eventual monotonicity of the radii gives
\begin{align*}
 H_{N_{k+1}}(u_{N_{k+1}})-H_{N_k}(u_{N_{k+1}})
 &\leq C_{N_{k+1}}-C_{N_k}\\
 &\leq H_{N_{k+1}}(u_{N_k})-H_{N_k}(u_{N_k}).
\end{align*}
The four frozen-count errors are
\(o(N_kw(u_{N_k}))\). Summing them is \(o(M_{N_K})\):
Potter bounds with positive exponent \(1-\xi\) dominate past normalized
scales by a summable geometric sequence, while their relative errors vanish.
The two frozen mass sums bracket \(M_{N_K}\); their relative difference
is \(O(a^\xi-1)+o(1)\), by regular variation of
\(\overline\mu([0,u])\sim C_*w(u)\).
Let \(a\downarrow1\). Increasing-count interpolation and regular
variation of \(M_N\) give \(C_N/M_N\to1\) for all \(N\).

\subsection{Deep entrances and infinite upper limits}\label{sec:neutral-deep-entrances}
For any choices in Lemma~\ref{lem:neutral-long-rewards} with
\(\sum_jp_j=\infty\), choose consecutive blocks with expected counts
at least \(C2^{k}\). By the variance bound of that lemma and Chebyshev's
inequality, the events \(D_k\) that a block contains some \(E_j\) satisfy
\(\nu(D_k^c)\leq2^{-k}\). The first Borel--Cantelli lemma gives
infinitely many \(E_j\), almost surely.

Suppose \(g=1\), \(J(u)=\infty\) and \(I_n\to\infty\).
Fix \(D>\gamma\) and integer
\(K\geq1\). Set \(n_j=\lceil2Dj\rceil\),
\(L_j=\lceil KM_{n_j}\rceil\), and
\(E_j=\mathcal F^{-j}A_{u_{n_j},L_j}\).
Then \(L_j/a_{u_{n_j}}\asymp n_ju_{n_j}\to\infty\) and
\(\log L_j\sim\alpha\log n_j\). Thus
\eqref{eq:neutral-deep-profile} applies. Its averaged version follows
from the band \([\alpha,\alpha+1/\log L_j]\): on its good event
the normalized inverse threshold is bounded. On the exceptional event,
it is at most \(u_{n_j}L_j^r/b(\log L_j)^r\), whose product with
the probability bound \(Ce^{-cL_j/b(\log L_j)}\) tends to zero.
Lemma~\ref{lem:neutral-long-rewards} therefore gives
\begin{equation*}
 p_j\asymp L_j^{-r}b(\log L_j)^r
 \asymp\frac1{jL(j)^\zeta}.
\end{equation*}
Indeed \(b(\log L_j)/b(\log n_j)\to\alpha^{\kappa_0}\):
the slowly varying factors contributed by the parameter law to the mean
and to the reward tail cancel.
Uniform convergence compares this series with \(J(u)\), so it diverges.
For logarithmic schedules with \(-1\leq\kappa\leq0\), the compact
profile \eqref{eq:neutral-profile}, bounded by \(q/u\leq1\), instead gives
\(p_j\asymp u_{n_j}\asymp j^{-1}(\log j)^\kappa\), also divergent:
\(L_j/a_{u_{n_j}}\) tends to a positive constant when \(\kappa=0\),
and tends to zero when \(-1\leq\kappa<0\).
In either case,
since \(\mathsf T_j\leq Dj\) eventually and \(L_j=o(j)\),
\(\mathsf T_j+L_j\leq n_j\). On \(E_j\), the next \(L_j\)
observations are below \(u_{n_j}\); monotonicity of the radii makes
them actual moving-target hits. At \(N=\mathsf T_j+L_j\),
\(C_N\geq L_j\geq KM_{n_j}\geq KM_N\).
Intersecting over integers \(K\) proves the infinite upper limit.

\subsection{The annealed \texorpdfstring{\(L^1\)}{L1} criterion}\label{sec:neutral-l1}
For \(g<1\), the strong law just proved gives \(C_N/M_N\to1\)
in stationary probability. For \(g=1\), Theorem~\ref{thm:neutral-occupation},
proved below from the reward profiles, normalization and
\eqref{eq:neutral-weighted-perturbation}, without using the present
\(L^1\) conclusion, gives this convergence
when \(I_N\to\infty\). In either case stationarity gives
\(\mathbb E(C_N/M_N)=1\), so the bounded negative part tends to zero and
\(\mathbb E|C_N/M_N-1|=2\mathbb E(1-C_N/M_N)^+\to0\).
Conversely, if \(I_N\) does not tend to infinity at \(g=1\), there is a
deterministic subsequence on which \(I_N\to\iota<\infty\).
The occupation theorem gives zero or the nondegenerate \(Z_\iota(1)\),
not the constant one, ruling out \(L^1\) convergence. This argument covers
nonmonotone radii.

\subsection{The positive lower envelope}\label{sec:neutral-positive-envelope}
First suppose \(g=1\) and \(I_n\geq c_0\log\log n\) eventually,
where \(c_0>0\). Fix \(\Lambda_-<\gamma<\Lambda_+\) with
\(D=2/\Lambda_+-1/\Lambda_->0\). At \(N_k=2^k\) use a fixed number
\(J\) of physical windows \([t_iN_k,2t_iN_k]\),
\(t_i=2^{i-J}\), \(0\leq i<J\). Their complete return segments are
\begin{equation*}
 s_i=\lceil t_iN_k/\Lambda_-\rceil,\qquad
 e_i=\lfloor2t_iN_k/\Lambda_+\rfloor,\qquad
 l_i=e_i-s_i,\qquad v_i=u_{\lfloor2t_iN_k\rfloor}.
\end{equation*}
For large \(k\), these return-index segments are disjoint and
\(l_iv_i/I_{N_k}\to D/2\). Define the low-reward events on
\([s_i,e_i)\) by
\(E_i=\{\sum_{j=s_i}^{e_i-1}V_{v_i}\circ\mathcal F^j<l_im_{v_i}/2\}\).
Every later event is measurable at the earlier endpoint. Conditioning
out the earliest event using \eqref{eq:neutral-endpoint-lower}, then
iterating over the disjoint segments, gives
\begin{equation*}
 \nu\Bigl(\bigcap_{i<J}E_i\Bigr)
 \leq\exp\Bigl(-c_{1/2}\sum_{i<J}l_iv_i\Bigr)
 \leq\exp(-c_{1/2}JD I_{N_k}/4).
\end{equation*}
Since \(\log\log N_k=\log k+O(1)\), choose \(J\) so that
\(c_{1/2}JDc_0/4>2\), before letting \(k\to\infty\).
First Borel--Cantelli makes at least one segment good for every large
\(k\). On the full-measure event where \(\mathsf T_j/j\to\gamma\), each
complete segment lies inside its physical window; strict margins absorb any finite entry offset.
The frozen radius then gives actual moving-target hits. Moreover,
\begin{equation*}
 \frac{l_im_{v_i}}{2M_{N_k}}
 \longrightarrow\frac{\alpha\gamma D}{2}\,2^{\alpha-1}t_i^\alpha>0.
\end{equation*}
The minimum over fixed \(J\) windows is positive. Increasing-count
interpolation and \(M_{N_{k+1}}/M_{N_k}\to2^\alpha\) prove a positive
lower envelope.

Suppose \(g=1\) and \(I_n/\log\log n\to\infty\). This includes
all logarithmic schedules with \(\kappa>0\) and all strong-law-threshold
schedules in Theorem~\ref{thm:neutral-log}. Fix \(a>1\), put \(N_k=\lfloor a^k\rfloor\),
and choose rational \(\Lambda_-<\gamma<\Lambda_+\)
with \(a/\Lambda_+>1/\Lambda_-\). The return-index segment from
\(s_k=\lceil N_k/\Lambda_-\rceil\) through
\(e_k=\lfloor N_{k+1}/\Lambda_+\rfloor-1\) has length \(l_k\asymp N_k\).
Use \eqref{eq:neutral-lower} at radius \(u_{N_{k+1}}\).
Its failure probabilities are bounded by \(\exp(-c_\varepsilon I_{N_k})\).
Since \(\log\log N_k=\log k+O(1)\), the assumed intensity condition
makes this bound at most \(k^{-2}\) eventually, hence summable.
Eventually its reward sum is at least \((1-\varepsilon)l_km_{u_{N_{k+1}}}\).
The clock gives \(\mathsf T_{s_k}>N_k\), \(\mathsf T_{e_k+1}<N_{k+1}\);
these entire excursions lie inside the physical block and give actual hits.
Their length factor is
\(\beta_a=\gamma(a/\Lambda_+-1/\Lambda_-)/(a-1)\).
Endpoint masses have ratio tending to \(a^{1-\alpha}\).
Summing the lower frozen block masses and interpolating all \(N\) gives
\(\liminf C_N/M_N\geq(1-\varepsilon)\beta_a/a\).
Take the rational margins to \(\gamma\), then countably
\(a\downarrow1\) and \(\varepsilon\downarrow0\). Thus the lower limit
is at least one.

Under the stationary physical law \(\mathbb EC_N=M_N\); Fatou therefore
forces this lower limit to equal one. The \(L^1\) conclusion was proved above.

\subsection{Conditional survival and deficient counts}\label{sec:neutral-survival}
Here \(g=1\). We prove a strictly sub-one lower envelope when
\(I_n/\log\log n\) has finite upper limit, and zero both when this
ratio tends to zero and when \(I_n\) has a bounded subsequence.
Finite changes of the schedule are harmless, so assume the decreasing
radii are all below \(1/2\). Fix \(0<\lambda<\gamma\), and put
\(w_j=u_{\lfloor\lambda j\rfloor+1}\),
\(n_N=\lfloor N/\lambda\rfloor+1\).
On the full-measure event where \(\mathsf T_j/j\to\gamma\), we have \(\mathsf T_j\geq\lambda j\) eventually, including
any finite entry offset. Every later neutral observation in excursion
\(j\) occurs at an integer time at least \(\lfloor\lambda j\rfloor+1\),
and \(K_N\leq N/\lambda\) eventually. Thus, for a pathwise finite
\(B_0\) and all large \(N\),
\begin{equation}\label{eq:neutral-complete-domination}
 C_N\leq B_0+\sum_{j<n_N}V_{w_j}\circ\mathcal F^j.
\end{equation}
The finite \(B_0\) absorbs entry observations and complete rewards of the
finitely many exceptional early excursions.
Complete rewards on the right overcount the unfinished excursion.

Let \(\mathcal G\) be the full return-skew sigma-algebra. For
\(G_{w,h}=\prod_{i<h}\mathbf1_{B_w^c}\circ\mathcal F^i\),
the identity \(\mathcal P[H(g\circ\mathcal F)]=(\mathcal PH)g\)
and Lemma~\ref{lem:neutral-rewards} give, inductively,
\begin{equation}\label{eq:neutral-survival}
 \mathcal P^hG_{w,h}\geq(1-C_Bw)^h\geq\exp(-2C_Bhw)
 \quad(C_Bw\leq1/2).
\end{equation}
For successive disjoint return-index blocks \([a_k,b_k)\), let
\(E_k=\mathcal F^{-a_k}H_k\). If
\(\mathcal P^{b_k-a_k}\mathbf1_{H_k}\geq p_k\), then
\begin{equation}\label{eq:neutral-reverse-block}
 \nu\Bigl(\bigcap_{k=m}^nE_k^c\Bigr)\leq\prod_{k=m}^n(1-p_k).
\end{equation}
Indeed, every later factor is \(\mathcal F^{-b_m}\mathcal G\)-measurable,
and conditioning the earliest one gives
\(\mathbb E[\mathbf1_{E_m}\mid\mathcal F^{-b_m}\mathcal G]
=(\mathcal P^{b_m-a_m}\mathbf1_{H_m})\circ\mathcal F^{b_m}\geq p_m\).
Iterate. Thus \(\sum_kp_k=\infty\) implies \(E_k\) infinitely often.

Choose an integer \(p>4/\alpha\). In either iterated-logarithm regime take
\(N_k=(k!)^p\); in the bounded-intensity case choose a deterministic
subsequence with \(I_{N_k}\leq I_0\) and \(N_k/N_{k-1}\geq k^p\).
Put \(a_k=n_{N_{k-1}}\), \(b_k=n_{N_k}\). By regular-variation summation,
\(\sum_{j<n_N}m_{w_j}\sim(\gamma/\lambda)M_N\).
The old prefix \(A_k=\sum_{j<a_k}V_{w_j}\circ\mathcal F^j\)
therefore has expectation asymptotic to
\((\gamma/\lambda)M_{N_{k-1}}\). Potter bounds give
\begin{equation*}
 \frac{\nu(A_k)}{M_{N_k}}\leq C(N_{k-1}/N_k)^{\alpha/2}
 \leq Ck^{-p\alpha/2}.
\end{equation*}
Markov and first Borel--Cantelli imply \(A_k/M_{N_k}\to0\) almost surely.

Fix \(0<t<1\), and put
\(L_k=\lfloor tN_k/\lambda\rfloor\), \(l_k=L_k-a_k\),
\(h_k=b_k-L_k\). Eventually \(l_k>0\). In local block coordinates put
\(S_k=\sum_{j=a_k}^{L_k-1}V_{w_j}\circ\mathcal F^{j-a_k}\).
Every transfer exponent is positive, so \eqref{eq:neutral-reward-gap} gives
\begin{equation*}
 \begin{aligned}
 \mathcal P^{l_k}S_k
 &=\sum_{j=a_k}^{L_k-1}\mathcal P^{L_k-j}V_{w_j}
 =\sum_{j=a_k}^{L_k-1}m_{w_j}+\Delta_k,\\
 \|\Delta_k\|_\infty&\leq C\sum_{j=a_k}^{L_k-1}
 \vartheta^{L_k-j-1}m_{w_j}=o(M_{N_k}).
 \end{aligned}
\end{equation*}
Indeed, terms with \(j<L_k/2\) have exponentially small weight times
total mean \(O(M_{N_k})\); the remaining means are bounded by
\(Cm_{w_{L_k}}=O(M_{N_k}/N_k)\), and their geometric weights sum to
a constant. Karamata and the negligible old mean therefore give
\(\mathcal P^{l_k}S_k/M_{N_k}\to A_t=(\gamma/\lambda)t^\alpha\)
uniformly. For \(q>A_t\), pointwise Markov implies
\begin{equation*}
 \mathcal P^{l_k}\mathbf1_{S_k\leq qM_{N_k}}\geq\eta_0>0
\end{equation*}
eventually, for any fixed \(\eta_0<1-A_t/q\). Define
\begin{equation*}
 H_k=\{S_k\leq qM_{N_k}\}
 \cap\mathcal F^{-l_k}\{G_{w_{L_k},h_k}=1\}.
\end{equation*}
On \(E_k=\mathcal F^{-a_k}H_k\), every reward with
\(L_k\leq j<b_k\) vanishes, since \(w_j\leq w_{L_k}\).
Transfer the prefix indicator first,
multiply by survival, then transfer the suffix. Positivity and
\eqref{eq:neutral-survival} yield
\begin{equation}\label{eq:neutral-block-small}
 \begin{aligned}
 \mathcal P^{b_k-a_k}\mathbf1_{H_k}
 &\geq\eta_0\exp(-2C_Bh_kw_{L_k})=:p_k,\\
 \frac{h_kw_{L_k}}{I_{N_k}}\longrightarrow\frac{1-t}{\lambda t}.
 \end{aligned}
\end{equation}
Suppose first \(I_n\leq c_1\log\log n\) eventually, with \(c_1<\infty\).
Choose \(t\) close to one, then \(\gamma t^\alpha<\lambda<\gamma\)
and \(A_t<q<1\), so that
\(A_1=2C_Bc_1(1-t)/(\lambda t)<1\). For example,
\(\lambda/\gamma=(1+t^\alpha)/2\) allows both inequalities as
\(t\uparrow1\). Since \(\log\log N_k=\log k+\log\log k+O(1)\),
for any \(A_1<v<1\) we have \(p_k\geq\eta_0k^{-v}\) eventually.
The conditional block argument gives infinitely many successes.
Together with \eqref{eq:neutral-complete-domination} and the negligible
old prefix, this yields \(\liminf C_N/M_N\leq q<1\).

For zero lower limits, fix any \(q>0\), then \(\lambda=\gamma/2\)
and \(t\) with \(2t^\alpha<q\). If \(I_n/\log\log n\to0\), the
same factorial grid gives \(p_k\geq\eta_0k^{-1/2}\) eventually.
On a bounded-intensity subsequence, \(p_k\) has a fixed positive lower
bound. The block argument gives \(\liminf C_N/M_N\leq q\) in both
cases; intersect over countably \(q\downarrow0\). An absolutely continuous
initial law inherits these full events after its finite first entry,
and the finite \(B_0\) vanishes relative to \(M_N\).

For \(L(n)=(\log n)^\kappa\), \(\kappa<0\),
Theorem~\ref{thm:neutral-occupation} gives \(C_N/M_N\to0\)
in stationary annealed probability.

If \(\kappa<-1\), entrance probabilities have finite sum. The first Borel--Cantelli lemma gives
finitely many entrances. A left step cannot enter a decreasing neutral target
from outside, and each positive excursion is finite; hence only finitely many
visits occur.

\subsection{The iterated-logarithm lower-envelope constant}\label{sec:neutral-critical-constant}
Let the schedule be eventually equal to
\(u_n=c\log\log n/n\). First note
that \(H_r'(q)=\int_0^\infty e^{-qz}(1+z)^{-r}\,dz\) decreases
strictly from \(1/(r-1)\) to zero and \(H_r(q)=\log q+O(1)\)
as \(q\to\infty\). Thus \(h_r\) is continuous and strictly decreasing
from infinity to zero on \((0,1)\); its optimizer is unique and positive.
Let \(S_N^{\lambda,a,b_0}\) denote the sum in
Lemma~\ref{lem:neutral-critical-rate} with \(j_N=0\), and define
\begin{equation*}
 D_\lambda(c)=\frac\gamma\lambda
 h_r^{-1}\left(\frac{\lambda\alpha}{ck_B}\right).
\end{equation*}
We claim \(\liminf S_N^{\lambda,a,b_0}/M_N=D_\lambda(c)\), \(\nu\)-almost surely.
For \(x<D_\lambda(c)\), \eqref{eq:neutral-critical-rate} makes
\(S_{N_k}\leq xM_{N_k}\) summable on \(N_k=\lfloor\rho^k\rfloor\),
since \(I_{N_k}/\log k\to c\) and \(cJ_\lambda(x)>1\).
Increasing-count interpolation gives a lower bound \(x\rho^{-\alpha}\).
Choose countably \(x\uparrow D_\lambda(c)\), \(\rho\downarrow1\).

For \(x>D_\lambda(c)\), choose
\(D_\lambda(c)<x'<\min(x,\gamma/\lambda)\) and integer \(p>4/\alpha\).
On \(N_k=(k!)^p\), use consecutive return blocks
\([a_k,b_k)=[n_{N_{k-1}},n_{N_k})\). The old prefix mean is
\(O(M_{N_{k-1}})\); Potter bounds, Markov's inequality and the first
Borel--Cantelli lemma make its ratio to
\(M_{N_k}\) tend to zero almost surely. Here \(a_k=o(N_k)\).
Since \(cJ_\lambda(x')<1\), \eqref{eq:neutral-critical-rate} bounds the
new block's conditional success probability below by \(k^{-v}\),
for some \(v<1\). In local coordinates let \(H_k\) be the event that
this new-block reward is at most \(x'M_{N_k}\). The endpoint lower bound
survives composition at its global start \(\mathcal F^{a_k}\).
Apply \eqref{eq:neutral-reverse-block} to
\(E_k=\mathcal F^{-a_k}H_k\): successes occur infinitely often. Adding the
negligible old prefix proves the claimed upper bound; let \(x\downarrow D_\lambda(c)\).

Fix \(0<\Lambda_-<\gamma<\Lambda_+\) and set \(a_+=\lceil\Lambda_+\rceil+2\).
Ordinary clock margins give, for a pathwise finite \(B_0\),
\begin{equation}\label{eq:neutral-critical-clocks}
 S_N^{\Lambda_+,a_+,-3}-B_0\leq C_N
 \leq S_N^{\Lambda_-,1,1}+B_0.
\end{equation}
Indeed \(\mathsf T_j\geq\Lambda_-j\) bounds every later hit by the upper
radius; the whole possible unfinished reward is included. On the lower
side \(\mathsf T_{j+1}\leq\Lambda_+(j+1)\), and every included
\(j\leq\lfloor N/\Lambda_+\rfloor-4\) has its target index later
than the excursion endpoint but at most
\(N-4\Lambda_++\lceil\Lambda_+\rceil+2<N\), since \(\Lambda_+>\gamma\geq1\).
Those complete rewards are actual hits. Finite entry offsets and early
excursions are absorbed in \(B_0\).
Continuity of \(D_\lambda(c)\) at \(\gamma\) and countable clock margins
therefore yield \(\liminf C_N/M_N=D_\gamma(c)=d_\pi(c)\).
Its defining equation is the one in the statement. Since \(h_r\)
decreases continuously from infinity to zero, \(d_\pi\) increases
continuously and strictly from zero to one.

\subsection{Comparison without target monotonicity}\label{sec:neutral-comparison}
Use the smooth equivalent \(\widetilde L\) and decreasing comparison radii
\(v_n^\pm=(1\pm\epsilon)cn^{-g}\widetilde L(n)\) from the normalizer proof,
for rational \(\epsilon\in(0,1)\). Their counts satisfy
\(C_N^- -O_\epsilon(1)\leq C_N\leq C_N^+ +O_\epsilon(1)\).
Their means satisfy \eqref{eq:neutral-sandwich-means}; integral convergence,
\(I_n\to\infty\), bounded-intensity subsequences, and the
hypotheses on \(I_n/\log\log n\) are preserved.
If \(J(u)<\infty\), the two comparison strong laws squeeze
\(C_N/M_N\) between \((1\pm\epsilon)^{1-\alpha}\); intersect countably many
full-measure sets and let \(\epsilon\downarrow0\). On the divergent side with
\(I_n\to\infty\), one lower comparison gives infinite upper ratio.
If \(I_n\) does not tend to infinity, one upper comparison gives zero
lower ratio. The same upper comparison gives zero when
\(I_n/\log\log n\to0\). When this ratio tends to infinity, countably
many lower comparisons give lower limit at least one; stationary Fatou
again makes it equal one.
Positive lower critical intensity passes to any lower comparison and
gives a positive original lower envelope. For bounded upper critical
intensity, first fix \(0<t<1\), \(\gamma t^\alpha<\lambda<\gamma\)
and \(A_t<q<1\) for
\(cn^{-1}\widetilde L(n)\). Choose a sufficiently small upper factor
\(1+\epsilon\) so that the same survival exponent remains below one
and \((1+\epsilon)^{1-\alpha}q<1\). Its conditional prefix mean ratio
\(A_t\) is unchanged, so the same \(q\) bounds the comparison's
lower ratio. The mean factor then gives an original lower limit below one.
For a fixed \(u_n\sim c\log\log n/n\), compare instead with exact
critical schedules of coefficients \(c(1\pm\epsilon)\). Their mean
ratios tend to \((1\pm\epsilon)^{1-\alpha}\), so the lower envelope
is bracketed by these factors times \(d_\pi(c(1\pm\epsilon))\).
Continuity and countably \(\epsilon\downarrow0\) give the stated equality,
even for nonmonotone radii. The same mean comparison gives \(Q_N/M_N\to1\).

The logarithmic and strong-law-threshold arguments above use their asymptotic
sizes: these are unchanged by smooth equivalence and fixed radius factors.
In particular, their geometric-grid lower-tail series remain summable,
their compact-profile entrance series have the same convergence status,
and finite entrance sums stay finite. The same sandwich therefore
preserves all the stated finer classifications. The general integral
equivalence follows from the convergent side and the two failure cases.

Finally \(\overline\mu\) is equivalent to Lebesgue measure.
Disintegration over future symbols transfers all stationary full-measure
assertions to almost every environment and Lebesgue-almost every point;
Lemma~\ref{lem:neutral-normalizer} replaces \(M_N\) by \(Q_N\).
\end{proof}

\subsection{The occupation-process limit}\label{sec:neutral-occupation}
\begin{proof}[Proof of Theorem~\ref{thm:neutral-occupation}]
We first compute transforms at deterministic return horizons, compare
physical prefix counts, and then prove scalar \(M_1\) tightness.
Write \(a_N=a_{u_N}\). Endpoint normalization and Karamata give
\begin{equation}\label{eq:occupation-mean-scale}
 M_N\sim\frac{\overline cD_*}{\alpha(r-1)}I_Na_N,\qquad
 \frac{M_{\lfloor Nt\rfloor}}{M_N}\to t^\alpha\quad(t>0).
\end{equation}
In particular \(M_N\to\infty\), since its regular-variation index is
\(\alpha>0\), even along subsequences where \(I_N\to0\).
We may use a smooth decreasing equivalent of \(u_n\).
Indeed the supremum of the difference between the two paths, both
normalized by \(M_N\), is at most the number of visits to the symmetric
differences of their targets up to \(\lfloor TN\rfloor\), divided by \(M_N\).
Under \(\overline\Lambda\), this number has expectation \(o(M_N)\), by
bracketing the two radii with factors \(1\pm\varepsilon\), regular
variation of the target mass, and then \(\varepsilon\downarrow0\).
The two mean and intensity ratios tend to one. Absolute continuity
transfers this negligible error to each fixed initial law. Finite
prefix changes are harmless; take all comparison radii below \(1/2\).

\emph{Return-horizon transforms.}
Let \(W_0\) have tail \(\mathbb P(W_0>w)=(1+w)^{-r}\).
The profile and integrable mean calculation in the proof of
Lemma~\ref{lem:neutral-critical-rate} give
\(u^{-1}\nu(e^{-qV_u/a_u}-1)\to k_B(\mathbb Ee^{-qD_*W_0}-1)\),
uniformly for \(q\) in compact subsets of \([0,\infty)\).
Use \eqref{eq:neutral-weighted-perturbation} and the separating contours there,
not that lemma's critical-schedule conclusion. With \(\rho_0<1\),
\begin{equation*}
 \eta_{u,q}=1+\nu(e^{-qV_u/a_u}-1)+O(u^2),\quad
 \Pi_{u,q}=\Pi+O(u),\quad
 \mathcal A_{u,q}^K1
 =\eta_{u,q}^K\{\Pi_{u,q}1+O((\rho_0/\eta_{u,q})^K)\}.
\end{equation*}
Here \(\Pi f=\nu(f)1\), and the bounds are in base BV.
If \(K\to\infty\), \(u\to0\) and \(Ku\to a\in(0,\infty)\), then
\begin{equation}\label{eq:occupation-frozen}
 \log(\mathcal A_{u,q}^K1)
 \to ak_B(\mathbb Ee^{-qD_*W_0}-1)
 \quad\text{uniformly outside a null set in }Y.
\end{equation}
The projector error is \(O(u)\), \(Ku^2\to0\), and the relative
remainder vanishes. Divergence of \(Ku\) is not required.

Suppose \(I_N\to\iota\in(0,\infty)\) along the chosen subsequence.
Set \(n_N(0)=0\), \(n_N(t)=\lfloor Nt/\gamma\rfloor+1\) for \(t>0\),
\(v_j=u_{\lfloor\gamma j\rfloor+1}\), and
\(S_N(t)=\sum_{j<n_N(t)}V_{v_j}\circ\mathcal F^j\).
Fix \(0=t_0<\cdots<t_m=T\) and nonnegative increment tilts \(\theta_i\).
On a finite partition of \([\delta,T]\) containing the positive \(t_i\),
freeze the radii at block endpoints with the single denominator \(a_N\).
At fraction \(s\), radius divided by \(u_N\) tends to \(1/s\),
residence scale divided by \(a_N\) to \(s^\alpha\), and block length
times radius to \(\iota\Delta s/(\gamma s)\).
Apply \eqref{eq:occupation-frozen} to the two bounding powers.
The factors are ordered from earliest on the right to latest on the left.
If each block image of \(1\) is bounded by \(L_i1,U_i1\), positivity
bounds the product by \((\prod L_i)1,(\prod U_i)1\); commutation is
unnecessary. Take \(N\) first for the fixed partition, then refine it.
For the initial sector \(j<d=\lfloor\delta N/\gamma\rfloor\),
conditional Jensen and \(\mathcal P^{d-j}V_{v_j}\leq Cm_{v_j}\) bound
its transform below by
\(\exp[-C\max_i\theta_i\sum_{j<d}m_{v_j}/a_N]\) and above by one.
The exponent is \(O(\delta^\alpha)\); all transfer exponents are positive.
Let \(\delta\downarrow0\) after the preceding limits. The joint
increment transform of \(S_N/a_N\) tends to
\begin{equation}\label{eq:occupation-increments}
 \exp\left[\iota\overline c\sum_i\int_{t_{i-1}}^{t_i}
  (\mathbb Ee^{-\theta_iD_*s^\alpha W_0}-1)\frac{ds}{s}\right].
\end{equation}
Since \(M_N/a_N\to\beta_\iota D_* /(r-1)\), the multivariate
Laplace transform identifies the ideal prefix vector with that of
\(Z_\iota\). Independent increments emerge from this transform.

\emph{Physical prefixes.}
For \(0<\Lambda_-<\gamma<\Lambda_+\), put
\(a_+=\lceil\Lambda_+\rceil+2\) and
\begin{equation*}
 S_N^{\lambda,a,b_0}(t)
 =\sum_{j<\max(0,\lfloor Nt/\lambda\rfloor+b_0)}
       V_{u_{\lfloor\lambda j\rfloor+a}}\circ\mathcal F^j.
\end{equation*}
The ordinary return clock gives, for each fixed positive \(t\), eventually
\begin{equation*}
 S_N^{\Lambda_+,a_+,-3}(t)-B_0\leq C_{\lfloor Nt\rfloor}
 \leq S_N^{\Lambda_-,1,1}(t)+B_0,
\end{equation*}
with \(B_0\) pathwise finite. The upper bound follows from
\(\mathsf T_j\geq\Lambda_-j\) and includes the unfinished excursion.
For the lower bound, \(\mathsf T_{j+1}\leq\Lambda_+(j+1)\);
the chosen radius index is after that endpoint, and the \(-3\) offset
keeps every included excursion complete before \(\lfloor Nt\rfloor\).
The finitely many exceptional excursions contribute to \(B_0\).
The ideal \(S_N(t)\) is also between these two deterministic-index sums.
By \(\nu\)-invariance and \(m_u=\gamma\overline\mu([0,u])\),
\begin{equation*}
 \frac{\nu(S_N^{\lambda,a,b_0}(t))}{M_N}\to
 \frac\gamma\lambda t^\alpha.
\end{equation*}
Markov bounds the nonnegative upper-minus-lower error, whose limiting
normalized mean is \(\gamma(\Lambda_-^{-1}-\Lambda_+^{-1})t^\alpha\).
For a finite prefix vector take \(N\), then countable margins to \(\gamma\).
The ordinary clock full event and \(B_0/M_N\to0\) give closeness in
probability of the physical and ideal vectors under \(\nu\).
Only after this comparison take vector differences to obtain increments.

\emph{Initial laws.}
The physical skew map \(\mathcal S\) is ergodic: an invariant event
restricts to a \(\mathcal F\)-invariant event on \(Y\), and every
positive point enters \(Y\) in finite time. Also
\(\nu=\gamma\mathbf1_Y\overline\Lambda\), including the unused future.
Reindexing counts for decreasing targets and telescoping the annulus means
gives \(\mathbb E_{\overline\Lambda}|C_N\circ\mathcal S-C_N|\leq2\).
Every bounded prefix-vector Laplace test \(G_N\) is consequently
asymptotically \(\mathcal S\)-invariant in \(L^1\).
For the physical transfer \(P_{\mathcal S}\), Ces\`aro means
\(H^{-1}\sum_{j<H}P_{\mathcal S}^jf\to\int f\,d\overline\Lambda\) in \(L^1\):
apply Hilbert mean ergodicity to bounded \(L^2\) functions and then truncate.
Indeed \(P_{\mathcal S}f=f\) implies
\(\langle f,f\circ\mathcal S\rangle=\|f\|_2^2\), hence
\(f=f\circ\mathcal S\), and ergodicity makes fixed functions constant.
For fixed \(H\), asymptotic invariance and truncation of a probability
density \(f\) show that replacing \(f\) by this Ces\`aro mean changes
\(\int G_Nf\,d\overline\Lambda\) by \(o(1)\).
Take \(N\) first, then \(H\). Thus the finite-dimensional distributions
converge to the same limits under \(\nu\), \(\overline\Lambda\),
and every fixed \(\rho\ll\overline\Lambda\).

\emph{Tightness and intensity endpoints.}
Under \(\overline\Lambda\), the nonnegative increasing paths
\(X_N(t)=C_{\lfloor Nt\rfloor}/M_N\) satisfy
\(\mathbb EX_N(t)=M_{\lfloor Nt\rfloor}/M_N\to t^\alpha\).
Their interior \(M_1\) oscillation modulus is zero. Terminal means bound
their sup norms, and Markov gives vanishing initial, terminal and
fixed-time nearby increments in the order \(N\to\infty\), then
\(\delta\downarrow0\).
Here is a direct compactness argument. Parametrize the completed graph
of an increasing path by \(t+y\), rescaled to \([0,1]\).
When the terminal value is at most \(K\), both coordinates are Lipschitz
with constant at most \(T+K\). Uniform graph limits exist by
Arzel\`a--Ascoli; bounds \(X_N(\delta_l)\leq\varepsilon_l\downarrow0\)
exclude a spurious initial vertical segment. Terminal vertical segments
are allowed. Choose \(\delta_l\) with summable failure probabilities
from the initial modulus bound, shrinking it below \(1/N\) for the
finitely many exceptional indices. A union bound gives compact sets
carrying arbitrarily high probability for the entire sequence.
For any subsequential limit, choose neighboring times outside its
at-most-countable fixed-discontinuity set. Convergence there and the
nearby-increment bound exclude a jump at each deterministic time,
including \(T\). Dense finite-dimensional laws therefore identify the
path limit. Absolute continuity transfers the compact-set bounds to
each fixed \(\rho\). This proves scalar \(M_1\) convergence for finite \(\iota\).

If \(I_N\to0\), a hit between \(\lfloor\delta N\rfloor\) and
\(\lfloor TN\rfloor\) requires occupation at the first time or a new
right entrance. A left step cannot enter a decreasing target from outside.
Since \(h\leq2\) on \(Y\), its probability is at most
\(\overline\mu([0,u_{\lfloor\delta N\rfloor}])+
\sum_{\delta N<n\leq TN}u_n\to0\). Uniform convergence for slow variation
bounds the sum by \(O_{\delta,T}(I_N)\).
Early normalized occupation has limiting mean \(\delta^\alpha\).
Markov followed by \(\delta\downarrow0\) proves uniform convergence
in probability to zero.

If \(I_N\to\infty\), repeat the finite partition with denominator \(M_N\).
For fixed \(\delta>0\) and a fixed partition of \([\delta,T]\),
each frozen radius \(v\) has \(v/u_N\) and \(a_v/a_N\) bounded above
and away from zero. Since \(M_N/a_N\asymp I_N\), its tilt
\(q=\theta a_v/M_N\) tends to zero.
To justify linearization with this moving tilt, let \(Y_v\) have the
conditional law of \(V_v/a_v\) given \(V_v>0\).
The profile and mean limits in the proof of
Lemma~\ref{lem:neutral-critical-rate}, together with
\(\nu(V_v>0)\sim k_Bv\), give
\(Y_v\Rightarrow D_*W_0\) and
\(\mathbb EY_v\to D_* /(r-1)\). These nonnegative laws are therefore
uniformly integrable as \(v\downarrow0\).
For any \(q_v>0\) tending to zero, split at \(Y_v=K\): on \([0,K]\),
\((1-e^{-q_v y})/q_v\to y\) uniformly, while the tail error is bounded
by \(\mathbb E[Y_v\mathbf1_{Y_v>K}]\). Taking \(v\downarrow0\),
then \(K\to\infty\), proves
\begin{equation*}
 \nu(1-e^{-q_vV_v/a_v})
 =\frac{q_v m_v}{a_v}(1+o(1)).
\end{equation*}
For \(\theta>0\) this yields
\(\nu(e^{-\theta V_v/M_N}-1)=-(\theta/M_N)m_v(1+o(1))\)
at every endpoint of the fixed finite partition; zero tilts give equality.
The total eigenvalue error is \(O(Nu_N^2)=O(I_N^2/N)=o(1)\), since
\(I_N\) is slowly varying. Projector errors and relative remainders vanish.
The ordered transform tends to
\(\exp[-\sum_i\theta_i(t_i^\alpha-t_{i-1}^\alpha)]\).
Early Jensen now costs \(O(\delta^\alpha)\) on the \(M_N\) scale;
the prefix squeeze, change of law and tightness prove the deterministic
limit. Only the intensity along the chosen subsequence enters.

The Poisson integral exists because \(\mathbb EW_0=1/(r-1)\) and
\(\int_0^T s^\alpha ds/s<\infty\). Its small-time mass tends to zero
in \(L^1\); positive time strips have finite jump intensity. Scaling and
disjoint strips give self-similarity and independent increments.
Tonelli and layer cake give \eqref{eq:neutral-occupation-transform};
the raw terminal L\'evy density is \(\beta_\iota(1+z)^{-r}dz/z\),
whose positive mass proves nondegeneracy.
The residence scale \(D_*\) cancels on division by \(M_N\).
Prelimit jumps are at most \(1/M_N\), whereas \(Z_\iota\) has positive
jumps with positive probability; \(J_1\) cannot merge these clusters.
Finally past measurability projects \(m\) to \(\overline\Lambda\).
Lemma~\ref{lem:neutral-normalizer} gives \(Q_N/M_N\to1\) \(m\)-almost surely,
so Slutsky with one scalar multiplier over the path proves the asserted
\(Q_N\) replacement, also under every fixed law absolutely continuous
with respect to \(m\).
\end{proof}

\section{Excursion marks and residence clocks}\label{sec:marks}

The strong laws concern cumulative counts along a fixed divergent schedule.
We now hold the entrance interval fixed at size \(\varepsilon\) and
observe it for order \(\varepsilon^{-1}\) steps. Its rare entrances
generate long neutral excursions. The following result identifies both the
arrival process and the residence carried by each arrival.
Unlike Theorem~\ref{thm:neutral-occupation}, which uses moving targets
and the single-horizon denominator \(M_N\) or \(Q_N\) under annealed
laws, this section uses a fixed target, physical time of order
\(\varepsilon^{-1}\), and a separate normalization for residence times,
with the environment fixed.

For deterministic LSV maps, Freitas, Freitas, Todd and Vaienti
\cite{FFTV16} showed that the usual neutral-point normalization gives
degenerate limits and obtained nondegenerate limits with an adapted
normalization in a restricted parameter range. Zweim\"uller \cite{Zwe19}
proved the exponential first-hitting law at scale
\(h_\alpha(1/2)\varepsilon/2\) for every \(\alpha\in(0,1)\), where
\(h_\alpha\) is the invariant probability density of \(T_\alpha\).
The result below concerns the joint return and residence process under
ergodic forcing, rather than only the first hitting time.

Let \(\sigma\) be an invertible ergodic probability-preserving map
of \((\Omega,\mathbb P)\), let \(\lambda:\Omega\to[\alpha,\beta]\)
be measurable, where \(0<\alpha\leq\beta<1\). Write
\(T_\omega=T_{\lambda(\omega)}\) and
\(T_\omega^n=T_{\sigma^{n-1}\omega}\circ\cdots\circ T_\omega\).
The cone estimates in Section~\ref{sec:sequential} and the coupling in
Lemma~\ref{lem:coupling} give a canonical
equivariant cone probability \(\mu_\omega=h_\omega\,dx\).
Indeed, the pullbacks
\begin{equation*}
 h_\omega^{(m)}=
 \mathcal L_{\lambda(\sigma^{-1}\omega)}\cdots
 \mathcal L_{\lambda(\sigma^{-m}\omega)}\mathbf1
\end{equation*}
are \(L^1\)-Cauchy: comparing two past lengths uses their common last
\(m\) maps and two normalized cone inputs, whose difference after those
maps is at most \(C(1+m)^{-(1/\beta-1)}\). Local equicontinuity and
the common integrable envelope \(Ax^{-\beta}\) preserve normalization
and the cone. The limit is measurable and equivariant. The same comparison
gives uniqueness among normalized equivariant cone families. This construction
uses invertibility but not independence or mixing of the base.

Put
\begin{equation*}
 c(\omega)=\tfrac12h_\omega(1/2+),\qquad
 \bar c=\mathbb E c>0,\qquad p=\mathbb P(\lambda=\alpha).
\end{equation*}
When \(p>0\), put \(D_p=(p\alpha2^\alpha)^{-1}\).
The cone estimates give uniform positive bounds on \(c\), and uniform
Lipschitz bounds on the fibre densities restricted to \(Y\).
At an entrance \(y\in B_\varepsilon=[1/2,(1+\varepsilon)/2]\)
at time \(n\), define
\begin{equation*}
 \kappa_n(y)=\inf\{j\geq1:T_{\sigma^n\omega}^{j}y\in Y\},
 \qquad
 \zeta_{n,\varepsilon}(y)=
 \inf\{k\geq0:T_{\sigma^{n+1}\omega}^{k}(2y-1)>\varepsilon\}.
\end{equation*}
The full-return clock includes the right-branch step. The residence clock
starts after it; for almost every entrance the neutral visits occur at
times \(n+1,\ldots,n+\zeta_{n,\varepsilon}\).

\begin{theorem}[Joint excursion and residence process]\label{thm:marks}
If \(p>0\), then for almost every fixed
\(\omega\) and every fixed probability \(\rho\ll\operatorname{Leb}\),
with \(x\) distributed according to \(\rho\),
\begin{equation}
\label{eq:h-process}
 \sum_{n\geq0}\mathbf1_{B_\varepsilon}(T_\omega^n x)
 \delta_{(\bar c\varepsilon n,
            \varepsilon^\alpha\kappa_n(T_\omega^n x),
            \varepsilon^\alpha\zeta_{n,\varepsilon}(T_\omega^n x))}
 \ \Longrightarrow\ \operatorname{PRM}(dt\,\nu_p(dz,dw)),
\end{equation}
where convergence is in distribution for locally finite random measures
on \([0,\infty)\times[0,\infty)^2\) with the vague topology.
Here \(\operatorname{PRM}\) denotes a Poisson random measure with
the displayed intensity. Time zero and zero marks are included in the
state space; the limiting intensity has no atom there.
The measure \(\nu_p\) is the law of
\begin{equation*}
 (Z,W)=(D_p U^{-\alpha},D_p(U^{-\alpha}-1)),
 \qquad U\sim\operatorname{Uniform}(0,1).
\end{equation*}
The same limit holds with the reference time coordinate
\(\sum_{j<n}\mu_{\sigma^j\omega}(B_\varepsilon)\).
For every fixed \(T>0\),
\begin{equation}
\label{eq:h-occupation}
 V_\varepsilon(t)=\varepsilon^\alpha
  \sum_{0\leq j<\lfloor t/(\bar c\varepsilon)\rfloor}
       \mathbf1_{[0,\varepsilon]}(T_\omega^j x)
 \ \Longrightarrow\ V(t)=\sum_{S_i\leq t}W_i
 \quad\text{in }(D([0,T]),M_1),
\end{equation}
where the \(S_i\) form a rate-one Poisson process and its iid residence
marks have tail \(\mathbb P(W>w)=(1+w/D_p)^{-1/\alpha}\).
The environment may be nonmixing, and its parameter values may accumulate
at \(\alpha\); no positive gap above the minimum is required.

Alternatively, let \(\pi=\lambda_*\mathbb P\),
\(F(v)=\pi([\alpha,v])\), and suppose
\(F(\alpha+t)\sim C_0t^{\kappa_0}L_0(1/t)\),
where \(C_0>0\), \(\kappa_0\geq0\), and \(L_0\) is positive,
measurable and slowly varying. Assume there is a fixed
\(\delta\in(0,1)\) such that, for almost every \(\omega\) and every
integer \(d\geq1\), finite constants \(C_d(\omega),m_d(\omega)\) satisfy
\begin{equation}\label{eq:mark-discrepancy}
 \sup_{\substack{k\in\mathbb Z,\ |k|\leq m^d\\v\in[\alpha,\beta]}}
 \left|\sum_{i=0}^{m-1}
 \mathbf1_{\{\lambda(\sigma^{k+i}\omega)\leq v\}}-mF(v)\right|
 \leq C_d(\omega)m^{1-\delta},\qquad m\geq m_d(\omega).
\end{equation}
No moment assumption on these constants is imposed. Set
\begin{equation*}
 h_\varepsilon=
 \frac{\varepsilon^{-\alpha}(\log(1/\varepsilon))^{\kappa_0}}
 {\alpha2^\alpha C_0\Gamma(\kappa_0+1)L_0(\log(1/\varepsilon))}.
\end{equation*}
Then the same conclusions hold with \(\varepsilon^\alpha\) replaced
by \(h_\varepsilon^{-1}\) in both processes, and with \(\nu_p\)
replaced by the law of
\((U^{-\alpha},U^{-\alpha}-1)\). The limiting residence tail is
\((1+w)^{-1/\alpha}\). This alternative permits atomless minimum
laws, including index zero. The conclusions still hold for every fixed
absolutely continuous initial distribution and with the reference time
coordinate defined above. Every iid parameter
law satisfies \eqref{eq:mark-discrepancy}; dependent nonmixing examples
are given below.

\end{theorem}
\subsection{A complete-history Poisson lemma}

General stationary and nonstationary marked rare-event criteria are
established in \cite{FFM18,FFMV20}.
Joint hitting-time and position limits based on short-time regularization,
without classical mixing assumptions, are developed in \cite{Zwe22}.
Here full-return regularity gives the following physical-time history
argument directly. It will identify the joint physical return/residence
kernel once the intermittent clock is determined.

\begin{lemma}[Complete-history Poisson limit]\label{lem:history}
Fix a parameter sequence \((s_n)_{n\geq0}\) in
\([\alpha,\beta]\subset(0,1)\) and cone probabilities \(\mu_n\)
with \((T_{s_n})_*\mu_n=\mu_{n+1}\). Use the composition
\(T_{n,h}\) defined in the sequential section.
Suppose \(A_{n,\varepsilon}\subset Y\) are unions of full
first-return atoms, have length at most \(C\varepsilon\), and all selected
atoms return within \(s_\varepsilon\) steps, where
\(\varepsilon s_\varepsilon\to0\). Give each atom a constant mark via
\(\xi_{n,\varepsilon}:A_{n,\varepsilon}\to E\), where \(E\) is a compact
metric space. Let increasing deterministic timestamps
\(t_{n,\varepsilon}\) have at most \(C_T/\varepsilon\) observations on
\([0,T]\). Define the deterministic intensity measure by
\begin{equation*}
 \mathfrak I_\varepsilon(I\times G)=
 \sum_{n:t_{n,\varepsilon}\in I}
 \mu_n(A_{n,\varepsilon}\cap\{\xi_{n,\varepsilon}\in G\}).
\end{equation*}
Assume \(\mathfrak I_\varepsilon\) converges on time/mark
continuity rectangles to \(dt\,\nu\), for a finite mark measure
\(\nu\). With \(X_n=T_{0,n}x\) and \(x\) distributed according
to \(\mu_0\), the random measure
\begin{equation*}
 \sum_{n\geq0}\mathbf1_{A_{n,\varepsilon}}(X_n)
   \delta_{(t_{n,\varepsilon},\xi_{n,\varepsilon}(X_n))}
\end{equation*}
converges to that Poisson random measure. Physical gaps in the proof
refer to the integer observation indices, not the rescaled timestamps.

\end{lemma}
\begin{proof}
Use the uniform full-return regularity, its preservation under positive
mixtures, and the coupling in Section~\ref{sec:sequential}, in particular
Lemma~\ref{lem:coupling}. These
hold with one fixed log-Lipschitz constant \(K\), and loss
\(e(h)\to0\). The interpolation argument gives, for every regular
base probability \(\eta\) at time \(n\) and \(J\subset Y\),
\begin{equation}
\label{eq:h-local}
 |\eta(T_{n,h}^{-1}J)-\mu_{n+h}(J)|
       \leq C\operatorname{Leb}(J)\sqrt{e(h)}.
\end{equation}
Consider a history of positive hits. After its last hit, split the
conditional law into its full-atom components. Before a component returns
it gives no mass to another base hit; at return its normalized density is
regular, and afterwards its base restriction remains regular or zero.
At the next hit restrict each positive base component to the new full
atoms, normalize each nonzero piece, and reweight by its mass. Full-return
regularity is scale invariant. Therefore neither the regularity constant
nor the base sup-norm bound increases with the history length. Zero-mass
histories can be omitted. Only positive-hit histories are conditioned on;
factorial moments do not prescribe intermediate non-hits.

For any \(k\) distinct ordered observation times, let \(E_i\) denote
the corresponding marked hit events in the initial probability space.
This induction yields
\begin{equation*}
 \mu_0\Big(\bigcap_{i=1}^kE_i\Big)\leq(C\varepsilon)^k.
\end{equation*}
If successive physical gaps are at least \(s_\varepsilon+g_\varepsilon\),
every last-hit component has returned and mixed for at least
\(g_\varepsilon\) steps. By \eqref{eq:h-local} the next conditional probability
differs from its marginal \(p_i\) by at most
\(C\varepsilon\sqrt{e(g_\varepsilon)}\). Iteration gives
\begin{equation*}
 \left|\mu_0\Big(\bigcap_{i=1}^kE_i\Big)-\prod_i p_i\right|
 \leq k(C\varepsilon)^k\sqrt{e(g_\varepsilon)}.
\end{equation*}
Take \(g_\varepsilon\to\infty\) with \(\varepsilon g_\varepsilon\to0\).
Among \(H_\varepsilon\leq C_T/\varepsilon\) observations, the number of
\(k\)-tuples with a short consecutive gap is at most
\(C_kH_\varepsilon^{k-1}(s_\varepsilon+g_\varepsilon)\). Their total contribution,
and that of the corresponding marginal products, is
\(O_k(\varepsilon(s_\varepsilon+g_\varepsilon))=o(1)\). Summed factorization
errors are \(O_k(\sqrt{e(g_\varepsilon)})\); repeated-time diagonals cost
\(O_k(\varepsilon)\). Expanding mixed factorial moments of counts on
disjoint rectangles now gives the products of the limiting intensities.
Moreover the total count \(N_\varepsilon\) on a bounded window satisfies
\(\mathbb E(N_\varepsilon)_k\leq C_T^k\), so its nonnegative factorial
expansion gives \(\mathbb E(1+z)^{N_\varepsilon}\leq\exp(C_Tz)\)
for \(z\geq0\). These bounds ensure uniform integrability and moment
determinacy. Mixed moments identify independent Poisson counts, and the
first-moment bound gives local tightness. Approximation of compactly
supported continuous functions by continuity rectangles identifies the
limiting random measure.
\end{proof}

\Needspace{5\baselineskip}
\subsection{Ergodic inverse clocks}

\begin{proof}[Proof of Theorem~\ref{thm:marks}]
We first assume \(p>0\) and take \(\rho=\mu_\omega\). The discrepancy
alternative and general fixed absolutely continuous initial laws are
treated below. It remains to verify the intensities and
then transfer residence marks to the occupation process. Put
\(f=\mathbf1_{\{\lambda=\alpha\}}\), and, for rational \(\varsigma>0\),
\(g_\varsigma=\mathbf1_{\{\alpha<\lambda\leq\alpha+\varsigma\}}\),
\(v_\varsigma=\mathbb E g_\varsigma\downarrow0\) as \(\varsigma\downarrow0\).
For rational \(0<\eta<p/2\) and integer \(K\), let
\(D_{K,\eta,\varsigma}\) be the set of phases for which some backward
prefix of length \(m\geq K\) has either \(f\)-frequency differing
from \(p\) by more than \(\eta\), or \(g_\varsigma\)-frequency differing
from \(v_\varsigma\) by more than \(\eta\). Backward Birkhoff gives
\(\mathbb P(D_{K,\eta,\varsigma})\to0\) as \(K\to\infty\), with
\(\eta,\varsigma\) fixed. Take one full-measure set on which forward Birkhoff holds
for \(c\) and all these countably many bad-set indicators.
For every fixed such set \(D\) and nonnegative integer shift \(j_N=o(N)\),
\begin{equation}
\label{eq:h-shift}
 \left|\frac1N\sum_{n<N}\mathbf1_D(\sigma^{n+j_N}\omega)
       -\frac1N\sum_{n<N}\mathbf1_D(\sigma^n\omega)\right|
 \leq\frac{2j_N}{N}.
\end{equation}

Define inverse endpoints
\begin{equation*}
 b_j(n)=\phi_{\lambda(\sigma^n\omega)}\circ\cdots\circ
         \phi_{\lambda(\sigma^{n+j-1}\omega)}(1/2).
\end{equation*}
Their actual inverse application order starts at the terminal phase
\(\sigma^{n+j}\omega\) and reads backwards. At a good terminal
phase, writing \(q_i\) for successive inverse coordinates, the fast
count after \(i\geq K\) steps is between \((p-\eta)i\) and
\((p+\eta)i\). Taylor expansion gives uniformly
\begin{equation*}
 q_i^{-\alpha}-q_{i-1}^{-\alpha}=
 \begin{cases}
 \alpha2^\alpha+O(q_i^\alpha),&\lambda=\alpha,\\
 O(q_i^{\lambda-\alpha}),&\lambda>\alpha.
 \end{cases}
\end{equation*}
All increments are bounded and fast ones have a positive uniform lower
bound. Hence \(q_i\leq Ci^{-1/\alpha}\) for \(i\geq K\).
The near-minimum slow symbols contribute at most
\(C(v_\varsigma+\eta)j\); the remaining slow symbols contribute at most
\(C\sum_{i\leq j}i^{-\varsigma/\alpha}\). The fast error is
\(O(\log(1+j))\), and the initial prefix costs \(O(K)\). Thus
\begin{equation}
\label{eq:h-inverse}
 |b_j(n)^{-\alpha}-p\alpha2^\alpha j|
 \leq C(\eta+v_\varsigma)j+CK+C\log(1+j)
                    +C\sum_{i\leq j}i^{-\varsigma/\alpha}.
\end{equation}
The last sum is \(o(j)\) for each fixed \(\varsigma>0\). Constants multiplying
\(\eta+v_\varsigma\) can be chosen uniformly as \(\varsigma\downarrow0\);
the remainder convergence need not be uniform in \(\varsigma\).

Fix \(t,T>0\), put \(N_\varepsilon(T)=\lfloor T/(\bar c
\varepsilon)\rfloor\), \(\ell_\varepsilon=\lfloor t
\varepsilon^{-\alpha}\rfloor\), and define the clipped threshold
\begin{equation*}
 U_{n,\varepsilon}(t)
   =\min\{1,b_{\ell_\varepsilon-1}(n+1)/\varepsilon\}.
\end{equation*}
Its word ends at \(\sigma^{n+\ell_\varepsilon}\omega\), whose
shift is \(o(N_\varepsilon(T))\) since \(\alpha<1\). On good
phases \eqref{eq:h-inverse} gives discrepancy from
\(u(t)=\min\{1,(D_p/t)^{1/\alpha}\}\) at most
\(C_t(\eta+v_\varsigma)+o_{K,\eta,\varsigma,t}(1)\). On bad phases it is at most
one. Equation \eqref{eq:h-shift} therefore gives
\begin{equation}
\label{eq:h-clipped}
 \limsup_{\varepsilon\downarrow0}\frac1{N_\varepsilon(T)}
 \sum_{n<N_\varepsilon(T)}|U_{n,\varepsilon}(t)-u(t)|
 \leq\mathbb P(D_{K,\eta,\varsigma})+C_t(\eta+v_\varsigma).
\end{equation}
First let \(\varepsilon\to0\), then \(K\to\infty\),
\(\eta\to0\), and finally \(\varsigma\to0\). The average discrepancy
is zero. The clipping is needed because unclipped inverse endpoints
need not converge in mean.

Uniform base Lipschitz regularity gives
\begin{equation*}
 \mu_{\sigma^n\omega}(B_\varepsilon\cap\{
                   \varepsilon^\alpha\kappa_n>t\})
 =\varepsilon c(\sigma^n\omega)U_{n,\varepsilon}(t)
       +O(\varepsilon^2).
\end{equation*}
Use boundedness of \(c\), \eqref{eq:h-clipped}, and Birkhoff for \(c\):
the sum up to \(N_\varepsilon(T)\) tends to \(Tu(t)\). This
does not assume independence of \(c\) and the parameter word. Rational
\(t,T\), continuity, monotonicity and differences of prefixes give
all time/mark continuity rectangles. The same argument without a survival
threshold gives total intensity \(T\). A finite grid and monotonicity
give uniform convergence of the compensator to \(\bar c\varepsilon n\)
on bounded time windows.

\subsection{Full-return and residence marks}

Fix \(L>D_p\). Replace \(B_\varepsilon\) by its inner full-atom
union and retain returns \(\kappa\leq L\varepsilon^{-\alpha}\).
The single-cut estimate from the sequential proof changes expected counts
by \(O(\varepsilon^\alpha)\) on \(O(\varepsilon^{-1})\) observations;
equivariance and a union bound give the same bound for the probability of
any changed hit. Return marks are constant on full atoms, and
\(\varepsilon s_\varepsilon\leq L\varepsilon^{1-\alpha}\to0\).
The complete-history lemma now yields the truncated full-return process.
The expected number of omitted long marks tends to
\(T(D_p/L)^{1/\alpha}\); Markov's inequality and \(L\to\infty\)
remove the cutoff after \(\varepsilon\to0\).

To identify residence jointly, let \(q_m^\varepsilon(n+1)\) be the
same inverse word started at \(\varepsilon\), and take
\(m=\lfloor w\varepsilon^{-\alpha}\rfloor\), \(w>0\).
All intermediate coordinates are at most \(\varepsilon\). The
increment expansion on a good terminal phase yields
\begin{equation*}
 \varepsilon^\alpha(q_m^\varepsilon)^{-\alpha}
 =1+p\alpha2^\alpha w+O_w(\eta+v_\varsigma)+o_{K,\varsigma,w}(1).
\end{equation*}
The scaled fast error is \(O(m\varepsilon^{2\alpha})\), near-minimum
slow error is \(O(m\varepsilon^\alpha(v_\varsigma+\eta))\), and separated
slow error is \(O(m\varepsilon^{\alpha+\varsigma})\); the fixed initial
prefix and integer rounding vanish. Its terminal shift \(m+1\) is
again \(o(\varepsilon^{-1})\). The ratio
\(q_m^\varepsilon/\varepsilon\) is already bounded by one, so the
same bad-set argument gives averaged convergence to
\((1+w/D_p)^{-1/\alpha}\).

The return and residence survival sets have the same left endpoint zero
after the right branch. Their joint upper endpoint is
\(\min\{\varepsilon,b_{\ell_\varepsilon-1}(n+1),
q_m^\varepsilon(n+1)\}\). Since minimum is Lipschitz in bounded
ratios, their joint cumulative intensity tends to
\begin{equation}
\label{eq:h-joint}
 T\min\{1,(D_p/t)^{1/\alpha},(1+w/D_p)^{-1/\alpha}\}.
\end{equation}
This is precisely the joint tail of the displayed \((Z,W)\).
Total intensity and the two marginal tails give tightness of the finite
joint-intensity measures, and rational upper rectangles identify their
weak limit. The kernel has no atom at \(Z=L\), and the prelimit
inequality \(\zeta\leq\kappa-1\) bounds the residence coordinate
on \(Z\leq L\). For every \(\eta>0\), the limiting measure of
the closed set
\(\{Z\leq L,\ |W-(Z-D_p)_+|\geq\eta\}\) is zero. Portmanteau
therefore gives
\begin{equation*}
 \sum_{n<N_\varepsilon(T)}\mu_{\sigma^n\omega}
 \left(B_\varepsilon\cap\{\varepsilon^\alpha\kappa_n\leq L,
 |\varepsilon^\alpha\zeta_{n,\varepsilon}
       -(\varepsilon^\alpha\kappa_n-D_p)_+|\geq\eta\}\right)
 \longrightarrow0.
\end{equation*}
By equivariance and a union bound, no retained hit has such a discrepancy
with probability tending to one. Tightness of the number of hits then
transfers the residence coordinate from the continuous function
\((z-D_p)_+\) of the full-return mark. This proves the joint process
\eqref{eq:h-process}.

\subsection{Occupation as compressed residence clusters}

The left branch increases the coordinate. Thus every neutral visit belongs
to a right-branch entrance cluster, unless the orbit starts in
\([0,\varepsilon]\), an event of probability
\(O(\varepsilon^{1-\beta})\). For returns bounded by
\(L\varepsilon^{-\alpha}\), each cluster spans at most
\(\bar cL\varepsilon^{1-\alpha}\) in rescaled time. A cluster
crossing a fixed endpoint \(T\) must enter in the preceding physical
window of length \(L\varepsilon^{-\alpha}+1\); uniform base bounds
make the probability of such an entrance
\(O(L\varepsilon^{1-\alpha}+\varepsilon)\).

Collapse each retained complete cluster to one jump at its entrance, of
height \(\varepsilon^\alpha\zeta\). For fixed \(L\), joint marked
convergence gives convergence of this finite-jump process: limit times
are distinct and avoid prescribed endpoints, so a time change matches
the finite sets of times and heights. The actual occupation staircase
and its collapsed vertical segment have the same ordered height ranges.
Parametrize their completed graphs by matching heights, pausing when
necessary; their time coordinates differ by at most the maximal cluster
duration plus a mesh step. Hence their \(M_1\) distance is at most
\(\bar cL\varepsilon^{1-\alpha}+O(\varepsilon)\) outside the
initial/final exceptional events. Remove long marks after taking
\(\varepsilon\to0\), and obtain \eqref{eq:h-occupation}.

\subsection{Discrepancy-controlled inverse clocks}

For the alternative case, write \(b(v)=v^{\kappa_0}/L_0(v)\),
\(a_u=u^{-\alpha}b(\log(1/u))\), and use the constants
\(A=\alpha^{\kappa_0+1}2^\alpha C_0\Gamma(\kappa_0+1)\),
\(D_*=(\alpha2^\alpha C_0\Gamma(\kappa_0+1))^{-1}\).
Extend \(b\) to be bounded above and away from zero on finite argument
intervals. We prove the two word profiles needed here directly, without
using the iid conditional-mean or integrable-envelope assertions of
Lemma~\ref{lem:neutral-normalizer}.

On the full environment set of
\eqref{eq:mark-discrepancy}, for each fixed polynomial phase exponent
\(D\), we prove the compact limit \eqref{eq:mark-compact-profile}
starting at \(u\), uniformly on compact scaled time intervals, and the
deep limit \eqref{eq:mark-deep-profile} whenever its initial-coordinate
condition holds. The phase windows are respectively \(a_u^D\) and \(L^D\).

Read an inverse word in application order as \(v_1,v_2,\ldots\), and
put \(q_i=\phi_{v_i}(q_{i-1})\), \(z_i=q_i^{-\alpha}\).
The deterministic inverse expansion gives
\begin{equation*}
 z_i-z_{i-1}=f_{z_{i-1}}(v_i)+O(z_{i-1}^{-1}),\qquad
 f_z(s)=\alpha2^s z^{-(s-\alpha)/\alpha},\qquad
 0<z_i-z_{i-1}\leq C.
\end{equation*}
For large \(z\), the function \(f_z\) is decreasing in \(s\),
has uniformly bounded variation, and
\(\sup_s|\partial_z f_z(s)|\leq C/(z\log z)\).
To see the last bound, write \(s=\alpha+v\) and use
\(\sup_{v\geq0}ve^{-vt}\leq1/(et)\), with
\(t=\log z/\alpha-\log2\). The spatial average
\begin{equation*}
 M(z)=\int f_z(s)\,d\pi(s)\sim A/b(\log z)
\end{equation*}
follows from the Laplace--Stieltjes calculation in the neutral
normalization proof. It is not a conditional next-symbol mean.

Freeze the coordinate at the start of a block of \(m\) symbols.
Integration by parts against its empirical CDF and
\eqref{eq:mark-discrepancy} bound its frozen error by
\(C C_d(\omega)m^{1-\delta}\). Both measures have mass \(m\),
so there is no endpoint mass error, including for atomic \(\pi\).
If the coordinates are at least \(Z_{\min}\), bounded increments
and the derivative estimate give adaptive error
\(Cm^2/(Z_{\min}\log Z_{\min})\) per block, for both \(f_z\)
and \(M(z)\). Splitting \(N\) steps into complete blocks and one
remainder therefore yields, uniformly for every prefix,
\begin{equation}\label{eq:mark-freezing}
 \begin{split}
 z_N-z_0&=\sum_{i=1}^N M(z_{i-1})+E_N,\\
 |E_N|&\leq C_\omega\left(Nm^{-\delta}
 +\frac{Nm}{Z_{\min}\log Z_{\min}}+\frac{N}{Z_{\min}}+m\right).
 \end{split}
\end{equation}
For fixed window and block exponents (and fixed \(\varpi,\eta\) below),
choose one sufficiently large \(d\). Constants and cutoffs may depend
on these choices, \(\omega\) and \(\pi\), but not on the radius,
word length or terminal phase inside that window. Endpoint parameter
limits can be taken along countable rational choices on the same event.
The third term is the inverse Taylor remainder and the last term is the
incomplete block. Reverse blocks are contiguous blocks of the original
two-sided word. If the total word length is at most \(R\) and its
terminal phase has absolute value at most \(R^D\), every block start
has absolute value at most \(R^D+R\). When \(m\) is a fixed
positive power of \(R\), up to subpower factors, one fixed sufficiently
large integer \(d\) makes \(R^D+R\leq m^d\). Thus the same
\(C_d(\omega)\) controls all blocks; it does not vary with the scale.

For compact residence windows, let \(q_0=u\),
\(z_0=u^{-\alpha}\), \(N\leq T a_u\), and
\(m=\lfloor z_0^\theta\rfloor\), \(0<\theta<1\).
Then \(z_0\leq z_i\leq z_0+CN\), so
\(\log z_i/\log z_0\to1\) uniformly. The uniform-convergence
theorem for regular variation gives \(M(z_i)\sim A/b(\log z_0)\).
With \(Z_{\min}=z_0\), all four errors in
\eqref{eq:mark-freezing}, divided by \(z_0\), tend to zero:
\(b(\log(1/u))\) is subpower in \(z_0\).
Since \(b(\alpha v)/b(v)\to\alpha^{\kappa_0}\),
\begin{equation}\label{eq:mark-compact-profile}
 q_{\lfloor t a_u\rfloor}(u)/u\longrightarrow
       (1+t/D_*)^{-1/\alpha}.
\end{equation}
This is uniform for \(t\) in compact subsets of \([0,\infty)\)
and terminal phases of absolute value at most \(a_u^D\), for each
fixed \(D\). The estimates hold directly for all real radii.

For the deep profile, let \(L\to\infty\),
\(z_0=u^{-\alpha}=o(L/b(\log L))\), and bound terminal phases
by \(L^D\). Fix \(0<\varpi<1\) and \(\eta>0\) with
\(\alpha+\eta<1\), and put \(\chi=\alpha/(\alpha+\eta)\).
The fixed band \([\alpha,\alpha+\eta]\) has mass \(p_\eta>0\).
For each reverse prefix of length \(i\geq L^\varpi\), use
\eqref{eq:mark-discrepancy} with a fixed integer
\(d>\max(D,1)/\varpi\). Its fast-symbol count is at least
\(p_\eta i/2\), eventually and uniformly. Each such symbol
increases \(q^{-(\alpha+\eta)}\) by a fixed positive amount;
the other increments are nonnegative. Hence, on these prefixes,
\begin{equation*}
 c_\eta L^{\varpi\chi}\leq z_i\leq z_0+CL\leq C'L.
\end{equation*}
Discard the first \(\lfloor L^\varpi\rfloor\) steps, costing
\(O(L^\varpi)\), and set \(m=\lfloor L^\theta\rfloor\),
\(0<\theta<\varpi\chi\). With
\(Z_{\min}=c_\eta L^{\varpi\chi}\), each error in
\eqref{eq:mark-freezing}, multiplied by \(b(\log L)/L\), vanishes.
Uniform convergence on the logarithmic ratio interval
\([\varpi\chi+o(1),1+o(1)]\) bounds the scaled mean drift between
\(A(1-o(1))\) and \(A(\varpi\chi)^{-\kappa_0}(1+o(1))\).
The initial coordinate and discarded prefix are negligible. Thus
\begin{equation*}
 A\leq\liminf\frac{z_L b(\log L)}L
 \leq\limsup\frac{z_L b(\log L)}L
 \leq A(\varpi\chi)^{-\kappa_0}.
\end{equation*}
First take \(L\to\infty\), then \(\varpi\uparrow1\) and
\(\eta\downarrow0\). This proves
\begin{equation}\label{eq:mark-deep-profile}
 q_L(u)\sim A^{-1/\alpha}L^{-1/\alpha}b(\log L)^{1/\alpha},
 \qquad u^{-\alpha}=o(L/b(\log L)),
\end{equation}
uniformly in each polynomial terminal-phase window. Index zero is
included; no minimum atom or moment of the discrepancy constants was used.

Apply the compact and deep profiles just proved. With
\(b(v)=v^{\kappa_0}/L_0(v)\) and \(A,D_*\) as above,
\(h_\varepsilon=D_*\varepsilon^{-\alpha}b(\log(1/\varepsilon))\)
and \(AD_*=\alpha^{\kappa_0}\). For each fixed \(t>0\),
\(\log(th_\varepsilon)/\log(1/\varepsilon)\to\alpha\), hence
\(b(\log(th_\varepsilon))/b(\log(1/\varepsilon))
\to\alpha^{\kappa_0}\). The deep profile started at \(1/2\)
therefore gives \(b_{\lfloor th_\varepsilon\rfloor-1}(n+1)/
\varepsilon\to t^{-1/\alpha}\); the compact profile gives
\(q_{\lfloor wh_\varepsilon\rfloor}^\varepsilon(n+1)/\varepsilon
\to(1+w)^{-1/\alpha}\), for \(w\geq0\).
Both limits are uniform over \(n<N_\varepsilon(T)\): choose an
integer \(d>1/\alpha\) in the polynomial terminal-phase windows
of those profiles. Indeed \(h_\varepsilon\) is
\(\varepsilon^{-\alpha}\) times a subpower correction, so these
windows contain every phase \(n+O(h_\varepsilon)\). The single
right-branch step and integer rounding vanish on this scale.

The common landing endpoint and the uniform base trace expansion now give
\begin{equation*}
 \sum_{n<N_\varepsilon(T)}\mu_{\sigma^n\omega}
 \bigl(B_\varepsilon\cap\{\kappa_n/h_\varepsilon>t,
                  \zeta_{n,\varepsilon}/h_\varepsilon>w\}\bigr)
 \longrightarrow T\min\{1,t^{-1/\alpha},(1+w)^{-1/\alpha}\}.
\end{equation*}
The summed trace error is \(O(\varepsilon)\); Birkhoff applies to
the fixed bounded function \(c\), without independence from inverse
words. These joint intensities identify the displayed shared-\(U\) law.
Retain full returns \(\kappa\leq Lh_\varepsilon\), \(L>1\).
The inner-atom error remains \(O(\varepsilon^\alpha)\), and
\(\varepsilon\lceil Lh_\varepsilon\rceil\to0\), so
Lemma~\ref{lem:history} gives the truncated full-return process.
The omitted-return mean tends to \(TL^{-1/\alpha}\).
On retained marks, the joint-intensity limit assigns zero mass to
\(\{Z\leq L,\ |W-(Z-1)_+|\geq\eta\}\), \(\eta>0\).
The preceding Portmanteau and union-bound argument thus transfers actual
residence from the return mark. Cluster compression uses width at most
\(\bar cL\varepsilon h_\varepsilon\to0\), with endpoint-crossing
probability \(O(L\varepsilon h_\varepsilon+\varepsilon)\).
Letting \(L\to\infty\) after \(\varepsilon\to0\) proves the joint
and \(M_1\) occupation limits, including the reference time change.

For a general fixed \(\rho\ll\operatorname{Leb}\), intersect the good
environment set above with all its integer shifts. Approximation of the
density of \(\rho\) in \(L^1\) by normalized nonnegative H\"older densities,
the uniform memory bound in \cite{KL21}, and transfer
contraction give
\begin{equation*}
 d_k=\|(T_\omega^k)_*\rho-\mu_{\sigma^k\omega}\|_{\mathrm{TV}}
 \longrightarrow0\qquad(k\to\infty).
\end{equation*}
Fix \(k\) before sending \(\varepsilon\to0\). Mapping the initial
state to its entire future process does not increase total variation,
so the process started at time \(k\) differs in law by at most \(d_k\)
from the reference process at \(\sigma^k\omega\). The probability of
an entrance during the first \(k\) observations tends to zero, since
each fixed-time pushforward of \(\rho\) is absolutely continuous.
All later marks are unchanged by reindexing. Their timestamps shift by
\(\bar c\varepsilon k\), or by at most \(C\varepsilon k\) for the
reference compensator; these shifts vanish in the vague topology, using
local tightness of the reference counts. For occupation, removing the
first \(k\) and adding the last \(k\) observations changes the process
by at most \(2k\varepsilon^\alpha\), or \(2k/h_\varepsilon\) in
the discrepancy alternative, uniformly on \([0,T]\), and hence
negligibly in \(M_1\). First let \(\varepsilon\to0\), then \(k\to\infty\).
Both limits therefore hold under \(\rho\), on the same full environment set.
\end{proof}

The residence law has finite mean \(D/(1/\alpha-1)\), with
\(D=D_p\) in the ergodic atom normalization and \(D=1\) in the
discrepancy alternative. Its second
moment is finite precisely when \(\alpha<1/2\); this is a statement
about the limiting entrance mark, not convergence of prelimit moments.
The \(J_1\) topology cannot give the occupation limit: every prelimit
jump has height \(\varepsilon^\alpha\), or \(h_\varepsilon^{-1}\)
in the discrepancy alternative, tending to zero, whereas the limit has
positive jumps. The monotone cluster compression is exactly what \(M_1\)
allows. In particular \(\mathbb P(V(t)=0)=e^{-t}\), reflecting the
possibility of no entrance in a rare-event window. This limit, obtained by
shrinking a fixed target while enlarging its observation window, is distinct
from an almost-sure strong law for one fixed target schedule.

For a minimum atom of mass \(p\), take \(\kappa_0=0\),
\(C_0=p\), \(L_0=1\). Then the alternative scale is
\(h_\varepsilon=D_p\varepsilon^{-\alpha}\), recovering the first
normalization after multiplication by \(D_p\). For uniform noise on
\([\alpha,\beta]\), \(\alpha<\beta\), it is
\(h_\varepsilon=(\beta-\alpha)\varepsilon^{-\alpha}
\log(1/\varepsilon)/(\alpha2^\alpha)\).
The corresponding logarithmic inverse-endpoint correction is already
present in \cite{BBR}; the same clock here
normalizes the coupled return and residence processes.

To check \eqref{eq:mark-discrepancy} for iid laws, write the symbols as
generalized quantiles of iid uniforms. For a length-\(m\) uniform block,
the empirical count at \(j/m\) has deviation probability at most
\(2\exp(-2a^2/m)\). Indeed the logarithm of the moment generating
function of a centred Bernoulli variable has
second derivative at most \(1/4\); integration and exponential Markov
give this bound. Take \(a=m^{3/4}/4\), and union over \(m+1\)
grid points and \(|k|\leq m^d\). The resulting bound
\(Cm^{d+1}e^{-c\sqrt m}\) is summable. Monotonicity between grid
points costs at most one count unit. First Borel--Cantelli and the
countable intersection over \(d\) prove the assertion with
\(\delta=1/4\). Quantile threshold sets include atomic marginals.

\begin{example}[Nonmixing continuum of parameters]\label{ex:rotation}
Choose
\(0<\alpha<\beta<1\). On \(\mathbb R/\mathbb Z\)
let \(\sigma u=u+\vartheta\), where \(\vartheta\) is irrational, and set
\(\lambda(u)=\alpha+(\beta-\alpha)\max\{0,\cos(2\pi u)\}\).
Then \(p=1/2\), while the slow parameters take a continuum of values
accumulating at \(\alpha\). The theorem applies although the base is
not mixing: the centred correlation of \(\cos(2\pi u)\)
with its \(n\)-step translate is \(\tfrac12\cos(2\pi n\vartheta)\).
The parameter correlations may affect \(\bar c\),
although the normalized joint mark kernel depends only on \(\alpha,p\).

For an atomless example on the same base, instead take
\(\lambda(u)=\alpha+2(\beta-\alpha)\operatorname{dist}(u,0)\).
Its marginal is uniform on \([\alpha,\beta]\), and every parameter
threshold preimage is a circle arc. For every bounded-type angle, and
for Lebesgue-almost every irrational angle, arc counts have discrepancy
\(O_\vartheta((\log(2+m))^3)\), uniformly in the arc and phase.
Consequently \eqref{eq:mark-discrepancy} holds with \(\delta=1/2\),
and the theorem gives the unit paired-mark law at
\begin{equation*}
 h_\varepsilon=\frac{\beta-\alpha}{\alpha2^\alpha}
       \varepsilon^{-\alpha}\log(1/\varepsilon).
\end{equation*}
These limits hold for almost every fixed phase, with randomness in the
initial fibre point, not in the rotation orbit alone.

The discrepancy bound is the classical continued-fraction estimate,
as recorded in \cite[Section 2, equation (2.4)]{FH24}.
A convergent denominator block is within \(1/q\) of a translated rational
\(q\)-grid, so its arc-count error is bounded uniformly. Decomposing
\(m\) into denominator blocks bounds the total error by a constant times
the sum of its first \(K+1\) partial quotients, with \(K=O(\log m)\).
Bounded partial quotients give \(O(\log m)\). For the metric assertion,
invariance of the Gauss measure, which is equivalent to Lebesgue measure,
shows that the \(j\)-th partial quotient exceeds \(j^2\) with probability
\(O(j^{-2})\). The first Borel--Cantelli lemma bounds the partial
quotients by \(C(\vartheta)j^2\), and the Fibonacci lower bound for the
denominators gives the stated log-cubed bound.
\end{example}

The first example does not satisfy the full-skew mixing hypothesis used
in \cite{AHV}. Random intermittent Poisson limits in \cite{ANT} are
important antecedents; they concern time-collapsed counts and
heavy-observable marks at generic centres, whereas here the
expanding-branch target produces a full physical-time process with
coupled return and residence marks. Quenched compensator-time Poisson
limits without base mixing were proved earlier in \cite{CS} for
cylinders in random subshifts of finite type.
\section*{Concluding remarks}

The common thread of this paper is the difference between how often an
orbit enters a target and how long it then stays there. Away from the
neutral fixed point, correlations between visits are controlled by the
masses of the two targets involved, and this gives a strong law under
divergence of the total length alone, however slow. Near the fixed
point, a single entrance may contribute a substantial part of the count.
For regularly varying schedules this produces three separate thresholds:
the integral test for the strong law, the weaker condition
\(nu_n\to\infty\) for convergence in mean, and the scale
\(nu_n\sim c\log\log n\), \(c>0\), at which the normalized count
has an almost-sure lower limit strictly between zero and one.
Below this scale the lower limit is zero; above it the lower limit is
one. When \(nu_n\) has a finite positive limit, the
clusters remain visible after normalization, and the occupation paths
converge to an additive jump process. For a point-mass parameter law,
all of these statements concern a single intermittent map.

The same distinction organizes the remaining results. For irregular
neutral schedules, the entrance envelope of Theorem~\ref{thm:irregular},
rather than the total length, decides whether there are infinitely many
visits, and the observation-gap test shows that the support of the
parameter law and its minimum do not determine recurrence by themselves.
This is a hitting criterion, not a strong law. Under ergodic forcing,
Theorem~\ref{thm:marks} describes the return and residence marks of
excursions to a fixed target, an observation scheme different from the
moving targets above.

Several questions remain open. We do not know a necessary and sufficient
condition for the strong law along neutral schedules that are not
regularly varying; the entrance envelope only decides whether infinitely
many visits occur. The moving-target occupation limits of
Theorem~\ref{thm:neutral-occupation} are annealed, and we do not know
whether they hold for almost every fixed environment. Finally, the
neutral-count results use independent parameters, whereas the
excursion-mark limits allow ergodic forcing; whether the integral
criterion and the critical lower envelope persist under ergodic forcing
is also open.

\section*{Acknowledgements}

The author acknowledges Gewu Intelligence Lab for providing the collaborative
research environment in which this project was developed.

\section*{Competing interests}

The author declares none.

\section*{Generative AI declaration}

\href{https://openai.com/}{OpenAI GPT-5.6 Sol},
\href{https://www.anthropic.com/}{Anthropic Claude Fable 5.1},
OpenAI GPT-6.1 Sol, OpenAI GPT-6 Astra and Anthropic Claude Opus 5.5
contributed to the development and preparation of this paper during
September--October 2026.
Their combined work included developing and checking the inducing and
target-mass estimates, the strong-law and lower-envelope arguments,
the neutral-target criterion and second-moment obstruction, and the occupation and
excursion-mark limits.
They also searched the literature, checked references, and drafted and
revised the manuscript.


{\small
\begin{thebibliography}{99}
\setlength{\itemsep}{1pt}
\bibitem{ANT}
R. Aimino, M. Nicol and A. T\"or\"ok. \href{https://doi.org/10.1017/etds.2024.5}{Stable laws for random dynamical systems}. \emph{Ergod. Th. \& Dynam. Sys.} 44(11) (2024), 3041--3090.

\bibitem{ANV}
R. Aimino, M. Nicol and S. Vaienti. \href{https://doi.org/10.1007/s00440-014-0571-y}{Annealed and quenched limit theorems for random expanding dynamical systems}. \emph{Probab. Theory Related Fields} 162(1--2) (2015), 233--274.

\bibitem{AHV}
L. Amorim, N. Haydn and S. Vaienti. \href{https://doi.org/10.1016/j.spa.2024.104511}{Compound Poisson distributions for random dynamical systems using probabilistic approximations}. \emph{Stochastic Process. Appl.} 179 (2025), Paper No. 104511.

\bibitem{BB16}
W. Bahsoun and C. Bose. \href{https://doi.org/10.1088/0951-7715/29/4/1417}{Mixing rates and limit theorems for random intermittent maps}. \emph{Nonlinearity} 29(4) (2016), 1417--1433; \href{https://doi.org/10.1088/0951-7715/29/12/c4}{Corrigendum}, \emph{Nonlinearity} 29(12) (2016), C4.

\bibitem{BBD14}
W. Bahsoun, C. Bose and Y. Duan. \href{https://doi.org/10.1088/0951-7715/27/7/1543}{Decay of correlation for random intermittent maps}. \emph{Nonlinearity} 27(7) (2014), 1543--1554.

\bibitem{BBR}
W. Bahsoun, C. Bose and M. Ruziboev. \href{https://doi.org/10.1090/tran/7811}{Quenched decay of correlations for slowly mixing systems}. \emph{Trans. Amer. Math. Soc.} 372(9) (2019), 6547--6587.

\bibitem{BGT}
N. H. Bingham, C. M. Goldie and J. L. Teugels. \emph{\href{https://doi.org/10.1017/CBO9780511721434}{Regular Variation}} (Encyclopedia of Mathematics and its Applications, 27). Cambridge University Press, Cambridge, 1987.

\bibitem{BT26}
C. Bonanno and S. N. Tikekar. \href{https://doi.org/10.1090/proc/17925}{On the escape rate for intermittent maps with holes shrinking around the indifferent fixed point}. \emph{Proc. Amer. Math. Soc.}, to appear; published online 2026; \href{https://arxiv.org/abs/2601.15908v3}{arXiv:2601.15908v3}.

\bibitem{CN}
M. Carney and M. Nicol. \href{https://doi.org/10.1088/1361-6544/aa72c2}{Dynamical Borel--Cantelli lemmas and rates of growth of Birkhoff sums of non-integrable observables on chaotic dynamical systems}. \emph{Nonlinearity} 30(7) (2017), 2854--2870.

\bibitem{CS}
H. Crimmins and B. Saussol. \href{https://doi.org/10.1088/1361-6544/ac3d00}{Quenched Poisson processes for random subshifts of finite type}. \emph{Nonlinearity} 35(6) (2022), 3036--3058.

\bibitem{DMR}
J. Dedecker, F. Merlev\`ede and E. Rio. \href{https://doi.org/10.1007/978-3-030-96409-2_7}{Criteria for Borel--Cantelli lemmas with applications to Markov chains and dynamical systems}. \emph{S\'eminaire de Probabilit\'es LI} (Lecture Notes in Mathematics, 2301). Eds. C. Donati-Martin, A. Lejay and A. Rouault. Springer, Cham, 2022, pp. 189--238.

\bibitem{DHL}
D. Dragi\v{c}evi\'c, Y. Hafouta and J. Lepp\"anen. \href{https://doi.org/10.1112/jlms.70707}{Decay of correlations and limit theorems for random intermittent maps}. \emph{J. Lond. Math. Soc. (2)} 114(3) (2026), Paper No. e70707.

\bibitem{FL13}
C. Feng and L. Li. \href{https://doi.org/10.1556/SScMath.50.2013.2.1241}{On the M\'ori--Sz\'ekely conjectures for the Borel--Cantelli lemma}. \emph{Studia Sci. Math. Hungar.} 50(2) (2013), 280--285.

\bibitem{FLS09}
C. Feng, L. Li and J. Shen. \href{https://doi.org/10.1016/j.crma.2009.09.011}{On the Borel--Cantelli lemma and its generalization}. \emph{C. R. Math. Acad. Sci. Paris} 347(21--22) (2009), 1313--1316.

\bibitem{FFM18}
A. C. M. Freitas, J. M. Freitas and M. Magalh\~aes. \href{https://doi.org/10.4171/JEMS/808}{Convergence of marked point processes of excesses for dynamical systems}. \emph{J. Eur. Math. Soc. (JEMS)} 20(9) (2018), 2131--2179.

\bibitem{FFMV20}
A. C. M. Freitas, J. M. Freitas, M. Magalh\~aes and S. Vaienti. \href{https://doi.org/10.1007/s10955-020-02630-z}{Point processes of non stationary sequences generated by sequential and random dynamical systems}. \emph{J. Stat. Phys.} 181(4) (2020), 1365--1409.

\bibitem{FFTV16}
A. C. M. Freitas, J. M. Freitas, M. Todd and S. Vaienti. \href{https://doi.org/10.1016/j.spa.2016.05.001}{Rare events for the Manneville--Pomeau map}. \emph{Stochastic Process. Appl.} 126(11) (2016), 3463--3479.

\bibitem{FFV18}
A. C. M. Freitas, J. M. Freitas and S. Vaienti. \href{https://doi.org/10.1090/proc/13892}{Extreme value laws for sequences of intermittent maps}. \emph{Proc. Amer. Math. Soc.} 146(5) (2018), 2103--2116.

\bibitem{Fro21}
A. N. Frolov. \href{https://doi.org/10.1016/j.jmaa.2021.125425}{On strong forms of the Borel--Cantelli lemma and intermittent interval maps}. \emph{J. Math. Anal. Appl.} 504(2) (2021), Paper No. 125425.

\bibitem{Fro25}
A. N. Frolov. \href{https://doi.org/10.1134/S1063454125700517}{On the Borel--Cantelli lemma for intermittent interval maps}. \emph{Vestnik St. Petersburg Univ. Math.} 58(4) (2025), 504--510.

\bibitem{FH24}
L. Fr\"uhwirth and M. Hauke. \href{https://doi.org/10.5802/ahl.200}{On Birkhoff sums that satisfy no temporal distributional limit theorem for almost every irrational}. \emph{Ann. H. Lebesgue} 7 (2024), 251--265.

\bibitem{Gou04}
S. Gou\"ezel. \href{https://doi.org/10.1007/BF02787541}{Sharp polynomial estimates for the decay of correlations}. \emph{Israel J. Math.} 139 (2004), 29--65.

\bibitem{Gou07}
S. Gou\"ezel. \href{https://doi.org/10.1088/0951-7715/20/6/010}{A Borel--Cantelli lemma for intermittent interval maps}. \emph{Nonlinearity} 20(6) (2007), 1491--1497.

\bibitem{GNO}
C. Gupta, M. Nicol and W. Ott. \href{https://doi.org/10.1088/0951-7715/23/8/010}{A Borel--Cantelli lemma for nonuniformly expanding dynamical systems}. \emph{Nonlinearity} 23(8) (2010), 1991--2008.

\bibitem{HNPV13}
N. Haydn, M. Nicol, T. Persson and S. Vaienti. \href{https://doi.org/10.1017/S014338571100099X}{A note on Borel--Cantelli lemmas for non-uniformly hyperbolic dynamical systems}. \emph{Ergod. Th. \& Dynam. Sys.} 33(2) (2013), 475--498.

\bibitem{HKKP24}
M. Holland, M. Kirsebom, P. Kunde and T. Persson. \href{https://doi.org/10.1090/tran/9102}{Dichotomy results for eventually always hitting time statistics and almost sure growth of extremes}. \emph{Trans. Amer. Math. Soc.} 377(6) (2024), 3927--3982.

\bibitem{HT25}
M. Holland and M. Todd. \href{https://doi.org/10.1088/1361-6544/ade044}{On distributional limit laws for recurrence}. \emph{Nonlinearity} 38(7) (2025), Paper No. 075028.

\bibitem{KL09}
G. Keller and C. Liverani. \href{https://doi.org/10.1007/s10955-009-9747-8}{Rare events, escape rates and quasistationarity: some exact formulae}. \emph{J. Stat. Phys.} 135(3) (2009), 519--534.

\bibitem{Kim07}
D. H. Kim. \href{https://doi.org/10.3934/dcds.2007.17.891}{The dynamical Borel--Cantelli lemma for interval maps}. \emph{Discrete Contin. Dyn. Syst.} 17(4) (2007), 891--900.

\bibitem{KKP20}
M. Kirsebom, P. Kunde and T. Persson. \href{https://doi.org/10.1088/1361-6544/ab5160}{Shrinking targets and eventually always hitting points for interval maps}. \emph{Nonlinearity} 33(2) (2020), 892--914.

\bibitem{KL21}
A. Korepanov and J. Lepp\"anen. \href{https://doi.org/10.1007/s00220-021-04071-5}{Loss of memory and moment bounds for nonstationary intermittent dynamical systems}. \emph{Comm. Math. Phys.} 385(2) (2021), 905--935.

\bibitem{KL25}
A. Korepanov and J. Lepp\"anen. \href{https://doi.org/10.1016/j.physd.2025.134939}{Improved polynomial rates of memory loss for nonstationary intermittent dynamical systems}. \emph{Phys. D} 483 (2025), Paper No. 134939.

\bibitem{LSV99}
C. Liverani, B. Saussol and S. Vaienti. \href{https://doi.org/10.1017/S0143385799133856}{A probabilistic approach to intermittency}. \emph{Ergod. Th. \& Dynam. Sys.} 19(3) (1999), 671--685.

\bibitem{NTV}
M. Nicol, A. T\"or\"ok and S. Vaienti. \href{https://doi.org/10.1017/etds.2016.69}{Central limit theorems for sequential and random intermittent dynamical systems}. \emph{Ergod. Th. \& Dynam. Sys.} 38(3) (2018), 1127--1153.

\bibitem{PY22}
J. Pitman and Z. You. The range of a self-similar additive gamma process is a scale invariant Poisson point process. Preprint, 2022, \href{https://arxiv.org/abs/2111.09409v2}{arXiv:2111.09409v2}.

\bibitem{RS26}
A. Rodriguez Sponheimer. Limit laws for Poincar\'e recurrence and the shrinking target problem. Preprint, 2026, \href{https://arxiv.org/abs/2510.12596v3}{arXiv:2510.12596v3}.

\bibitem{Sera20}
T. Sera. \href{https://doi.org/10.1088/1361-6544/ab5ceb}{Functional limit theorem for occupation time processes of intermittent maps}. \emph{Nonlinearity} 33(3) (2020), 1183--1217.

\bibitem{Sera26}
T. Sera. \href{https://doi.org/10.1017/etds.2025.10246}{Large deviations for occupation and waiting times of infinite ergodic transformations}. \emph{Ergod. Th. \& Dynam. Sys.} 46(3) (2026), 805--844.

\bibitem{BCNote}
J. Xie and Q. Wang. \href{https://doi.org/10.1360/SSM-2024-0256}{A note on the first and the second Borel--Cantelli lemmas}. \emph{Sci. Sin. Math.} 55(7) (2025), 1261--1266 (in Chinese).

\bibitem{Zwe19}
R. Zweim\"uller. \href{https://doi.org/10.1016/j.spa.2018.05.011}{Hitting-time limits for some exceptional rare events of ergodic maps}. \emph{Stochastic Process. Appl.} 129(5) (2019), 1556--1567.

\bibitem{Zwe22}
R. Zweim\"uller. \href{https://doi.org/10.5802/ahl.154}{Hitting times and positions in rare events}. \emph{Ann. H. Lebesgue} 5 (2022), 1361--1415.

\end{thebibliography}
}
\end{document}